\documentclass[10.5pt]{article}
\usepackage[leqno]{amsmath}
\usepackage{amsfonts}
\usepackage{graphicx}
\usepackage{amsmath}
\usepackage{amssymb}
\usepackage{latexsym}
\usepackage{amsmath,amsfonts,amssymb,amsthm,euscript,makeidx,color,mathrsfs}

\def\sqr#1#2{{\vcenter{\vbox{\hrule height.#2pt
              \hbox{\vrule width.#2pt height#1pt \kern#1pt \vrule width.#2pt}
              \hrule height.#2pt}}}}
\def\signed #1{{\unskip\nobreak\hfil\penalty50
              \hskip2em\hbox{}\nobreak\hfil#1
              \parfillskip=0pt \finalhyphendemerits=0 \par}}
\def\endpf{\signed {$\sqr69$}}

\def\3n{\negthinspace \negthinspace \negthinspace }
\def\2n{\negthinspace \negthinspace }
\def\1n{\negthinspace }

\def\dbE{\mathbb{E}}
\def\dbF{\mathbb{F}}

\def\dbH{\mathbb{H}}

\def\dbP{\mathbb{P}}

\def\dbR{\mathbb{R}}
\def\dbS{\mathbb{S}}

\def\sN{\mathscr{N}}

\def\sR{\mathscr{R}}
\def\sS{\mathscr{S}}

\def\sU{\mathscr{U}}

\def\sX{\mathscr{X}}

\def\={\buildrel \triangle \over =}

\def\ds{\displaystyle}

\def\ns{\noalign{\ss}}
\def\a{\alpha}

\def\d{\delta}
\def\e{\varepsilon}
\def\z{\zeta}

\def\n{\nu}
\def\si{\sigma}

\def\f{\varphi}
\def\th{\theta}
\def\o{\omega}

\def\i{\infty}
\def\R{{\bf R}}
\def\G{\Gamma}
\def\D{\Delta}
\def\Th{\Theta}
\def\L{\Lambda}
\def\Si{\Sigma}

\def\OO{\Omega}

\def\cA{{\cal A}}

\def\cC{{\cal C}}

\def\cF{{\cal F}}

\def\cM{{\cal M}}

\def\cQ{{\cal Q}}

\def\cS{{\cal S}}

\def\BA{{\bf A}}
\def\BB{{\bf B}}
\def\BC{{\bf C}}
\def\BD{{\bf D}}

\def\BI{{\bf I}}
\def\BJ{{\bf J}}

\def\BP{{\bf P}}
\def\BQ{{\bf Q}}
\def\BR{{\bf R}}
\def\BS{{\bf S}}

\def\BY{{\bf Y}}
\def\BZ{{\bf Z}}
\def\BTh{{\bf\Theta}}
\def\BPhi{{\bf\Phi}}

\def\BSi{{\bf\Sigma}}

\def\ss{\smallskip}
\def\ms{\medskip}

\def\q{\quad}
\def\qq{\qquad}
\def\hb{\hbox}

\def\limsup{\mathop{\overline{\rm lim}}}
\def\liminf{\mathop{\underline{\rm lim}}}

\def\Ra{\mathop{\Rightarrow}}

\def\lan{\mathop{\langle}}
\def\ran{\mathop{\rangle}}

\def\h{\widehat}
\def\wt{\widetilde}

\def\cd{\cdot}

\def\ae{\hbox{\rm a.e.{}}}
\def\as{\hbox{\rm a.s.{}}}

\def\deq{\mathop{\buildrel\D\over=}}

\def\({\Big (}
\def\){\Big )}
\def\[{\Big[}
\def\]{\Big]}

\def\lan{\langle}
\def\ran{\rangle}
\def\llan{\left\langle}
\def\rran{\right\rangle}
\def\blan{\big\langle}
\def\bran{\big\rangle}

\def\rf{\eqref}
\def\bde{\begin{definition}\label}
\def\ede{\end{definition}}
\def\be{\begin{equation}}
\def\bel{\begin{equation}\label}
\def\ee{\end{equation}}
\def\bt{\begin{theorem}\label}
\def\et{\end{theorem}}
\def\bc{\begin{corollary}\label}
\def\ec{\end{corollary}}
\def\bl{\begin{lemma}\label}
\def\el{\end{lemma}}
\def\bp{\begin{proposition}\label}
\def\ep{\end{proposition}}
\def\bas{\begin{assumption}}
\def\eas{\end{assumption}}
\def\br{\begin{remark}\label}
\def\er{\end{remark}}
\def\bex{\begin{example}\label}
\def\ex{\end{example}}

\newtheorem{theorem}{Theorem}[section]
\newtheorem{corollary}[theorem]{Corollary}

\newtheorem{lemma}[theorem]{Lemma}
\newtheorem{proposition}[theorem]{Proposition}

\theoremstyle{definition}
\newtheorem{definition}[theorem]{Definition}
\newtheorem{example}[theorem]{Example}
\newtheorem{remark}[theorem]{Remark}

\def\ba{\begin{array}}
\def\ea{\end{array}}
\def\ben{\begin{enumerate}}
\def\een{\end{enumerate}}
\def\ed{\end{document}}
\def\les{\leqslant}
\def\ges{\geqslant}

\def\square#1{\vbox{\hrule\hbox{\vrule height#1%
     \kern#1\vrule}\hrule}}
\def\rectangle#1#2{\vbox{\hrule\hbox{\vrule height#1%
     \kern#2\vrule}\hrule}}

\font\tenbb=msbm10 \font\sevenbb=msbm7 \font\fivebb=msbm5

\newfam\bbfam
\scriptscriptfont\bbfam=\fivebb \textfont\bbfam=\tenbb
\scriptfont\bbfam=\sevenbb

\newtheorem{assumption}{Assumption}[section]

\makeatletter
   
   \@addtoreset{equation}{section}
\makeatother

\begin{document}
\title{\bf Mean-Field Linear-Quadratic Stochastic\\ Differential Games in an Infinite Horizon}

\author{Xun Li\footnote{Department of Applied Mathematics, Hong Kong Polytechnic University, Hong Kong, China. This author was supported in part by Research Grants Council of Hong Kong under Grant 15215319.}, \quad
Jingtao Shi\footnote{School of Mathematics, Shandong University, Jinan 250100, China. This author was financially supported by National Key R$\&$D Program of China under Grant 2018YFB1305400, National Natural Science Funds of China under Grant 11971266, 11831010, 11571205, China Scholarship Council and Shandong Provincial Natural Science Foundations under Grant ZR2020ZD24 and ZR2019ZD42.}, \quad
Jiongmin Yong\footnote{Department of Mathematics, University of Central Florida, Orlando, FL 32816, USA. This author was supported in part by NSF under Grant DMS-1812921.}}

\maketitle

\begin{abstract}
This paper is concerned with two-person mean-field linear-quadratic non-zero sum stochastic differential games in an infinite horizon. Both open-loop and closed-loop Nash equilibria are introduced. Existence of an open-loop Nash equilibrium is characterized by the solvability of a system of mean-field forward-backward stochastic differential equations in an infinite horizon and the convexity of the cost functionals, and the closed-loop representation of an open-loop Nash equilibrium is given through the solution to a system of two coupled non-symmetric algebraic Riccati equations. The existence of a closed-loop Nash equilibrium is  characterized by the solvability of a system of two coupled symmetric algebraic Riccati equations. Two-person mean-field linear-quadratic zero-sum stochastic differential games in an infinite time horizon are also considered. Both the existence of open-loop and closed-loop saddle points are characterized by the solvability of a system of two coupled generalized algebraic Riccati equations with static stabilizing solutions. Mean-field linear-quadratic stochastic optimal control problems in an infinite horizon are discussed as well, for which it is proved that the open-loop solvability and closed-loop solvability are equivalent.

\end{abstract}


\ms


\section{Introduction}
Let $(\OO,\cF,\dbP,\dbF)$ be a complete filtered probability space, on which a one-dimensional standard Brownian motion $W(\cd)$ is defined with $\dbF\equiv\{\cF_t\}_{t\ges0}$ being its natural filtration augmented by all the $\dbP$-null sets in $\cF$, and $\dbE[\,\cd\,]$ denotes the expectation with respect to $\dbP$. Throughout this paper, we let $\dbR^{n\times m}$ and $\dbS^n$ be the set of all $(n\times m)$ (real) matrices and $(n\times n)$ symmetric (real) matrices. We denote $\dbR^n=\dbR^{n\times1}$. For a Euclidean space $\dbH$, say, $\dbH=\dbR^n,\dbR^{n\times m}$, let $C([0,\i);\dbH)$ denote the space of $\dbH$-valued continuous functions $\f:[0,\i)\to\dbH$, $L^2(\dbH)$ denote the space of $\dbH$-valued functions $\f:[0,\i)\rightarrow\dbH$ with $\int_0^\i|\f(t)|^2dt<\i$, and $L_\dbF^2(\dbH)$ denote the space of $\dbF$-progressively measurable processes $\f:[0,\i)\times\Omega\rightarrow\dbH$ with $\dbE\int_0^\i|\f(t)|^2dt<\i$.

Consider the following controlled linear {\it mean-field stochastic differential equation} (MF-SDE, for short) on an infinite horizon $[0,\i)$:
\bel{state-LQ-game}\left\{\2n\ba{ll}
\ds dX(t)=\big\{AX(t)+\bar A\dbE[X(t)]+B_1u_1(t)+\bar B_1\dbE[u_1(t)]+B_2u_2(t)+\bar B_2\dbE[u_2(t)]+b(t)\big\}dt\\
\ns\ds\qq\qq+\big\{CX(t)\1n+\1n\bar C\dbE[X(t)]\1n+\1n D_1u_1(t)\1n+\1n\bar D_1\dbE[u_1(t)]\1n+\1n D_2u_2(t)\1n+\1n\bar D_2\dbE[u_2(t)]\1n+\1n\si(t)\big\}dW(t),~ t\ges0,\\
\ns\ds X(0)=x.\ea\right.\ee
In the above, $X(\cd)$ is the state process taking values in $\dbR^n$ with $x$ being the initial state. For $i=1,2$, $u_i(\cd)$ is the control process of Player $i$, taking values in $\dbR^{m_i}$, respectively. The coefficients $A,\bar A,C,\bar C\in\dbR^{n\times n}$, $B_1,\bar B_1,D_1,\bar D_1\in\dbR^{n\times m_1}$, $ B_2,\bar B_2,D_2,\bar D_2\in\dbR^{n\times m_2}$ are given constant matrices, and the non-homogenous terms $b(\cd),\si(\cd)\in L_\dbF^2(\dbR^n)$. We introduce the following spaces:
$$\left\{\ba{ll}
\ds\sX[0,T]=\Big\{X:[0,\i)\times\Omega\to\dbR^n\bigm| X(\cd)\mbox{ is }\dbF\mbox{-adapted}, t\mapsto X(t,\o)\hb{ is continuous, and}\\
\ns\ds\qq\qq\qq\qq\qq\qq\qq\qq\dbE\Big[\sup_{t\in[0,T]}|X(t)|^2\Big]<\i\Big\},\q \hb{for }T>0,\\
\ns\ds\sX_{loc}[0,\i)=\bigcap_{T>0}\sX[0,T],\qq
\sX[0,\i)=\Big\{X(\cd)\in\sX_{loc}[0,\i)\bigm|\dbE\int_0^\i|X(t)|^2dt<\i\Big\}.\ea\right.$$
By a standard argument using contraction mapping theorem, one can show that for any initial state $x\in\dbR^n$ and control pair $(u_1(\cd),u_2(\cd))\in L_\dbF^2(\dbR^{m_1})\times L_\dbF^2(\dbR^{m_2})$, state equation \rf{state-LQ-game} admits a unique strong solution $X(\cd)\equiv X(\cd;x,u_1(\cd),u_2(\cd))\in\sX_{loc}[0,\i)$. Next, for $i=1,2$, we introduce the following cost functionals:
\bel{cost-LQ-game}J_i(x;u_1(\cd),u_2(\cd))=\dbE\int_0^\i g_i\big(t,X(t),u_1(t),u_2(t),\dbE[X(t)],\dbE[u_1(t)],\dbE[u_2(t)]\big)dt,\ee
with
\bel{g_i}\ba{ll}
\ns\ds g_i(t,x,u_1,u_2,\bar x,\bar u_1,\bar u_2)=\llan\2n
\begin{pmatrix}Q_i&S_{i1}^\top&S_{i2}^\top\\ S_{i1}&R_{i11}&R_{i12}\\ S_{i2}&R_{i21}&R_{i22}\end{pmatrix}
\begin{pmatrix}x\\ u_1\\ u_2\end{pmatrix},\begin{pmatrix}x\\ u_1\\ u_2\end{pmatrix}\2n\rran
+2\llan\2n\begin{pmatrix}q_i(t)\\ \rho_{i1}(t)\\ \rho_{i2}(t)\end{pmatrix},
\begin{pmatrix}x\\ u_1\\ u_2\end{pmatrix}\2n\rran\\
\ns\ns\ds\qq\qq\qq\qq\qq\qq+\llan\2n
\begin{pmatrix}\bar Q_i&\bar S_{i1}^\top&\bar S_{i2}^\top\\ \bar S_{i1}&\bar R_{i11}&\bar R_{i12}\\ \bar S_{i2}&\bar R_{i21}&\bar R_{i22}\end{pmatrix}
\begin{pmatrix}\bar x\\ \bar u_1\\ \bar u_2\end{pmatrix},\begin{pmatrix}\bar x\\ \bar u_1\\ \bar u_2\end{pmatrix}
\2n\rran,\ea\ee
where
\begin{equation*}\left\{\2n\ba{lll}
\ds Q_i,~\bar Q_i\in\dbS^n, & S_{i1},~\bar S_{i1}\in\dbR^{m_1\times n}, & S_{i2},~\bar S_{i2}\in\dbR^{m_2\times n},\\
\ns\ds R_{i11},~\bar R_{i11}\in\dbS^{m_1}, & R_{i22},~\bar R_{i22}\in\dbS^{m_2}, & R_{i12}=R_{i21}^\top,~\bar R_{i12}=\bar R_{i21}^\top\in\dbR^{m_1\times m_2},\\
\ns\ds q_i(\cd)\in L_\dbF^2(\dbR^n), & \rho_{i1}(\cd)\in L_\dbF^2(\dbR^{m_1}), & \rho_{i2}(\cd)\in L_\dbF^2(\dbR^{m_2}).\ea\right.
\end{equation*}
Note that for $(x,u_1(\cd),u_2(\cd))\in\dbR^n\times L_\dbF^2(\dbR^{m_1})\times L_\dbF^2(\dbR^{m_2})$, the solution $X(\cd)\equiv X(\cd\,;x,u_1(\cd),u_2(\cd))$ to \rf{state-LQ-game} might just be in $\sX_{loc}[0,\i)$ in general. Therefore, in order the cost functionals $J_i(x;u_1(\cd),u_2(\cd)),i=1,2$ to be defined, the control pair $(u_1(\cd),u_2(\cd))$ has to be restricted in the following set of {\it admissible control pairs}:
\bel{Uad-LQ-game}
\sU_{ad}(x)\1n=\Big\{(u_1(\cd),u_2(\cd))\in L_\dbF^2(\dbR^{m_1})\times L_\dbF^2(\dbR^{m_2})\bigm|X(\cd\,;x,u_1(\cd),u_2(\cd))\in\sX[0,\i)\Big\},\q x\in\dbR^n.
\ee
Note that $\sU_{ad}(x)$ depends on the initial state $x$. For any $(u_1(\cd),u_2(\cd))\in\sU_{ad}(x)$, the corresponding $X(\cd)\equiv X(\cd\,;x,u_1(\cd),u_2(\cd))$ is called an {\it admissible state process} for the initial state $x$. Then we can loosely formulate the following problem.

\ms

{\bf Problem (MF-SDG).}\q For any initial state $x\in\dbR^n$, Player $i$ ($i=1,2$) wants to find a control $u_i^*(\cd)$ so that $(u_1^*(\cd),u_2^*(\cd))\in\sU_{ad}(x)$ such that the cost functionals $u_1(\cd)\mapsto J_1(x;u_1(\cd),u^*_2(\cd))$ and $u_2(\cd)\mapsto J_2(x;u_1^*(\cd),u_2(\cd))$ are minimized, for all $(u_1(\cd),u_2^*(\cd)),
(u_1^*(\cd),u_2(\cd))\in\sU_{ad}(x)$, subject to \rf{state-LQ-game}.

\ms

We refer to the above problem as a {\it mean-field linear-quadratic {\rm(LQ, for short)} two-person (non-zero sum) stochastic differential game in an infinite horizon}. In the special case where $b(\cd),\si(\cd),q_i(\cd),\rho_{ij}(\cd)$ are all zero, we denote the corresponding problem and cost functionals by Problem (MF-SDG)$^0$ and $J_i^0(x;u_1(\cd),u_2(\cd))$, respectively. On the other hand, when
\bel{zero-sum-LQ-game}
J_1(x;u_1(\cd),u_2(\cd))+J_2(x;u_1(\cd),u_2(\cd))=0,\q \forall x\in\dbR^{n},\ \forall(u_1(\cd),u_2(\cd))\in\sU_{ad}(x),
\ee
the corresponding Problem (MF-SDG) is called a {\it mean-field LQ two-person zero-sum stochastic differential game in an infinite horizon} and is denoted by Problem (MF-SDG)$_0$. To guarantee \rf{zero-sum-LQ-game}, one usually lets
\bel{zero-zero}\left\{\2n\ba{lll}
\ds Q_1+Q_2=0, & \bar Q_1+\bar Q_2=0, & q_1(\cd)+q_2(\cd)=0,\\
\ns\ds S_{1j}+S_{2j}=0, & \bar S_{1j}+\bar S_{2j}=0, & \rho_{1j}(\cd)+\rho_{2j}(\cd)=0,\qq j=1,2,\\
\ns\ds R_{1jk}+R_{2jk}=0, & \bar R_{1jk}+\bar R_{2jk}=0, & \q j,k=1,2.\ea\right.\ee

One may feel that cost functional \rf{cost-LQ-game} could be a little more general by including terms like $\lan\bar q_i(s),\dbE[X(s)]\ran$. However, it is not hard to see that (as long as the integrals exist)
$$\dbE\int_0^\infty\blan\bar q_i(s),\dbE[X(s)]\bran ds=\dbE\int_0^\infty\blan\dbE[\bar q_i(s)],X(s)\bran ds,$$
which can be absorbed by replacing $q_i(\cd)$ by $q_i(\cd)+\dbE[\bar q_i(\cd)]$ in the cost functionals. Likewise, terms like $\blan\bar\rho_{ij}(s),\dbE[u_j(s)]\bran$ are not necessarily included.

\ms

We will introduce proper stabilizability conditions for the system so that $\sU_{ad}(x)\ne\varnothing$ for all $x\in\dbR^n$. It is not hard to see that without stabilizabity conditions, one might only has $X(\cd\,;x,u_1(\cd),u_2(\cd))\in\sX_{loc}[0,\infty)$ so that $J_i(x;u_1(\cd),u_2(\cd))$ might not be well-defined. Then, instead of considering the cost functionals of form \rf{cost-LQ-game}, one might naturally consider the following ergodic type cost functionals:
\bel{ergodic}\wt J_i(x;u_1(\cd),u_2(\cd))=\limsup_{T\to\i}{1\over T}\int_0^Tg_i\big(t,X(t),u_1(t),u_2(t),\dbE[X(t)],\dbE[u_1(t)],\dbE[u_2(t)]\big)dt.\ee
The above type cost functionals are normally used for the case that the running cost rate is bounded so that the right-hand side of the above is always finite. However, if no stabilizability conditions are assumed, the state process $X(\cd)$ could be of exponential growth. In that case, the above type cost functionals are still not useful. Therefore, instead of \rf{ergodic}, we prefer to study our problems under certain stabilizability condition with cost functionals of for \rf{cost-LQ-game} restricted on $\sU_{ad}(x)$.

\ms

The theory of MF-SDEs can be traced back to the work of Kac \cite{Kac1956} in the middle of 1950s, where a stochastic toy model for the Vlasov type kinetic equation of plasma was presented. Its rigorous study was initiated by McKean \cite{McKean1966} in 1966, which is now known as McKean-Vlasov stochastic differential equations. Since then, many researchers have made contributions to the related topics and their applications; see, for example, Scheutzow \cite{Scheutzow1988}, Chan \cite{Chan1994}, Huang--Malhame--Caines \cite{Huang-Malhame-Caines2006}, Buckdahn--Li--Peng \cite{Buckdahn-Li-Peng2009}, Carmona and Delarue \cite{Carmona-Delarue2013}, Bensoussan--Yam--Zhang \cite{Bensoussan-Yam-Zhang2015}, Buckdahn--Li--Peng--Rainer \cite{Buckdahn-Li-Peng-Rainer2017}, etc.

\ms

Optimal control and differential game problems of MF-SDEs have drawn enormous researchers' attention recently. See Ahmed--Ding \cite{Ahmed-Ding2001}, Lasry--Lions \cite{Lasry-Lions2007}, Andersson--Djehiche\cite{Andersson-Djehiche2011}, Buckdahn--Djehiche-Li \cite{Buckdahn-Djehiche-Li2011}, Li \cite{Li2012}, Meyer-Brandis--Oksendal--Zhou \cite{MeyerBrandis-Oksendal-Zhou2012}, Hosking \cite{Hosking2012}, Bensoussan--Sung--Yam \cite{Bensoussan-Sung-Yam2013}, Djehiche--Tembine-Tempone \cite{Djehiche-Tembine-Tempone2015}, Bensoussan--Sung--Yam--Yung \cite{Bensoussan-Sung-Yam-Yung2016}, Djehiche--Huang \cite{Djehiche-Huang2016}, Huang--Li--Wang \cite{Huang-Li-Wang2016}, Yong \cite{Yong2017}, Buckdahn--Li--Ma \cite{Buckdahn-Li-Ma2016,Buckdahn-Li-Ma2017}, Pham--Wei \cite{Pham-Wei2017,Pham-Wei2018}, Li--Sun--Xiong \cite{Li-Sun-Xiong2018}, Miller--Pham \cite{Miller-Pham2019}, Moon \cite{Moon2020}, and the references therein. Next, let us mention a few recent pieces of literature related to our present paper. In Yong \cite{Yong2013}, an LQ optimal control problem for MF-SDEs in a finite horizon was introduced and investigated. The optimality system of a linear {\it mean-field forward-backward stochastic differential equation} (MF-FBSDE, for short) is derived, and two Riccati differential equations are obtained to present the feedback representation of the optimal control. Huang--Li--Yong \cite{Huang-Li-Yong2015} generalized the results in \cite{Yong2013} to the infinite horizon case, and the feedback representation of the optimal control is derived via two {\it algebraic Riccati equations} (AREs, for short). Note that in \cite{Huang-Li-Yong2015}, some notions of stabilizability for controlled MF-SDEs are introduced, which are interestingly different from the classical ones, due to the presence of the terms $\dbE[X(\cd)]$ and $\dbE[u(\cd)]$. Sun \cite{Sun2017} continued to investigate the LQ optimal control problem for MF-SDEs in the finite horizon with additional nonhomogeneous terms and concluded that the uniform convexity of the cost functional is sufficient for the open-loop solvability of the LQ optimal control problems for MF-SDEs. Moreover, the uniform convexity of the cost functional is equivalent to the solvability of two coupled differential Riccati equations, and the unique open-loop optimal control admits a state feedback representation in the case that the cost functional is uniformly convex. Li--Sun--Yong \cite{Li-Sun-Yong2016} studied the closed-loop solvability of the corresponding problem, which is characterized by the existence of a regular solution to the coupled two generalized Riccati equations, together with some constraints on the adapted solution to a linear {\it backward stochastic differential equation} (BSDE, for short) and a linear terminal value problem of an {\it ordinary differential equation} (ODE, for short). Li--Li--Yu \cite{Li-Li-Yu2020} analyzed the indefinite mean-field type LQ stochastic optimal control problems, where they introduced a relaxed compensator to characterize the open-loop solvability of the problem. Tian--Yu--Zhang \cite{Tian-Yu-Zhang2020} considered an LQ zero-sum stochastic differential game with mean-field type, proposed the notions of explicit and implicit strategy laws, and established the closed-loop formulation for saddle points in the mixed-strategy-law. Very recently, Sun--Wang--Wu \cite{Sun-Wang-Wu2021} studied a two-person zero-sum mean-field LQ stochastic differential game over the finite horizon by a Hilbert space method introduced by Mou--Yong \cite{Mou-Yong2006}. It is shown that the associated two Riccati equations admit unique and strongly regular solutions under the sufficient condition for the existence of an open-loop saddle point when the open-loop saddle point can be represented as linear feedback of the current state. When only the necessary condition for the existence of an open-loop saddle point is satisfied, we can construct an approximate sequence by solving a family of Riccati equations and closed-loop systems. The approximate sequence's convergence turns out to be equivalent to the open-loop solvability of the game, and its limit exactly equals an open-loop saddle point, provided that the game is open-loop solvable.

\ms

Ait Rami--Zhou \cite{AitRami-Zhou2000} and Ait Rami--Zhou--Moore \cite{AitRami-Zhou-Moore2000} considered stochastic LQ problems in an infinite horizon, with indefinite control weighting matrices. They introduced a generalized ARE, involving a matrix pseudo-inverse and two additional algebraic equality/inequality constraints, and proved that the problem's attainability is equivalent to the existence of a static stabilizing solution to the generalized ARE. In particular, the associated AREs can be solved by linear matrix inequality and semidefinite programming techniques. In addition to the statements in the previous paragraph about \cite{Huang-Li-Yong2015}, the authors discussed the solvabilities of AREs, by linear matrix inequalities. Sun--Yong \cite{Sun-Yong2018} first found that both the open-loop and closed-loop solvabilities of the stochastic LQ problems in the infinite horizon are equivalent to the existence of a static stabilizing solution to the associated generalized ARE. We refer the readers to the recent monographs by Bensoussan--Frehse--Yam \cite{Bensoussan-Frehse-Yam2013} and by Sun--Yong \cite{Sun-Yong2020a, Sun-Yong2020b} for more details and references cited therein.

\ms

In this paper, we consider two-person mean-field LQ non-zero sum stochastic differential games in an infinite horizon. Both open-loop and closed-loop Nash equilibria are introduced. The existence of an open-loop Nash equilibrium is characterized by the solvability of a system of MF-FBSDEs in an infinite horizon and the convexity of the cost functionals. The closed-loop representation of an open-loop Nash equilibrium is given through the solution to a system of two coupled non-symmetric AREs. The existence of a closed-loop Nash equilibrium is characterized by the solvability of a system of two coupled symmetric AREs. Two-person mean-field LQ zero-sum stochastic differential games in an infinite horizon are also considered. The existence of open-loop and closed-loop saddle points is characterized by the solvability of a system of two coupled generalized AREs with static stabilizing solutions. As special cases, mean-field LQ stochastic optimal control problems in an infinite horizon are discussed as well, for which it is proved that the open-loop solvability and closed-loop solvability are equivalent. The results obtained in this paper enrich the theory of optimal control and differential games of mean-field type.

\ms

Let us briefly highlight the major novelty of this paper:

\ms

(i) For MF-SDEs with quadratic performance indexes in an infinite-horizon, problems of two-person non-zero sum differential games, two-person zero-sum differential games and optimal control (which is a single player differential game) are treated in a unified framework. Among other results, the most significant one is the discovery of the system of coupled algebraic Riccati equations (\ref{AE-LQ-game-closed}) which is used to characterize the closed-loop Nash equilibrium. From this point of view, the current paper can be regarded as a complementary or a continuation of Sun--Yong \cite{Sun-Yong2019}.

\ms

(ii) For MF-SDE LQ optimal control problems in an infinite horizon, we have established the equivalence among the solvability of a system of coupled algebraic Riccati equations, open-loop solvability, and closed-loop solvability. This covers the relevant results found in Sun--Yong \cite{Sun-Yong2018} where mean-field terms were absent.

\ms

(iii) For two-person zero-sum differential games of MF-SDEs with quadratic performance index in $[0,\i)$, we have proved that if an open-loop saddle point admits a closed-loop representation, and the closed-loop saddle point exists, then the open-loop saddle point must be the outcome of the closed-loop saddle point. It is also shown that such a property fails for non-zero sum differential games. From this angle, the current paper is an extension of Sun--Yong--Zhang \cite{Sun-Yong-Zhang2016} where the mean-field terms did not appear.

\ms

The rest of the paper is organized as follows. In Section 2, we present some preliminary results about mean-field LQ stochastic optimal control problems in an infinite horizon. Section 3 aims to give results on mean-field LQ non-zero sum stochastic differential games, including open-loop Nash equilibria and their closed-loop representation, and closed-loop Nash equilibria with algebraic Riccati equations. In Section 4, the open-loop and closed-loop saddle points for mean-field LQ zero-sum stochastic differential games are investigated. Some examples are presented in Section 5 illustrating the results developed in the earlier sections. In Section 6, the detailed proof of the main result in Section 2 is given. Finally, some concluding remarks are collected in Section 7.

\section{Preliminaries}
Throughout this paper, besides the notation introduced in the previous section, we let $I$ be the identity matrix or operator with a suitable size. We will use $\lan\cd\,,\cd\ran$ for inner products in possibly different Hilbert spaces, and denote by $|\cd|$ the norm induced by $\lan\cd\,,\cd\ran$. Let $M^\top$ and $\sR(M)$ be the transpose and range of a matrix $M$, respectively. For $M,N\in\dbS^n$, we write $M\ges N$ (respectively, $M>N$) for $M-N$ being positive semi-definite (respectively, positive definite). Let $M^\dag$ denote the pseudo-inverse of a matrix $M\in\dbR^{m\times n}$, which is equal to the inverse $M^{-1}$ of $M\in\dbR^{n\times n}$ if it exists. See Penrose \cite{Penrose1955} or Anderson--Moore \cite{Anderson-Moore1989} for some basic properties of the pseudo-inverse. We define the inner product in $L^2_\dbF(\dbH)$ by
$\lan\f,\phi\ran=\dbE\int_0^\i\lan\f(t),\phi(t)\ran dt$ so that $L^2_\dbF(\dbH)$ is a Hilbert space.

\ms

We now consider the following controlled linear MF-SDE over $[0,\infty)$:
\bel{state-LQ-control}\left\{\2n\ba{ll}
\ds dX(t)=\big\{AX(t)+\bar A\dbE[X(t)]+Bu(t)+\bar B\dbE[u(t)]+b(t)\big\}dt\\
\ns\ds\qq\qq\q+\big\{CX(t)+\bar C\dbE[X(t)]+Du(t)+\bar D\dbE[u(t)]+\si(t)\big\}dW(t),
\q t\ges0,\\
\ns\ds X(0)=x,\ea\right.\ee
with quadratic cost functional
\bel{cost-LQ-control}\ba{ll}
\ds J(x;u(\cd))=\dbE\int_0^\i\bigg[\bigg\langle\2n
\begin{pmatrix}Q&S^\top\\ S&R\end{pmatrix}
\begin{pmatrix}X(t)\\ u(t)\end{pmatrix},\begin{pmatrix}X(t)\\ u(t)\end{pmatrix}\2n\bigg\rangle
+2\bigg\langle\2n\begin{pmatrix}q(t)\\ \rho(t)\end{pmatrix},\begin{pmatrix}X(t)\\ u(t)\end{pmatrix}\2n\bigg\rangle\\
\ns\ds\qq\qq\qq\qq+\bigg\langle\2n
\begin{pmatrix}\bar Q&\bar S^\top\\ \bar S&\bar R\\\end{pmatrix}
\begin{pmatrix}\dbE[X(t)]\\ \dbE[u(t)]\end{pmatrix},\begin{pmatrix}\dbE[X(t)]\\ \dbE[u(t)]\end{pmatrix}
\2n\bigg\rangle\bigg]dt,\ea\ee
where $A,\bar A,C,\bar C\in\dbR^{n\times n}$, $B,\bar B,D,\bar D\in\dbR^{n\times m}$, $Q,\bar Q\in\dbS^n$, $S,\bar S\in\dbR^{m\times n}$, $R,\bar R\in\dbS^m$ are given
constant matrices, and $b(\cd),\si(\cd),q(\cd)\in L_\dbF^2(\dbR^n)$, $\rho(\cd)\in L_\dbF^2(\dbR^m)$ are stochastic processes. For any initial state $x\in\dbR^n$ and control $u(\cd)\in L_\dbF^2(\dbR^m)$, equation \rf{state-LQ-control} admits a unique strong solution $X(\cd)\equiv X(\cd\,;x,u(\cd))\in\sX_{loc}[0,\i)$. We define the {\it admissible control set} as
\bel{Uad-LQ-control}
\sU_{ad}(x)\1n=\big\{u(\cd)\in L_\dbF^2(\dbR^m)\bigm|X(\cdot)\equiv X(\cd;x,u(\cd))\in\sX[0,\i)\big\}.
\ee
In general, $\sU_{ad}(x)$ depends on $x\in\dbR^n$. Let us pose the following optimal control problem.

\ms

{\bf Problem (MF-SLQ).}\q For any initial state $x\in\dbR^n$, find a control $u^*(\cd)\in\sU_{ad}(x)$ such that the cost functional $J(x;u(\cd))$  of
\rf{cost-LQ-control} is minimized, subject to \rf{state-LQ-control}. That is to say,
\bel{value-LQ-control}
J(x;u^*(\cd))=\inf\limits_{u(\cd)\in\sU_{ad}(x)}J(x;u(\cd))\equiv V(x).
\ee

Any $u^*(\cd)\in\sU_{ad}(x)$ satisfying \rf{value-LQ-control} is called an {\it open-loop optimal control} of Problem (MF-SLQ), and the corresponding $X^*(\cdot)\equiv X(\cd\,;x,u^*(\cd))$ is called an {\it open-loop optimal state process}. The function $V(\cd)$ is called the {\it value function} of Problem (MF-SLQ). In the special case where $b(\cd),\si(\cd),q(\cd),\rho(\cd)$ are all zero, we denote the corresponding problem by Problem (MF-SLQ)$^0$, the cost functional by $J^0(x;u(\cd))$ and the value function by $V^0(x)$, respectively.

\ms

In order Problem (MF-SLQ) to be meaningful, we need to find conditions under which  $\sU_{ad}(x)$ is non-empty and admits an accessible characterization. For this target, let us first look at the following uncontrolled non-homogeneous linear system on $[0,\infty)$:
\bel{MF-SDE}\left\{\2n\ba{ll}
\ds dX(t)=\big\{AX(t)+\bar A\dbE[X(t)]+b(t)\big\}dt+\big\{CX(t)+\bar C\dbE[X(t)]+\si(t)\big\}dW(t),\q t\ges0,\\
\ns\ds X(0)=x.\ea\right.\ee
When $b(\cd)=\si(\cd)=0$, the system is said to be {\it homogeneous} and denoted by $[A,\bar A,C,\bar C]$. For simplicity, we also denote $[A,C]=[A,0,C,0]$ (the linear SDE without mean-fields), and $A=[A,0]=[A,0,0,0]$ (the linear {\it ordinary differential equation}, ODE, for short). The following notions can be found in \cite{Huang-Li-Yong2015}.
\bde{D-2.1} \rm
(i)\ System $[A,\bar A,C,\bar C]$ is said to be $L^2$-{\it globally integrable}, if for any $x\in\dbR^n$, the solution $X(\cd)\equiv X(\cd\,;x)$ of \rf{MF-SDE} with $b(\cd)=\si(\cd)=0$ is in $\sX[0,\i)$.

\ms

(ii)\ System $[A,\bar A,C,\bar C]$ is said to be $L^2$-{\it asymptotically stable}, if for any $x\in\dbR^n$, the solution $X(\cd)\equiv X(\cd\,;x)\in\sX_{loc}[0,\i)$ of \rf{MF-SDE} with $b(\cd)=\si(\cd)=0$ satisfies $\ds\lim_{t\to\i}\dbE|X(t)|^2=0$.
\ede

According to \cite{Huang-Li-Yong2015}, and via a similar argument proving Theorem 3.3 of \cite{Sun-Yong-Zhang2016}, we have the following result.
\bp{P-2.2} \sl For any $x\in\dbR^n$ and $b(\cd),\si(\cd)\in L_\dbF^2(\dbR^n)$, linear MF-SDE \rf{MF-SDE} admits a unique solution $X(\cd)\in\sX_{loc}[0,\i)$. Further, if $[A,\bar A,C,\bar C]$ is $L^2$-asymptotically stable and system $[A,C]$ is $L^2$-globally integrable, then $X(\cd)\in\sX[0,\i)$, with
$$\dbE\int_0^\i |X(t)|^2dt\les K\bigg[|x|^2+\dbE\int_0^\i\big(|b(t)|^2+|\si(t)|^2\big)dt\bigg],$$
for some constant $K>0$. On the other hand, for any $\f(\cd)\in L_\dbF^2(\dbR^n)$, the following linear MF-BSDE:
\bel{MF-BSDE}
-dY(t)=\big\{A^\top Y(t)+\bar A^\top\dbE[Y(t)]+C^\top Z(t)+\bar C^\top\dbE[Z(t)]+\f(t)\big\}dt-Z(t)dW(t),\q t\ges0,
\ee
admits a unique adapted solution $(Y(\cd),Z(\cd))\in\sX[0,\i)\times L^2_\dbF(\dbR^n)$.
\ep

For general theory of MF-BSDEs and MF-FBSDEs in a finite horizon, see \cite{Buckdahn-Li-Peng2009,Carmona-Delarue2013,Carmona-Delarue2015}.

\ms

Now we return to \rf{state-LQ-control}. Similar to the above, when $b(\cd)=\si(\cd)=0$, the system is said to be {\it homogeneous} and denote it by $[A,\bar A,C,\bar C;B,\bar B,D,\bar D]$.

\bde{feedback} \rm (i) For any $\BTh\equiv(\Th,\bar\Th)\in\dbR^{m\times2n}$ and $v(\cd)\in L^2_\dbF(\dbR^m)$,
\bel{closed-loop control}u(\cd)=u^{\BTh,v}(\cd)\equiv\Th\big\{X(\cd)-\dbE[X(\cd)]\big\}
+\bar\Th\dbE[X(\cd)]+v(\cd)\equiv\BTh\begin{pmatrix}X(\cd)-\dbE[X(\cd)]\\
\dbE[X(\cd)]\end{pmatrix}+v(\cd)\ee
is called a {\it feedback control}. Under such a
control, the state equation \rf{state-LQ-control} becomes
\bel{closed-loop-system}\left\{\ba{ll}
\ds dX(t)=\big\{A_\Th X(t)+\bar A_\BTh\dbE[X(t)]+Bv(t)+\bar B\dbE[v(t)]+b(t)\big\}dt\\
\ns\ds\qq\qq\q+\big\{C_\Th X(t)+\bar C_\BTh\dbE[X(t)]+Dv(t)+\bar D\dbE[v(t)]+\si(t)\big\}dW(t),\qq t\ges0,\\
\ns\ds X(0)=x,\ea\right.\ee
where
\bel{A-Th}
A_\Th=A+B\Th,\q \bar A_\BTh=\bar A+\bar B\bar\Th+B(\bar\Th-\Th),\q C_\Th=C+D\Th,\q \bar C_\BTh=\bar C+\bar D\bar\Th+D(\bar\Th-\Th).
\ee

(ii) System $[A,\bar A,C,\bar C;B,\bar B,D,\bar D]$ is said to be MF-$L^2$-{\it stabilizable}, if there exists a $\BTh\equiv(\Th,\bar\Th)\in\dbR^{m\times2n}$ such that system $[A_\Th,\bar A_\BTh,C_\Th,\bar C_\BTh]$ is $L^2$-asymptotically stable and system $[A_\Th,C_\Th]$ is $L^2$-globally integrable.
In this case, $\BTh\equiv(\Th,\bar\Th)$ is called an MF-$L^2$-{\it stabilizer}
of $[A,\bar A,C,\bar C;B,\bar B,D,\bar D]$. The set of all MF-$L^2$-stabilizers of $[A,\bar A,C,\bar C;B,\bar B,D,\bar D]$ is denoted by $\sS[A,\bar A,C,\bar C;B,\bar B,D,\bar D]$.

\ms

(iii) Any pair $(\BTh,v(\cd))\in\sS[A,\bar A,C,\bar C;B,\bar B,D,\bar D]\times L^2_\dbF(\dbR^m)$ is called a {\it closed-loop strategy} of Problem (MF-SLQ).
The solution $X(\cd)\equiv X^{\BTh,v}(\cd)$ of \rf{closed-loop-system} is called the {\it closed-loop state process} corresponding to $(\BTh,v(\cd))$.
The control $u(\cd)$ defined by \rf{closed-loop control} is called the {\it outcome} of $(\BTh,v(\cd))$, or a {\it closed-loop control} for the initial state $x\in\dbR^n$.

\ede

Note that in the above, $A_\Th$ and $C_\Th$ only depend on $\Th$, and $\bar A_{\BTh}$ and $\bar C_{\BTh}$ depend on $\BTh=(\Th,\bar\Th)$. Also, one sees that the corresponding coefficients $B,\bar B,D,\bar D$ of the control process, as well as the nonhomogeneous terms $b(\cd),\si(\cd)$ are unchanged under \rf{closed-loop control}. See \cite{Huang-Li-Yong2015} for a relevant presentation.

\ms

We introduce the following assumption.

\ms

{\bf (H1)} System $[A,\bar A,C,\bar C;B,\bar B,D,\bar D]$ is MF-$L^2$-stablizable, i.e.,
$\sS[A,\bar A,C,\bar C;B,\bar B,D,\bar D]\ne\varnothing$.

\ms

We have the following result.

\bp{P-2.4} \sl Let {\rm(H1)} hold. Then for any $x\in\dbR^n$, $\sU_{ad}(x)\ne\varnothing$ and $u(\cd)\in\sU_{ad}(x)$ if and only if $u(\cd)=u^{\BTh,v}(\cd)$ given by \rf{closed-loop control} for some $\BTh\equiv(\Th,\bar\Th)\in\sS[A,\bar A,C,\bar C;B,\bar B,D,\bar D]$ and $v(\cd)\in L^2_\dbF(\dbR^m)$, where $X(\cd)\equiv X^{\BTh,v}(\cd)$ is the solution to the closed-loop system \rf{closed-loop-system}.

\ep
\proof
Sufficiency. Let $v(\cd)\in L^2_\dbF(\dbR^m)$ and $X(\cd)$ be the solution to \rf{closed-loop-system}. Since system $[A_\Th,\bar A_\BTh,C_\Th,\bar C_\BTh]$ is $L^2$-asymptotically stable and system $[A_\Th,C_\Th]$ is $L^2$-globally integrable, by Proposition \ref{P-2.2}, the solution $X(\cd)$ to \rf{closed-loop-system} is in $\sX[0,\i)$. Hence, setting $u(\cd)$ by \rf{closed-loop control}, we see that $u(\cd)\in L^2_\dbF(\dbR^m)$. By the uniqueness of the solutions, one has that  $X(\cd)\equiv X^{\BTh,v}(\cd)$ also solves \rf{state-LQ-control}. Therefore, $u(\cd)\in\sU_{ad}(x)$.

\ms

Necessity. Assume that $u(\cd)\in\sU_{ad}(x)$. Let $\BTh\equiv(\Th,\bar\Th)\in\sS[A,\bar A,C,\bar C;B,\bar B,D,\Bar D]$ and the corresponding $X(\cd)\in\sX[0,\i)$ be the solution to \rf{state-LQ-control}. Set
$v(\cd)\triangleq u(\cd)-\Th\big\{X(\cd)-\dbE[X(\cd)]\big\}-\bar\Th\dbE[X(\cd)]\in L^2_\dbF(\dbR^m)$.
By the uniqueness of the solutions again, $X(\cd)$ coincides with the solution to \rf{closed-loop-system}.
Thus, $u(\cd)$ admits a representation of the form \rf{closed-loop control}. The proof is complete. \endpf

\ms

From the above, we can easily show that under (H1), $\sU_{ad}(x)=L^2_\dbF(\dbR^m)$ which is independent of $x$. Hence, hereafter, once (H1) is assumed, we will denote $\sU_{ad}(x)=\sU_{ad}$. Now, we introduce the following definitions concerning Problem (MF-SLQ).

\bde{D-2.5}
(i)\ Problem (MF-SLQ) is said to be {\it finite} at $x\in\dbR^n$ if $V(x)>-\i$, and Problem (MF-SLQ) is said to be finite if it is finite at every $x\in\dbR^n$.

(ii)\ An element $u^*(\cd)\in\sU_{ad}(x)$ is called an {\it open-loop optimal control} of Problem (MF-SLQ) for the initial state $x\in\dbR^n$ if
\bel{open-loop control}
J(x;u^*(\cd))\les J(x;u(\cd)),\qq\forall u(\cd)\in\sU_{ad}(x).\ee
If an open-loop optimal control (uniquely) exists for $x\in\dbR^n$, Problem (MF-SLQ) is said to be (uniquely) {\it open-loop solvable} at $x$. Problem (MF-SLQ) is said to be (uniquely) {\it open-loop solvable} if it is (uniquely) open-loop solvable at all $x\in\dbR^n$.

\ms

(iii) A pair $(\BTh^*,v^*(\cd))\in\sS[A,\bar A,C,\bar C;B,\bar B,D,\bar D]\times L^2_\dbF(\dbR^m)$ is called a {\it closed-loop optimal strategy} if
\bel{closed-loop strategy}\ba{ll}
\ns\ds J\big(x;\Th^*\big\{X^*(\cd)-\dbE[X^*(\cd)]\big\}+\bar\Th^*\dbE[X^*(\cd)]+v^*(\cd)\big)\les J\big(x;\Th\big\{X(\cd)-\dbE[X(\cd)]\big\}+\bar\Th\dbE[X(\cd)]+v(\cd)\big),\\
\ns\ds\qq\qq\qq\qq\qq\forall(\BTh,v(\cd))\in\sS[A,\bar A,C,\bar C;B,\bar B,D,\bar D]\times L^2_\dbF(\dbR^m),~x\in\dbR^n,\ea\ee
where $X^*(\cd)\equiv X^{\BTh^*,v^*}(\cd)$, $X(\cd)\equiv X^{\BTh,v}(\cd)$ are the closed-loop state processes corresponding to $(x,\BTh^*,v^*(\cd))$ and $(x,\BTh,v(\cd))$, respectively.
If an optimal closed-loop strategy (uniquely) exists, Problem (MF-SLQ) is said to be (uniquely) {\it closed-loop solvable}.
\ms

(iv) An open-loop optimal controls $u^*(\cd\,;x)\in\sU_{ad}$ of Problem (MF-SLQ), parameterized by $x\in\dbR^n$, admits a {\it closed-loop representation}, if there exists a pair $(\BTh^*,v^*(\cd))\in\sS[A,\bar A,C,\bar C;B,\bar B$, $D,\bar D]\times L^2_\dbF(\dbR^m)$ such that for any initial state $x\in\dbR^n$,
\bel{closed-loop optimal control}
u^*(\cd)\equiv\Th^*\big\{X^*(\cd)-\dbE[X^*(\cd)]\big\}+\bar\Th^*\dbE[X^*(\cd)]+v^*(\cd)
\ee
where $X^*(\cd)\equiv X^{\BTh^*,v^*}(\cd)\in\sX[0,\i)$ is the solution to the closed-loop system \rf{closed-loop-system} corresponding to $(\BTh^*,v^*(\cd))$.
\ede

Similar to Proposition 2.5 of \cite{Li-Sun-Yong2016}, we have that $(\BTh^*,
v^*(\cd))\in\sS[A,\bar A,C,\bar C;B,\bar B,D,\bar D]\times L^2_\dbF(\dbR^m)$ is an optimal closed-loop strategy, if and only if the following condition holds:
\bel{2.12}
J\big(x;\Th^*\big\{X^*(\cd)-\dbE[X^*(\cd)]\big\}
+\bar\Th^*\dbE[X^*(\cd)]+v^*(\cd)\big)\les J\big(x;\Th^*\big\{X(\cd)-\dbE[X(\cd)]\big\}+\bar\Th^*\dbE[X(\cd)]+v(\cd)\big),\ee
for any $(x,v(\cd))\in\dbR^n\times L^2_\dbF(\dbR^m)$, where $X^*(\cd)\equiv X^{\BTh^*,v^*}(\cd)$ and $X(\cd)\equiv X^{\BTh^*,v}(\cd)$ are the closed-loop state processes corresponding to $(x,\BTh^*,v^*(\cd))$ and $(x,\BTh^*,v(\cd))$, respectively.
On the other hand, from Proposition \ref{P-2.4}, we see that under (H1), \rf{2.12} is equivalent to the following:
\bel{2.13}
J\big(x;\Th^*\big\{X^*(\cd)-\dbE[X^*(\cd)]\big\}+\bar\Th^*\dbE[X^*(\cd)]+v^*(\cd)\big)\les J\big(x;u(\cd)\big),\q\forall (x,u(\cd))\in\dbR^n\times\sU_{ad}.\ee
In general, an open-loop optimal control depends on the initial state $x\in\dbR^n$, whereas a closed-loop strategy is required to be independent of $x$. From \rf{2.13}, we see that the outcome $u^*(\cd)$ given by \rf{closed-loop optimal control}
for some closed-loop strategy $(\BTh^*,v^*(\cd))$ is an open-loop optimal control for the initial state $X^*(0)$. Hence, for Problem (MF-SLQ), the closed-loop solvability implies the open-loop solvability. The converse is also true for stochastic LQ optimal control problems in an infinite horizon without mean fields. That is to say, the open-loop and closed-loop solvabilities are equivalent (see \cite{Sun-Yong2018}).

\ms

It is natural for us to ask: Do we have such equivalence for Problem (MF-SLQ)? To answer this question, we first present the result concerning the characterization of open-loop and closed-loop solvabilities of Problem (MF-SLQ). To simplify notation,  in what follows, we denote
\bel{hat}\ba{ll}
\ds\h A=A\1n+\1n\bar A,\q\h B=B\1n+\1n\bar B,\q\h C=C\1n+\1n\bar C,\q\h D=D\1n+\1n\bar D,\q\h Q=Q\1n+\1n\bar Q,\q\h S=S\1n+\1n\bar S,\q\h R=R\1n+\1n\bar R.\ea
\ee
Note that for any $\BTh\equiv(\Th,\bar\Th)\in\sS[A,\bar A,C,\bar C;B,\bar B,D,\bar D]$ and $v(\cd)\in L^2_\dbF(\dbR^m)$, we have the closed-loop system \rf{closed-loop-system}. The cost functional \rf{cost-LQ-control} becomes
\bel{cost-LQ-control-closed-loop}\ba{ll}
\ds J^\BTh(x;v(\cd))\equiv J\big(x;\Th\big(X(\cd)-\dbE[X(\cd)]\big)+\bar\Th\dbE[X(\cd)]+v(\cd)\big)\\
\ns\ds=\1n\dbE\2n\int_0^\i\2n\bigg[\bigg\langle\2n
\begin{pmatrix}Q_\Th&S_\Th^\top\\ S_\Th&R\end{pmatrix}\2n
\begin{pmatrix}X\\ v\end{pmatrix}\1n,\1n
\begin{pmatrix}X\\ v\end{pmatrix}\2n\bigg\rangle
\1n+\1n2\bigg\langle\2n\begin{pmatrix}q_\BTh(t)\\ \rho(t)\end{pmatrix}\1n,\1n\begin{pmatrix}X\\ v\end{pmatrix}\2n\bigg\rangle\1n+\1n\bigg\langle\2n
\begin{pmatrix}\bar Q_\BTh&\bar S_\BTh^\top\\ \bar S_\BTh&\bar R\\\end{pmatrix}
\begin{pmatrix}\dbE[X]\\ \dbE[v]\end{pmatrix}\1n,\1n\begin{pmatrix}\dbE[X]\\ \dbE[v]\end{pmatrix}\2n\bigg\rangle\bigg]dt,
\ea\ee
where
\bel{Q-Th}\left\{\2n\ba{ll}
\ds Q_\Th=Q\1n+\1nS^\top\Th\1n+\1n\Th^\top S\1n+\1n\Th^\top R\Th,\q
\bar Q_\BTh=\bar Q+\h S^{\,\top}\bar\Th+\bar\Th^\top\h S+\bar\Th^\top\h R\bar\Th-S^\top\Th-\Th^\top S-\Th^\top R\Th,\\
\ns\ds S_\Th=S+R\Th,\qq\bar S_\BTh=\bar S+\h R\bar\Th-R\Th,\qq q_\BTh(\cd)=q(\cd)+\Th^\top\big(\rho(\cd)-\dbE[\rho(\cd)]\big)+\bar\Th^\top\dbE[\rho(\cd)].
\ea\right.\ee
We see that $Q_\Th,S_\Th$ depend on $\Th$, $\bar Q_\BTh,\bar S_\BTh,q_\BTh$ depend on $\BTh\equiv(\Th,\bar\Th)$, and $R(\cd),\bar R(\cd),\rho(\cd)$ are unchanged. Similar to \rf{hat}, we will denote
\bel{Hat}
\h A_\BTh=A_\Th+\bar A_\BTh,\q \h C_\BTh=C_\Th+\bar C_\BTh,\q \h Q_\BTh=Q_\Th+\bar Q_\BTh,\q \h S_\BTh=S_\Th+\bar S_\BTh.\ee

In later investigations, we will encounter the comparison between two closed-loop
strategies. Therefore, we need the following definition.

\bde{diff} \rm Let $(\BTh,v(\cd)),(\BTh',v'(\cd))\in\sS[A,\bar A,C,\bar C;B,\bar B,D,\bar D]\times L^2_\dbF(\dbR^n)$.

\ms

(i) We say that $(\BTh,v(\cd))$ and $(\BTh',v'(\cd))$ are {\it intrinsically different} if for some $x\in\dbR^n$, $X(\cd\;x,\BTh,v(\cd))\ne X(\cd\,;x,\BTh',v'(\cd))$.

\ms

(ii) We say that $(\BTh,v(\cd))$ and $(\BTh',v'(\cd))$ are {\it intrinsically the same} if for any $x\in\dbR^n$, $X(\cd\;x,\BTh,v(\cd))=X(\cd\,;x,\BTh',v'(\cd))$.

\ede

\br{different} The point that we would like to make here is that sometimes, $(\BTh,v(\cd))\ne(\BTh',v'(\cd))$. But they could be intrinsically the same. Here is such a situation. Let $(\BTh,v(\cd))\ne(\BTh',v'(\cd))$ and let $X(\cd)$ and $X'(\cd)$ be the corresponding state processes. Then (note \rf{A-Th} and \rf{hat})
$$\ba{ll}
\ns\ds A_\Th-A_{\Th'}=B(\Th-\Th'),\q\bar A_\BTh-\bar A_{\BTh'}=\h B(\bar\Th-\bar\Th')-B(\Th-\Th'),\\
\ns\ds C_\Th-C_{\Th'}=D(\Th-\Th'),\q\bar C_\BTh-\bar C_{\BTh'}=\h D(\bar\Th-\bar\Th')-D(\Th-\Th').\ea$$
Thus, if
\bel{in N(B)}\ba{ll}
\ns\ds\sR(\Th-\Th')\subseteq\sN(B)\cap\sN(D),\qq\sR(\bar\Th-\bar\Th')
\subseteq\sN(\h B)\cap\sN(\h D)\\
\ns\ds v(t)-v'(t)\in\sN(B)\cap\sN(D),\qq\dbE[v(t)-v'(t)]\in\sN(\bar B)\cap
\sN(\bar D),\qq \forall t\in[0,\i),\ea\ee
then the closed-loop systems under $(\BTh,v(\cd))$ and $(\BTh',v'(\cd))$ are the same. By definition, this means that if $(\BTh,v(\cd)),(\BTh',v'(\cd))$ are two closed-loop strategies such that \rf{in N(B)} holds, then they are intrinsically the same. If the above fails, then the two closed-loop strategies will be intrinsically different.

\ms

Note that there is another issue when we compare two closed-loop strategies, namely, the corresponding costs could be different. But, we prefer to concentrate on the difference of the corresponding state processes, which will be mainly used later.

\er

\bde{D-2.7}
The following is called a system of generalized AREs:
\bel{ARE-LQ-control}\left\{\2n\ba{ll}
\ds PA+A^\top P+C^\top PC+Q-(PB+C^\top PD+S^\top)\Si^\dag(B^\top P+D^\top PC+S)=0,\\
\ns\ds\h P\h A+\h A^{\,\top}\h P+\h C^\top P\h C+\h Q-\big(\h P\h B+\h C^\top P\h D+\h S^\top\big)\bar\Si^\dag\big(\h B^\top\h P+\h D^\top P\h C+\h S\,\big)=0,\\
\ns\ds\Si\equiv R+D^\top PD\ges0,\qq\bar\Si\equiv\h R+\h D^\top P\h D\ges0,\\
\ns\ds\sR(B^\top P+D^\top PC+S)\subseteq\sR(\Si),\qq\sR\big(\h B^\top\h P+\h D^\top P\h C+\h S\big)\subseteq\sR\big(\bar\Si\big),\ea\right.\ee
with the unknown $(P,\h P)\in\dbS^n\times\dbS^n$.
A solution pair $(P,\h P)$ to \rf{ARE-LQ-control} is said to be {\it static stabilizing} if there exists a pair $(\th,\bar\th)\in\dbR^{m\times2n}$, such that $\BTh\equiv
(\Th,\bar\Th)\in\sS[A,\bar A,C,\bar C;B,\bar B,D,\bar D]$, where
\bel{MF-stabilizer}\Th=-\Si^\dag(B^\top P+D^\top PC+S)+(I-\Si^\dag\Si)\th,\q\bar\Th=-\bar\Si^\dag\big(\h B^\top\h P+\h D^\top P\h C+\h S\,\big)+(I-\bar\Si^\dag\bar\Si)\bar\th.
\ee
\ede

We now state the following result, which is an extension of a result without mean-field terms in \cite{Sun-Yong2018} and will play an important role in the next section. The proof is postponed to Section 6.
\bt{T-2.9} \sl
Let {\rm(H1)} hold. Then the following are equivalent:

\ss

{\rm(i)}\ Problem {\rm(MF-SLQ)} is open-loop solvable.

\ss

{\rm(ii)}\ Problem {\rm(MF-SLQ)} is closed-loop solvable.

\ss

{\rm(iii)}\ The system \rf{ARE-LQ-control} admits a static stabilizing solution pair $(P,\h P)\in\dbS^n\times\dbS^n$, the BSDE on $[0,\i)$:
\bel{BSDE-eta-control}\ba{ll}
\ds-d\eta(t)\1n=\1n\big\{A^\top\1n\eta(t)\1n-\1n(B^\top\1n P\1n+\1n D^\top\1n PC\1n+\1n S)^\top\Si^\dag\big[B^\top\1n\eta(t)\1n
+\1n D^\top\1n\big(\z(t)\1n+\1n P\si(t)\big)\1n+\1n\rho(t)\big]\\
\ns\ds\qq\qq\qq+C^\top\1n\big[\z(t)\1n
+\1n P\si(t)\big]\1n+Pb(t)+q(t)\big\}dt-\z(t)dW(t),\qq t\ges0,\ea\ee
admits an adapted solution $(\eta(\cd),\z(\cd))\in\sX[0,\i)\times L^2_\dbF(\dbR^n)$ such that
\bel{range-eta}\ba{ll}
\ds B^\top\big[\eta(t)-\dbE[\eta(t)]\big]+D^\top\big[\z(t)-\dbE[\z(t)]\big]+D^\top P\big[\si(t)-\dbE[\si(t)]\big]+\rho(t)-\dbE[\rho(t)]\in\sR(\Si),\\
\ns\ds\qq\qq\qq\qq\qq\qq\qq\qq\qq\qq\qq\qq\qq\qq\qq\ae \ t\in[0,\i),~\as,
\ea\ee
and the ODE on $[0,\i)$:
\bel{BODE-bar-eta-control}\ba{ll}
\ds\dot{\bar\eta}(t)+\h A^{\,\top}\bar\eta(t)-(\h B^\top\h P+\h D^\top P\h C+\h S\,)^\top\bar\Si^\dag\big\{\h B^\top\bar\eta(t)
+\h D^\top\dbE\big[\z(t)+P\si(t)\big]+\dbE[\rho(t)]\big\}\\
\ns\ds\q +\h C^\top\dbE\big[\z(t)+P\si(t)\big]+\dbE\big[\h Pb(t)+q(t)\big]=0,\ea\ee
admits a solution $\bar\eta(\cd)\in L^2(\dbR^n)$ such that
\bel{range-bar-eta}\h B^\top\bar\eta(t)+\h D^\top\dbE[\z(t)]+\h D^\top P\dbE[\si(t)]+\dbE[\rho(t)]\in\sR\big(\bar\Si\big),\q \ae \ t\in[0,\i).
\ee
In the above case, the closed-loop optimal strategy $(\BTh^*,v^*(\cd))\equiv(
\Th^*,\bar\Th^*,v^*(\cd))$ is given by
\bel{optimal closed-loop strategy}\left\{\2n\ba{ll}
\ds\Th^*=-\Si^\dag(B^\top P+D^\top PC+S)+(I-\Si^\dag\Si)\th,\\
\ns\ds\bar\Th^*=-\bar\Si^\dag(\h B^\top\h P+\h D^\top P\h C+\h S)+(I-\bar\Si^\dag\bar\Si)\bar\th,\\
\ns\ds v^*(\cd)=\f(\cd)-\dbE[\f(\cd)]+\bar\f(\cd)+(I-\Si^\dag\Si)\big(\n(\cd)-\dbE[\n(\cd)]\big)+(I-\bar\Si^\dag\bar\Si)\bar\n(\cd),
\ea\right.\ee
with $\BTh^*\equiv(\Th^*,\bar\Th^*)\in\sS[A,\bar A,C,\bar C;B,\bar B,D,\bar D]$,
for some $\th,\bar\th\in\dbR^{m\times n}$, $\n(\cd),\bar\n(\cd)\in L^2(\dbR^m)$, and
\bel{phi-phi}\left\{\2n\ba{ll}
\ds\f(\cd)\deq-\Si^\dag\big\{B^\top\eta(\cd)+D^\top\big[\z(\cd)+P\si(\cd)\big]
+\rho(\cd)\big\},\\
\ns\ds\bar\f(\cd)\deq-\bar\Si^\dag\big\{\h B^\top\bar\eta(\cd)+\h D^\top\dbE\big[\z(\cd)+P\si(\cd)\big]+\dbE[\rho(\cd)]\big\}.\ea\right.\ee
Every open-loop optimal control $u^*(\cd)$ for the initial state $x\in\dbR^n$ admits a closed-loop representation \rf{closed-loop optimal control}, where $X^*(\cd)\in\sX[0,\i)$ is the solution to the closed-loop system \rf{closed-loop-system} under $(\BTh^*,v^*(\cd))$.
Further, the value function is given by
\bel{V(x)-control}\ba{ll}
\ds V(x)=\blan\h Px,x\bran+2\blan\bar\eta(0),x\bran+\dbE\int_0^\i\Big[\blan P\si(t),\si(t)\bran+2\blan\eta(t),b(t)-\dbE[b(t)]\bran+2\blan\bar\eta(t),\dbE[b(t)]\bran\\
\ns\ds\qq\qq\qq\qq\q+2\blan\zeta(t),\si(t)\bran-\blan\Si\big(\f(t)-\dbE[\f(t)]\big),\f(t)-\dbE[\f(t)]\bran-\blan\bar\Si\bar\f(t),\bar\f(t)\bran\Big]dt.
\ea\ee
\et

\section{Mean-Field LQ Non-Zero Sum Stochastic Differential Games}
We now return to our Problem (MF-SDG).

\subsection{Notions of Nash equilibria}

To simplify the notation, we let $m=m_1+m_2$ and denote (for $i=1,2$)
\bel{BBDD}\left\{\2n\ba{ll}
\ns\ds B=\big(B_1,B_2\big),\q \bar B=\big(\bar B_1,\bar B_2\big),\q D=\big(D_1,D_2\big),\q \bar D=\big(\bar D_1,\bar D_2\big),\\
\ns\ds S_i=\begin{pmatrix}S_{i1}\\ S_{i2}\end{pmatrix},\q \bar S_i=\begin{pmatrix}\bar S_{i1}\\ \bar S_{i2}\end{pmatrix},\q
            R_i=\begin{pmatrix} R_{i11}&R_{i12}\\ R_{i21}& R_{i22}\end{pmatrix}\equiv\begin{pmatrix} R_{i1}\\ R_{i2}\end{pmatrix},\\
\ns\ds \bar R_i=\begin{pmatrix}\bar R_{i11}&\bar R_{i12}\\ \bar R_{i21}&\bar R_{i22}\end{pmatrix}\equiv\begin{pmatrix} \bar R_{i1}\\ \bar R_{i2}\end{pmatrix},\q q_i(\cd)=\begin{pmatrix}q_{i1}(\cd)\\ q_{i2}(\cd)\end{pmatrix},\q
       \rho_i(\cd)=\begin{pmatrix}\rho_{i1}(\cd)\\ \rho_{i2}(\cd)\end{pmatrix},\q
            u(\cd)=\begin{pmatrix}u_1(\cd)\\u_2(\cd)\end{pmatrix}.
\ea\right.\ee
Then the state equation \rf{state-LQ-game} becomes
\bel{state-LQ-game-u}\left\{\2n\ba{ll}
\ds dX(t)=\big\{AX(t)+\bar A\dbE[X(t)]+Bu(t)+\bar B\dbE[u(t)]+b(t)\big\}dt\\
\ns\ds\qq\qq\q+\big\{CX(t)+\bar C\dbE[X(t)]+Du(t)+\bar D\dbE[u(t)]+\si(t)\big\}dW(t),
\qq t\ges0,\\
\ns\ds X(0)=x,\ea\right.\ee
which is of the same form as \rf{state-LQ-control}, and the cost functionals are, for $i=1,2$,
\bel{cost-LQ-game-u}\ba{ll}
\ds J_i(x;u(\cd))=\dbE\int_0^\i\bigg[\bigg\langle\2n
\begin{pmatrix}Q_i&S_i^\top\\ S_i&R_i\end{pmatrix}
\begin{pmatrix}X(t)\\ u(t)\end{pmatrix},\begin{pmatrix}X(t)\\ u(t)\end{pmatrix}\2n\bigg\rangle
+2\bigg\langle\2n\begin{pmatrix}q_i(t)\\ \rho_i(t)\end{pmatrix},\begin{pmatrix}X(t)\\ u(t)\end{pmatrix}\2n\bigg\rangle\\
\ns\ds\qq\qq\qq\qq\qq+\bigg\langle\2n
\begin{pmatrix}\bar Q_i&\bar S_i^\top\\ \bar S_i&\bar R_i\\\end{pmatrix}
\begin{pmatrix}\dbE[X(t)]\\ \dbE[u(t)]\end{pmatrix},\begin{pmatrix}\dbE[X(t)]\\ \dbE[u(t)]\end{pmatrix}
\2n\bigg\rangle\bigg]dt.\ea\ee
In order the game to make sense, we make a convention that both players at least want to keep the state $X(\cd)\equiv X(\cd\,;x,u_1(\cd),u_2(\cd))$ in $\sX[0,\infty)$ so that both cost functionals are well-defined. To guarantee this, similar to Problem (MF-SLQ), we introduce the following assumption.

\ms

{\bf(H2)} \rm System $[A,\bar A,C,\bar C;B,\bar B,D,\bar D]$  is MF-$L^2$-stablizable, i.e., $\sS[A,\bar A,C,\bar C;B,\bar B,D,\bar D]\neq\varnothing$.

\ms

Although (H2) looks the same as (H1), the meaning is different. Hypothesis (H2) provides the possibility for both players to make both cost functionals finite cooperatively. In fact, under (H2), for any $x\in\dbR^n$, by Proposition \ref{P-2.4}, the following set of all {\it admissible control pairs} is non-empty:
\bel{Uad-LQ-game-u}\sU_{ad}(x)=\Big\{u(\cd)=(u_1(\cd),u_2(\cd))\in L_\dbF^2(\dbR^m)\bigm|X(\cd)\equiv X(\cd\,;x,u(\cd))\in\sX[0,\i)\Big\}.\ee
Further, for any $u_2(\cd)\in L^2_\dbF(\dbR^{m_2})$, making use of Proposition \ref{P-2.4} again, the following holds:
\bel{Uad-LQ-game-1}\sU^1_{ad}(x)=\Big\{u_1(\cd)\in L_\dbF^2(\dbR^{m_1})\bigm|\exists u_2(\cd)\in L^2_\dbF(\dbR^{m_2}),~(u_1(\cd),u_2(\cd))\in\sU_{ad}(x)\Big\}\ne\varnothing.\ee
Likewise, for any $u_1(\cd)\in L^2_\dbF(\dbR^{m_1})$, one has
\bel{Uad-LQ-game-2}\sU^2_{ad}(x)=\Big\{u_2(\cd)\in L_\dbF^2(\dbR^{m_2})\bigm|\exists u_1(\cd)\in L^2_\dbF(\dbR^{m_1}),~(u_1(\cd),u_2(\cd))\in\sU_{ad}(x)\Big\}\ne\varnothing.\ee

We now present the following definition.
\bde{D-3.1} \rm
A $u^*(\cd)\equiv(u_1^*(\cd),u_2^*(\cd))\in\sU_{ad}(x)$ is called an {\it open-loop Nash equilibrium} of Problem (MF-SDG) for the initial state $x\in\dbR^n$ if
\bel{open-loop Nash equilibrium}\ba{ll}
\ds J_1(x;u_1^*(\cd),u_2^*(\cd))\les J_1(x;u_1(\cd),u_2^*(\cd)),\qq\forall u_1(\cd)\in\sU_{ad}^1(x),\\
\ns\ds J_2(x;u_1^*(\cd),u_2^*(\cd))\les J_2(x;u_1^*(\cd),u_2(\cd)),\qq\forall u_2(\cd)\in\sU^2_{ad}(x).
\ea\ee
\ede
For $\BTh_i\equiv(\Th_i,\bar\Th_i)\in\dbR^{m_i\times2n}, i=1,2$, we denote $\BTh\equiv\begin{pmatrix}\BTh_1\\ \BTh_2\end{pmatrix}\equiv\begin{pmatrix}
\Th_1&\bar\Th_1\\ \Th_2&\bar\Th_2\end{pmatrix}\in\dbR^{m\times 2n}$, and let
$$\ba{ll}
\ds\sS^1(\BTh_2)\deq\Big\{\BTh_1\in\dbR^{2m_1\times n}\bigm|\BTh\in\sS[A,\bar A,C,\bar C;B,\bar B,D,\bar D]\Big\},\\
\ns\ds\sS^2(\BTh_1)\deq\Big\{\BTh_2\in\dbR^{2m_2\times n}\bigm|\BTh\in\sS[A,\bar A,C,\bar C;B,\bar B,D,\bar D]\Big\}.\ea$$

It is clear that if $\BTh\in\sS[A,\bar A,C,\bar C;B,\bar B,D,\bar D]$, both $\sS^1(\BTh_2)$ and $\sS^2(\BTh_1)$ are nonempty. Similar to the optimal control problem case, any $(\BTh,v(\cd))\in\sS[A,\bar A,C,\bar C;B,\bar B,D,\bar D]\times L^2_\dbF(\dbR^m)$ is called a {\it closed-loop strategy} of Problem (MF-SDG). For any initial state $x\in\dbR^n$ and closed-loop strategy $(\BTh,v(\cd))\in\sS[A,\bar A,C,\bar C;B,\bar B,D,\bar D]\times L^2_\dbF(\dbR^m)$, we consider the following linear MF-SDE on $[0,\i)$ (recall \rf{A-Th}):
\bel{closed-loop-system-game}\left\{\2n\ba{ll}
\ds dX(t)=\big\{A_\Th X(t)+\bar A_\BTh\dbE[X(t)]+Bv(t)+\bar B\dbE[v(t)]+b(t)\big\}dt\\
\ns\ds\qq\qq+\big\{C_\Th X(t)+\bar C_\BTh\dbE[X(t)]\1n+\1n Dv(t)\1n+\1n\bar D\dbE[v(t)]\1n+\1n\si(t)\big\}dW(t),\qq t\ges0,\\
\ns\ds X(0)=x.\ea\right.\ee

By Proposition \ref{P-2.2}, \rf{closed-loop-system-game} admits a unique solution $X(\cd)\in\sX[0,\i)$. If we denote
\bel{closed-loop ui}
u_i(\cd)=\Th_i\big\{X(\cd)-\dbE[X(\cd)]\big\}+\bar\Th_i\dbE[X(\cd)]+v_i(\cd),\q i=1,2,
\ee
then \rf{closed-loop-system-game} coincides with the original state equation \rf{state-LQ-game-u}.
We call $(\BTh_i,v_i(\cd))$ a {\it closed-loop strategy} of Player $i$, and call \rf{closed-loop-system-game} the {\it closed-loop system} of the original system under closed-loop strategy $(\BTh,v(\cd))$.
Also, we call $u(\cd)\equiv(u_1(\cd),u_2(\cd))$, with $u_i(\cd)$ defined by \rf{closed-loop ui}, the {\it outcome} of the closed-loop strategy $(\BTh,v(\cd))$.

\ms

With the solution $X(\cd)\in\sX[0,\i)$ to \rf{closed-loop-system-game}, we denote, for $i=1,2$,
$$\ba{ll}
\ds J_i(x;\BTh,v(\cd))\equiv J_i(x;\BTh_1,v_1(\cd);\BTh_2,v_2(\cd))\equiv J_i\big(x;\Th\big\{X(\cd)-\dbE[X(\cd)]\big\}+\bar\Th\dbE[X(\cd)]+v(\cd)\big)\\
\ns\ds\equiv J_i\big(x;\Th_1\big\{X(\cd)-\dbE[X(\cd)]\big\}+\bar\Th_1\dbE[X(\cd)]+v_1(\cd);
\Th_2\big\{X(\cd)-\dbE[X(\cd)]\big\}+\bar\Th_2\dbE[X(\cd)]+v_2(\cd)\big).\ea$$
Similarly, we can define
$$\ba{ll}
\ds J_i\big(x;\BTh_1,v_1(\cd);u_2(\cd)\big)\equiv J_i\big(x;\Th_1\big\{X(\cd)-\dbE[X(\cd)]\big\}+\bar\Th_1\dbE[X(\cd)]+v_1(\cd);u_2(\cd)\big),\\
\ns\ds J_i\big(x;u_1(\cd);\BTh_2,v_2(\cd)\big)\equiv J_i\big(x;u_1(\cd);\Th_2\big\{X(\cd)-\dbE[X(\cd)]\big\}+\bar\Th_2\dbE[X(\cd)]+v_2(\cd)\big),\ea
\qq i=1,2.$$

We now introduce the following definition.
\bde{D-3.2} \rm
A closed-loop strategy $(\BTh^*,v^*(\cd))\in\sS[A,\bar A,C,\bar C;B,\bar B,D,\bar D]\times L^2_\dbF(\dbR^m)$ is called a {\it closed-loop Nash equilibrium} of Problem (MF-SDG) if
for any $\BTh_1\in\sS^1(\BTh_2^*)$, $\BTh_2\in\sS^2(\BTh_1^*)$, $v_1(\cd)\in L^2_\dbF(\dbR^{m_1})$ and $v_2(\cd)\in L^2_\dbF(\dbR^{m_2})$,
\bel{closed-loop Nash equilibrium}\left\{\2n\ba{ll}
\ds J_1\big(x;\BTh^*,v^*(\cd)\big)\les J_1\big(x;\BTh_1,v_1(\cd);\BTh_2^*,v_2^*(\cd)\big),\qq\forall x\in\dbR^n,\\
\ns\ds J_2\big(x;\BTh^*,v^*(\cd)\big)\les J_2\big(x;\BTh_1^*,v_1^*(\cd);\BTh_2,v_2(\cd)\big),\qq\forall x\in\dbR^n.
\ea\right.\ee
\ede

Note that on the left-hand sides of \rf{closed-loop Nash equilibrium}, the involved state is $X(\cd)\equiv X(\cd\,;x,\BTh^*,v^*(\cd))$, depending on $(\BTh^*,v^*(\cd))$. Whereas, on the right-hand sides of \rf{closed-loop Nash equilibrium}, the involved states are $X(\cd)\equiv X\big(\cd\,;x,\BTh_1,v_1(\cd);\BTh_2^*,\\v_2^*(\cd)\big)$ and $X(\cd)=X\big(\cd\,;\BTh_1^*,v_1^*(\cd);\BTh_2,v_2(\cd)\big)$ respectively, which are different in general. We emphasize that the open-loop Nash equilibrium $(u_1^*(\cd),u_2^*(\cd))$ usually depends on the initial state $x$, whereas a closed-loop Nash equilibrium $\big(\BTh^*,v^*(\cd)\big)$ is required to be independent of $x$.
It is easy to see that $\big(\BTh^*,v^*(\cd)\big)$ is a closed-loop Nash equilibrium of Problem (MF-SDG) if and only if one of the following hold:

\ms

(i)\ For any $v_1(\cd)\in L^2_\dbF(\dbR^{m_1})$ and $v_2(\cd)\in L^2_\dbF(\dbR^{m_2})$,
\bel{closed-loop Nash equilibrium-1}\ba{ll}
 J_1\big(x;\BTh^*,v^*(\cd)\big)\les  J_1\big(x;\BTh^*,v_1(\cd),v_2^*(\cd)\big),\q
 J_2\big(x;\BTh^*,v^*(\cd)\big)\les J_2\big(x;\BTh^*,v_1^*(\cd),v_2(\cd)\big);\ea\ee

(ii)\ For any $u_1(\cd)\in L^2_\dbF(\dbR^{m_1})$ and $u_2(\cd)\in L^2_\dbF(\dbR^{m_2})$,
\bel{closed-loop Nash equilibrium-2}\ba{ll}
 J_1\big(x;\BTh^*,v^*(\cd)\big)\les J_1\big(x;u_1(\cd);\BTh_2^*,v_2^*(\cd)\big),\q
 J_2\big(x;\BTh^*,v^*(\cd)\big)\les J_2\big(x;\BTh_1^*,v_1^*(\cd);u_2(\cd)\big).\ea\ee

If we denote (comparing with \rf{closed-loop Nash equilibrium})
\bel{closed-loop ui*}
u_i^*(\cd)=\Th_i^*\big\{X^*(\cd)-\dbE[X^*(\cd)]\big\}+\bar\Th_i^*\dbE[X^*(\cd)]+v_i^*(\cd),\q i=1,2,\ee
then \rf{closed-loop Nash equilibrium-2} becomes
\bel{closed-loop Nash equilibrium-3}\ba{ll}
\ds J_1\big(x;u_1^*(\cd);u_2^*(\cd)\big)\les J_1\big(x;u_1(\cd);\Th_2^*\big\{X^{u_1,v_2^*}(\cd)-\dbE[X^{u_1,v_2^*}(\cd)]\big\}
+\bar\Th_2^*\dbE[X^{u_1,v_2^*}(\cd)]+v_2^*(\cd)\big),\\
\ns\ds J_2\big(x;u_1^*(\cd);u_2^*(\cd)\big)\les J_2\big(x;\Th_1^*\big\{X^{v_1^*,u_2}(\cd)-\dbE[X^{v_1^*,u_2}(\cd)]\big\}
+\bar\Th_1^*\dbE[X^{v_1^*,u_2}(\cd)]+v_1^*(\cd);u_2(\cd)\big),\ea\ee
where $X^{u_1,v_2^*}(\cd)=X(\cd\,;x,u_1(\cd);\BTh_2^*,v_2^*(\cd))$ and $X^{v_1^*,u_2}(\cd)=X(\cd\,;x,\BTh_1^*,v_1^*(\cd);u_2(\cd))$.
Clearly, neither of the following holds in general:
$$\ba{ll}
\ns\ds u_1^*(\cd)=\Th_1^*\big\{X^{v_1^*,u_2}(\cd)-\dbE[X^{v_1^*,u_2}(\cd)]
\big\}+\bar\Th_1^*\dbE[X^{v_1^*,u_2}(\cd)]+v_1^*(\cd),\\
\ns\ds u_2^*(\cd)=\Th_2^*\big\{X^{u_1,v_2^*}(\cd)-\dbE[X^{u_1,v_2^*}(\cd)]
\big\}+\bar\Th_2^*\dbE[X^{u_1,v_2^*}(\cd)]+v_2^*(\cd).\ea$$
Hence, comparing this with \rf{open-loop Nash equilibrium}, we see that the outcome $(u_1^*(\cd),u_2^*(\cd))$ of the closed-loop Nash equilibrium $(\BTh^*,v^*(\cd))$ given by \rf{closed-loop ui*} is not necessarily an open-loop Nash equilibrium of Problem (MF-SDG) for $X^*(0)=x$.

On the other hand, if $\big(\BTh^*,v^*(\cd)\big)$ is a closed-loop Nash equilibrium of Problem (MF-SDG), we may consider the following stabilized state equation (recall \rf{A-Th}):
\bel{closed-loop-system-game-1}\left\{\2n\ba{ll}
\ds dX^{v_1,v_2}(t)=\big\{A_{\Th^*}X^{v_1,v_2}(t)+\bar A_{\BTh^*}\dbE[X^{v_1,v_2}(t)]+Bv(t)+\bar B\dbE[v(t)]+b(t)\big\}dt\\
\ns\ds\qq\qq\qq+\big\{C_{\Th^*}X^{v_1,v_2}(t)+\bar C_{\BTh^*}\dbE[X^{v_1,v_2}(t)]+Dv(t)+\bar D\dbE[v(t)]+\si(t)\big\}dW(t),\q t\ges0,\\
\ns\ds X^{v_1,v_2}(0)=x,\ea\right.\ee
with cost functionals
\bel{cost-LQ-game-closed}\ba{ll}
\ds J^{\BTh^*}_i\big(x;v_1(\cd),v_2(\cd)\big)\deq J_i\big(x;\Th_1^*\big\{X^{v_1,v_2}(\cd)-\dbE[X^{v_1,v_2}(\cd)]\big\}
+\bar\Th_1^*\dbE[X^{v_1,v_2}(\cd)]+v_1(\cd);\\
\ns\ds\qq\qq\qq\qq\qq\q\Th_2^*\big\{X^{v_1,v_2}(\cd)-\dbE[X^{v_1,v_2}(\cd)]\big\}
+\bar\Th_2^*\dbE[X^{v_1,v_2}(\cd)]+v_2(\cd)\big),\q i=1,2.\ea\ee
Then by \rf{closed-loop Nash equilibrium-1}, it is easy to see that $(v_1^*(\cd),v_2^*(\cd))$ is an open-loop Nash equilibrium of the corresponding mean-field LQ two-person non-zero sum stochastic differential games.

\ms

From the above, we see that Problems (MF-SDG) and (MF-SLQ) are essentially different in a certain sense, and we can only say that Problem (MF-SLQ) is formally a special case of Problem (MF-SDG).
\subsection{Open-loop Nash equalibria and their closed-loop representation}

In this section, we discuss the open-loop Nash equilibria for Problem (MF-SDG) in terms of MF-FBSDEs. We first have the following result.
\bt{T-3.3} \sl
Let {\rm(H2)} hold, and $x\in\dbR^n$. Then $u^*(\cd)\equiv(u_1^*(\cd),u_2^*(\cd))\in\sU_{ad}(x)$ is an open-loop Nash equilibrium of Problem (MF-SDG) for $x$ if and only if the following two conditions hold:

\ms

{\rm(i)}\ The adapted solution $(X^*(\cd),Y_i^*(\cd),Z_i^*(\cd))\1n\in\1n\sX[0,\i)\1n\times\1n\sX[0,\i)\1n\times L^2_\dbF(\dbR^n)$ to the following MF-FBSDE:
\bel{MF-FBSDE}\left\{\2n\ba{ll}
\ds dX^*(t)=\big\{AX^*(t)+\bar A\dbE[X^*(t)]+Bu^*(t)+\bar B\dbE[u^*(t)]+b(t)\big\}dt\\
\ns\ds\qq\qq\qq+\big\{CX^*(t)+\bar C\dbE[X^*(t)]+Du^*(t)+\bar D\dbE[u^*(t)]+\si(t)\big\}dW(t),\\
\ns\ds-dY_i^*(t)=\big\{A^\top Y_i^*(t)+\bar A^\top\dbE[Y_i^*(t)]+C^\top Z_i^*(t)+\bar C^\top\dbE[Z_i^*(t)]+Q_iX^*(t)+\bar Q_i\dbE[X^*(t)]\\
\ns\ds\qq\qq\qq+S_i^\top u^*(t)+\bar S_i^\top\dbE[u^*(t)]+q_i(t)\big\}dt-Z_i(t)dW(t),\qq t\ges0,\q i=1,2,\\
\ns\ds X^*(0)=x,\ea\right.\ee
satisfies the following stationarity condition:
\bel{stationarity condition}\ba{ll}
\ds B_i^\top Y_i^*(t)+\bar B_i^\top\dbE[Y_i^*(t)]+D_i^\top Z_i^*(t)+\bar D_i^\top\dbE[Z_i^*(t)]+S_{ii}X^*(t)+\bar S_{ii}\dbE[X^*(t)]\\
\ns\ds\qq +R_{ii}u^*(t)+\bar R_{ii}\dbE[u^*(t)]+\rho_{ii}(t)=0,\qq\ae \ t\in[0,\i),~\as,\q i=1,2.\ea\ee

{\rm(ii)}\ The maps $u_1(\cd)\mapsto J_1(x;u_1(\cd),u_2(\cd))$ and $u_2(\cd)\mapsto J_2(x;u_1(\cd),u_2(\cd))$ are convex, i.e.,
\bel{convexity condition}\ba{ll}
\ds\dbE\int_0^\i\Big[\blan Q_iX^0_i(t),X^0_i(t)\bran+2\blan S_{ii}X^0_i(t),u_i(t)\bran+\blan R_{iii}u_i(t),u_i(t)\bran+\blan\bar Q_i\dbE[X^0_i(t)],\dbE[X^0_i(t)]\bran\\
\ns\ds\qq+2\blan\bar S_{ii}\dbE[X^0_i(t)],\dbE[u_i(t)]\bran\1n+\1n\blan\bar R_{iii}\dbE[u_i(t)],\dbE[u_i(t)]\bran\Big]dt\1n\ges\1n0,\q\forall u(\cd)\1n\equiv\1n(u_1(\cd),u_2(\cd))\in\sU_{ad}(x),\ea\ee
where $X^0_i(\cd)\in\sX[0,\i)$ is the solution to the following homogeneous controlled MF-SDE:
\bel{homogeneous controlled MF-SDE}\left\{\2n\ba{ll}
\ds dX^0_i(t)=\big\{AX^0_i(t)+\bar A\dbE[X^0_i(t)]+B_iu_i(t)+\bar B_i\dbE[u_i(t)]\big\}dt\\
\ns\ds\qq\qq\q+\big\{CX^0_i(t)+\bar C\dbE[X^0_i(t)]+D_iu_i(t)+\bar D_i\dbE[u_i(t)]\big\}dW(t),\q t\ges0,\\
\ns\ds X^0_i(0)=0.\ea\right.\ee

\et
\proof
For given $x\in\dbR^n$ and $u^*(\cd)\in\sU_{ad}(x)$, let $(X^*(\cd),Y_1^*(\cd),Z_1^*(\cd))$  be the solution to \rf{MF-FBSDE} with $i=1$. For any $u_1(\cd)\in L^2_\dbF(\dbR^{m_1})$ and $\e\in\dbR$, let $X^\e(\cd)$ be the solution to the following perturbed state equation:
$$\left\{\2n\ba{ll}
\ds dX^\e=\big\{AX^\e+\bar A\dbE[X^\e]+B_1(u_1^*+\e u_1)+\bar B_1(\dbE[u_1^*]+\e \dbE[u_1])+B_2u_2^*+\bar B_2\dbE[u_2^*]+b\big\}dt\\
\ns\ds\qq\qq+\big\{CX^\e\2n+\1n\bar C\dbE[X^\e]\1n+\1n D_1(u_1^*\1n+\1n\e u_1)\1n+\1n\bar D_1(\dbE[u_1^*]\1n+\1n\e\dbE[u_1])\1n+\1n D_2u_2^*\1n+\1n\bar D_2\dbE[u_2^*]\1n+\1n\si\big\}dW(t),\q t\ges0,\\
\ns\ds X^\e(0)=x.\ea\right.$$
Obviously, we have $X^\e(\cd)=X^*(\cd)+\e X^0_1(\cd)$, where $X^0_1(\cd)$ is the solution to \rf{homogeneous controlled MF-SDE} with $i=1$. Thus,
$$\ba{ll}
\ds J_1(x;u_1^*(\cd)+\e u_1(\cd),u_2^*(\cd))-J_1(x;u_1^*(\cd),u_2^*(\cd))\\
\ns\ds=2\e\dbE\int_0^\i\Big[\blan Q_1X^*+S_1^\top u^*+q_1,X^0_1\bran+\blan S_{11}X^*+R_{11}u^*+\rho_{11},u_1\bran\Big]dt\\
\ns\ds\q+\e^2\dbE\int_0^\i\Bigg\langle\2n
\begin{pmatrix}Q_1&S_{11}^\top\\ S_{11}&R_{111}\end{pmatrix}
\begin{pmatrix}X^0_1\\ u_1\end{pmatrix},\begin{pmatrix}X^0_1\\ u_1\end{pmatrix}\2n\Bigg\rangle dt
+2\e\dbE\int_0^\i\Big[\blan\bar Q_1\dbE[X^*]+\bar S_1^\top\dbE[u^*],\dbE[X^0_1]\bran\\
\ns\ds\q+\blan\bar S_{11}\dbE[X^*]+\bar R_{11}\dbE[u^*],\dbE[u_1]\bran\Big]dt
+\e^2\dbE\int_0^\i\Bigg\langle\2n
\begin{pmatrix}\bar Q_1&\bar S_{11}^\top\\ \bar S_{11}&\bar R_{111}\end{pmatrix}
\begin{pmatrix}\dbE[X^0_1]\\ \dbE[u_1]\end{pmatrix},\begin{pmatrix}\dbE[X^0_1]\\ \dbE[u_1]\end{pmatrix}\2n\Bigg\rangle dt.\ea$$
Since $[A,\bar A,C,\bar C;B,\bar B,D,\bar D]$ is MF-$L^2$-stabilizable, by Proposition \ref{P-2.2}, the forward MF-SDE in \rf{MF-FBSDE} admits a unique solution $X^*(\cd)\in\sX[0,\i)$, and the MF-BSDEs in \rf{MF-FBSDE} admit unique solutions $(Y_i^*(\cd),Z_i^*(\cd))\in\sX[0,\i)\times L^2_\dbF(\dbR^n)$, $i=1,2$, respectively.
Applying It\^{o}'s formula to $\lan Y_1^*(\cd),X^0_1(\cd)\ran$, we have
$$\ba{ll}
\ds\dbE\big[\blan Y_1^*(s),X^0_1(s)\bran\big]=\dbE\int_0^s\Big[-\blan A^\top Y_1^*+\bar A\dbE[Y_1^*]+C^\top Z_1^*+\bar C^\top\dbE[Z_1^*]+Q_1X^*+\bar Q_1\dbE[X^*]\\
\ns\ds\qq\qq\qq\qq\qq\qq+S_1^\top u^*+\bar S_1^\top\dbE[u^*]+q_1,X^0_1\bran+\blan Y_1^*,AX^0_1+\bar A\dbE[X^0_1]+B_1u_1+\bar B_1\dbE[u_1]\bran\\
\ns\ds\qq\qq\qq\qq\qq\qq+\blan CX^0_1+\bar C\dbE[X^0_1]+D_1u_1+\bar D_1\dbE[u_1],Z_1^*\bran\Big]dt\\
\ns\ds\qq\qq\qq\qq=\dbE\int_0^s\Big[\blan B_1^\top Y_1^*+\bar B_1\dbE[Y_1^*]+D_1^\top Z_1^*+\bar D_1\dbE[Z_1^*],u_1\bran\\
\ns\ds\qq\qq\qq\qq\qq\qq-\blan Q_1X^*+\bar Q_1\dbE[X^*]+S_1^\top u^*+\bar S_1^\top\dbE[u^*]+q_1,X^0_1\bran\Big]dt,\qq\forall s\ges0.\ea$$
Noting that
$$\lim_{s\to\i}\big|\dbE\big[\blan Y_1^*(s),X^0_1(s)\bran\big]\big|^2\les\lim_{s\to\i}\dbE[|Y_1^*(s)|^2]\dbE[|X^0_1(s)|^2]=0,$$
then letting $s\to\i$ we have
$$\ba{ll}
\ds \dbE\int_0^\i\Big[\blan B_1^\top Y_1^*+\bar B_1\dbE[Y_1^*]+D_1^\top Z_1^*+\bar D_1\dbE[Z_1^*],u_1\bran\Big]dt\\
\ns\ds =\dbE\int_0^\i\Big[\blan Q_1X^*+\bar Q_1\dbE[X^*]+S_1^\top u^*+\bar S_1^\top\dbE[u^*]+q_1,X^0_1\bran\Big]dt.\ea$$
Combining the above equalities, we obtain
$$\ba{ll}
\ds J_1\big(x;u_1^*(\cd)+\e u_1(\cd),u_2^*(\cd)\big)-J_1\big(x;u_1^*(\cd),u_2^*(\cd)\big)=2\e\dbE\int_0^\i\Big[\blan B_1^\top Y_1^*+\bar B_1^\top\dbE[Y_1^*]
 +D_1^\top Z_1^*+\bar D_1^\top\dbE[Z_1^*]\\
\ns\ds\q+S_{11}X^*+\bar S_{11}\dbE[X^*]+R_{11}u^*+\bar R_{11}\dbE[u^*]+\rho_{11},u_1\bran\Big]dt+\e^2\dbE\int_0^\i\Big[\blan Q_1X^0_1,X^0_1\bran+2\blan S_{11}X^0_1,u_1\bran\\
\ns\ds\q+\blan R_{111}u_1,u_1\bran+\blan \bar Q_1\dbE[X^0_1],\dbE[X^0_1]\bran+2\blan \bar S_{11}\dbE[X^0_1],\dbE[u_1]\bran+\blan\bar R_{111}\dbE[u_1],\dbE[u_1]\bran\Big]dt.\ea$$
It follows that
$$J_1\big(x;u_1^*(\cd),u_2^*(\cd)\big)\les J_1\big(x;u_1^*(\cd)+\e u_1(\cd),u_2^*(\cd)\big),\q\forall u_1(\cd)\in L^2_\dbF(\dbR^{m_1}),\q\forall\e\in\dbR,$$
if and only if \rf{stationarity condition}--\rf{convexity condition} hold for $i=1$. Similarly,
$$J_2\big(x;u_1^*(\cd),u_2^*(\cd)\big)\les J_2\big(x;u_1^*(\cd),u_2^*(\cd)+\e u_2(\cd)\big),\q\forall u_2(\cd)\in L^2_\dbF(\dbR^{m_2}),\q \forall\e\in\dbR,$$
if and only if \rf{stationarity condition}--\rf{convexity condition} hold for $i=2$. Combining the above two cases, the theorem is proved.
\endpf

\ms

Note that \rf{MF-FBSDE} is a system of coupled MF-FBSDEs, with one forward equation and two backward equations, with the coupling through the relation \rf{stationarity condition}. Our next task is to investigate the solvability of \rf{MF-FBSDE}--\rf{stationarity condition}. To this end, we introduce the following notation (recall \rf{BBDD}):
\bel{BA}\ba{ll}
\ds\BA=\begin{pmatrix}A&0\\0&A\end{pmatrix},~\bar\BA=\begin{pmatrix}\bar A&0\\0&\bar A\end{pmatrix},~\BC=\begin{pmatrix}C&0\\0&C\end{pmatrix},~\bar\BC=\begin{pmatrix}\bar C&0\\0&\bar C\end{pmatrix}\in\dbR^{2n\times 2n};\\
\ns\ds \BB=\begin{pmatrix}B&0\\0&B\end{pmatrix},~\bar\BB=\begin{pmatrix}\bar B&0\\0&\bar B\end{pmatrix},~\BD=\begin{pmatrix}D&0\\0&D\end{pmatrix},~\bar\BD=\begin{pmatrix}\bar D&0\\0&\bar D\end{pmatrix}\in\dbR^{2n\times 2m};\\
\ns\ds\BQ=\begin{pmatrix}Q_1&0\\0&Q_2\end{pmatrix},~\bar\BQ=\begin{pmatrix}\bar Q_1&0\\0&\bar Q_2\end{pmatrix}\in\dbS^{2n};\q\BS=\begin{pmatrix}S_1&0\\0&S_2\end{pmatrix},~
\bar\BS=\begin{pmatrix}\bar S_1&0\\0&\bar S_2\end{pmatrix}\in\dbR^{2m\times 2n};\\
\ns\ds\BR=\begin{pmatrix}R_1&0\\0&R_2\end{pmatrix},~\bar\BR=\begin{pmatrix}\bar R_1&0\\0&\bar R_2\end{pmatrix}\in\dbS^{2m};~
 q(\cd)\1n=\1n\begin{pmatrix}q_1(\cd)\\q_2(\cd)\end{pmatrix}\1n\in\1n L^2_\dbF(\dbR^{2n});~
\rho(\cd)\1n=\1n\begin{pmatrix}\rho_1(\cd)\\\rho_2(\cd)\end{pmatrix}\1n\in\1n L^2_\dbF(\dbR^{2m}).\ea\ee
With the above notations, MF-FBSDEs \rf{MF-FBSDE} can be rewritten as
\bel{MF-FBSDE-2}\left\{\2n\ba{ll}
\ds dX^*(t)=\big\{AX^*+\bar A\dbE[X^*]+Bu^*+\bar B\dbE[u^*]+b\big\}dt\\
\ns\ds\qq\qq\qq+\big\{CX^*+\bar C\dbE[X^*]+Du^*+\bar D\dbE[u^*]+\si\big\}dW(t),\\
\ns\ds-d\BY^*(t)=\big\{\BA^\top\BY^*+\bar\BA^\top\dbE[\BY^*]+\BC^\top \BZ^*+\bar \BC^\top\dbE[\BZ^*]+\BQ\BI_n X^*+\bar\BQ\BI_n\dbE[X^*]\\
\ns\ds\qq\qq\qq+\BS^\top\BI_m u^*+\bar\BS^\top\BI_m\dbE[u^*]+q\big\}dt-\BZ^*dW(t),\q t\ges0,\\
\ns\ds X^*(0)=x,\ea\right.\ee
and the stationarity condition \rf{stationarity condition} can be written as
\bel{stationarity condition-BJ}\ba{ll}
\ds\BJ^\top\1n\big\{\BB^\top\BY^*\2n+\1n\bar\BB^\top\dbE[\BY^*]\1n+\1n
\BD^\top\BZ^*\1n+\1n\bar\BD^\top\dbE[\BZ^*]\1n+\1n\BS\BI_nX^*\1n+\1n\bar \BS\BI_n\dbE[X^*]\1n+\1n\BR\BI_mu^*\1n+\1n\bar\BR\BI_m\dbE[u^*]\1n+\1n\rho\big\}=0,\\
\ns\ds\qq\qq\qq\qq\qq\qq\qq\qq\qq\qq\ae \ t\in[0,\i),~\as,\ea\ee
where
\bel{IJ}\left\{\2n\ba{ll}
\ds\BY^*(\cd)=\begin{pmatrix}Y_1^*(\cd)\\Y_2^*(\cd)\end{pmatrix}\in\sX[0,\i)^2,\q
 \BZ^*(\cd)=\begin{pmatrix}Z_1^*(\cd)\\Z_2^*(\cd)\end{pmatrix}\in L^2_\dbF(\dbR^{2n}),\q \BI_n=\begin{pmatrix}I_{n\times n}\\I_{n\times n}\end{pmatrix}\in\dbR^{2n\times n},\\
\ns\ds\BI_m=\begin{pmatrix}I_{m\times m}\\I_{m\times m}\end{pmatrix}\in\dbR^{2m\times m},\q\BJ=\begin{pmatrix}I_{m_1\times m_1}&0\\0&0\\0&0\\0&I_{m_2\times m_2}\end{pmatrix}
 =\begin{pmatrix}I_{m_1\times m_1}&0_{m_1\times m_2}\\0_{m_2\times m_1}&0_{m_2\times m_2}\\0_{m_1\times m_1}&0_{m_1\times m_2}\\0_{m_2\times m_1}&I_{m_2\times m_2}\end{pmatrix}
 \in\dbR^{2m\times m}.\ea\right.\ee
%

%
Now, inspired by \cite{Sun-Yong2019, Yong2013}, we may obtain a closed-loop representation of open-loop Nash equilibria, by which we mean an open-loop Nash equilibrium $u^{**}(\cd)\equiv(u_1^{**}(\cd),u_2^{**} (\cd))$ admits a form \rf{closed-loop ui} for some $(\BTh^{**},v^{**}(\cd))\in\sS[A,\bar A,C,\bar C;B,\bar B,D,\bar D]\times L^2_\dbF(\dbR^m)$. We state the result here.
\bt{T-3.4} \sl
Let {\rm(H2)} hold. Then an open-loop Nash equilibria $u^{**}(\cd)\equiv(u_1^{**}(\cd),u_2^{**}(\cd))\in\sU_{ad}(x)$ of Problem {\rm(MF-SDG)} admits the following closed-loop representation
\bel{closed-loop ui**}
u_i^{**}(\cd)=\Th_i^{**}\big\{X(\cd)-\dbE[X(\cd)]\big\}+\bar\Th_i^{**}\dbE[X(\cd)]+v_i^{**}(\cd),\q i=1,2,\ee
with $\BTh^{**}\equiv\begin{pmatrix}\BTh_1^{**}\\ \BTh_2^{**}\end{pmatrix}\equiv\begin{pmatrix}\Th_1^{**}&\bar\Th_1^{**}\\ \Th_2^{**}&\bar\Th_2^{**}\end{pmatrix}\in\sS[A,\bar A,C,\bar C;B,\bar B,D,\bar D]$ and $v^{**}(\cd)\equiv\begin{pmatrix}v_1^{**}(\cd)\\ v_2^{**}(\cd)\end{pmatrix}\in L^2_\dbF(\dbR^m)$ if and only if the following hold:

\ms

{\rm(i)}\ The convexity condition (\ref{convexity condition}) holds for $i=1,2$.

\ms

{\rm(ii)}\ The solution pair $(\BP,\h\BP)\deq\begin{pmatrix}P_1&\h P_1\\ P_2&\h P_2\end{pmatrix}\in\dbR^{2n\times2n}$ to the system of coupled AREs:
\bel{ARE-LQ-game-BP}\left\{\2n\ba{ll}
\ds\BP A+\BA^\top\BP+\BC^\top\BP C+\BQ\BI_n-\big(\BP B+\BC^\top\BP D+\BS^\top\BI_m\big)\mathbf{\Si}^{-1}\BJ^\top\big(\BB^\top\BP+\BD^\top\BP C+\BS\BI_n\big)=0,\\
\ns\ds\h\BP\h A+\h\BA^\top\h\BP+\h\BC^\top\BP\h C+\h\BQ\BI_n-\big(\h\BP\h B+\h\BC^
\top\BP\h D+\h\BS^\top\BI_m\big)\bar\BSi^{-1}\BJ^\top\big(\h\BB^\top\h\BP+\h\BD^\top\BP\h C+\h\BS\BI_n\big)=0,
\ea\right.\ee
where $\h\BA=\BA+\bar\BA$ and $\h\BC,\h\BQ,\h\BS$ are defined similarly as \rf{hat} and
$$\BSi\deq\BJ^\top\big(\BR\BI_m+\BD^\top\BP D\big),\q \bar\BSi\deq\BJ^\top\big(\h\BR\BI_m+\h\BD^\top\BP\h D\big)\in\dbR^{m\times m}$$
are both invertible such that $\BTh^{**}\equiv(\Th^{**},\bar\Th^{**})\equiv\begin{pmatrix}\Th_1^{**}&\bar\Th_1^{**}\\ \Th_2^{**}&\bar\Th_2^{**}\end{pmatrix}\in\dbR^{m\times2n}$ defined by
\bel{Th-Th*-v*}\Th^{**}\deq-\BSi^{-1}\BJ^\top\big(\BB^\top\BP+\BD^\top\BP C+\BS\BI_n\big),\qq
\bar\Th^{**}\deq-\bar\BSi^{-1}\BJ^\top\big(\h\BB^\top\h\BP+\h\BD^\top\BP\h C+\h\BS\BI_n\big),\ee
stabilizes the system $[A,\bar A,C,\bar C;B,\bar B,D,\bar D]$. Also,
\bel{v*}\ba{ll}
\ns\ds v^{**}(\cd)=-\mathbf{\Si}^{-1}\BJ^\top\big\{\BB^\top(\eta-\dbE[\eta])+\BD^\top(\z-\dbE[\z])+\BD^\top\BP(\si-\dbE[\si])+\rho-\dbE[\rho]\big\}\\
\ns\ds\qq\qq-\bar{\mathbf{\Si}}^{-1}\BJ^\top\big\{\h\BB^\top\bar\eta+\h\BD^\top\dbE[\z]+\h\BD^\top\BP
\dbE[\si]+\dbE[\rho]\big\}\in L^2_\dbF(\dbR^m),\ea\ee
where  $\left(\eta(\cd)\deq\begin{pmatrix}\eta_1(\cd)\\\eta_2(\cd)\end{pmatrix},\z(\cd)\deq\begin{pmatrix}\zeta_1(\cd)\\\zeta_2(\cd)\end{pmatrix}\right)\in\sX[0,\i)^2\times L^2_\dbF(\dbR^{2n})$ is an adapted solution to BSDE:
\bel{BSDE-eta-game-BJ}\ba{ll}
\ds-d\eta(t)=\Big\{\BA^\top\eta(t)-\big(\BP B+\BC^\top\BP D+\BS^\top\BI_m\big)\BSi^{-1}
\BJ^\top\big\{\BB^\top\eta(t)+\BD^\top\big[\z(t)+\BP\si(t)\big]+\rho(t)\big\}\\
\ns\ds\qq\qq\q+\BC^\top\big[\z(t)+\BP\si(t)\big]+\BP b(t)+q(t)\Big\}dt-\zeta(t)dW(t),\qq t\ges0,
\ea\ee
and $\bar\eta(\cd)\deq\begin{pmatrix}\bar\eta_1(\cd)\\\bar\eta_2(\cd)\end{pmatrix}\in L^2(\dbR^{2n})$ is a solution to ODE:
\bel{BODE-bar-eta-BJ}\ba{ll}
\ns\ds \dot{\bar\eta}+\h\BA^\top\bar\eta+\h\BC^\top\dbE[\z+\BP\si]+\h\BP\dbE[b]+\dbE[q]\\
\ns\ds\qq -\big(\h\BP\h B+\h\BC^\top\BP\h D+\h\BS^\top\BI_m\big)\bar\BSi^{-1}\big\{\h\BB^\top\bar\eta+\h\BD^\top\dbE[\z+\BP\si]
+\dbE[\rho]\big\}=0,\qq t\ges0.\ea\ee
In such a case, the open-loop Nash equilibrium $u^{**}(\cd)$ admits a closed-loop representation \rf{closed-loop ui*} with $(\BTh^{**},v^{**}(\cd))$ given by the above.

\et

\proof We take the following {\it ansatz}:
\bel{BY}\BY^*(\cd)=\BP\big(X^*(\cd)-\dbE[X^*(\cd)]\big)+\h \BP\dbE[X^*(\cd)]+\eta(\cd)-\dbE[\eta(\cd)]+\bar\eta(\cd),\ee
where
\bel{BP}\BP\deq\begin{pmatrix}P_1\\ P_2\end{pmatrix},\q \h\BP\deq\begin{pmatrix}\h P_1\\\h P_2\end{pmatrix},\q
\eta(\cd)\triangleq\begin{pmatrix}\eta_1(\cd)\\\eta_2(\cd)\end{pmatrix},\q
\bar\eta(\cd)\triangleq\begin{pmatrix}\bar\eta_1(\cd)\\\bar\eta_2(\cd)\end{pmatrix}.\ee
Here, $P_i,\h P_i\in\dbR^{n\times n},i=1,2$, and $(\eta(\cd),\z(\cd))$ is an adapted solution to the following BSDE on $[0,\i)$:
\bel{BSDE-eta-game-form}-d\eta(t)=\a(t)dt-\z(t)dW(t),\q t\ges0,\ee
where $\a:[0,\i)\times\Omega\to\dbR^{2n}$ is undetermined, and $\bar\eta(\cd)$ is a deterministic differentiable function from $[0,\i)$ to $\dbR^n$. Applying It\^o's formula to \rf{BY}, noting \rf{MF-FBSDE-2}, we get
\bel{Ito-1}\ba{ll}
\ds-\big\{\BA^\top\big(\BY^*-\dbE[\BY^*]\big)+\h\BA^\top\dbE[\BY^*]+\BC^\top\big(\BZ^*-\dbE[\BZ^*]\big)
+\h\BC^\top\dbE[\BZ^*]+\BQ\BI_n\big(X^*-\dbE[X^*]\big)\\
\ns\ds\q+\h\BQ\BI_n\dbE[X^*]+\BS^\top\BI_m\big(u^{**}-\dbE[u^{**}]\big)+\h\BS^\top\BI_m\dbE[u^{**}]+ q\big\}dt+\BZ^*dW(t)=d\BY^*\\
\ns\ds=\big\{\BP A\big(X^*-\dbE[X^*]\big)+\BP B\big(u^{**}-\dbE[u^{**}]\big)+\BP\big(b-\dbE[b]\big)+\h\BP\h A\dbE[X^*]\\
\ns\ds\qq+\h\BP\h B\dbE[u^{**}]+\h\BP\dbE[b]+\dot{\bar\eta}-\big(\a-\dbE[\a]\big)\big\}dt+\big\{\BP C\big(X^*\2n-\1n\dbE[X^*]\big)\1n\\
\ns\ds\qq+\BP\h C\dbE[X^*]\1n+\1n\BP D\big(u^{**}\2n-\1n\dbE[u^{**}]\big)\1n+\1n\BP\h D\dbE[u^{**}]\1n+\1n\BP\si\1n+\1n\z\big\}dW(t).\ea\ee
Hence, by comparing the diffusion terms, we should have
\bel{BZ}\ba{ll}
\BZ^*=\BP C\big(X^*-\dbE[X^*]\big)+\BP\h C\dbE[X^*]+\BP D\big(u^{**}-\dbE[u^{**}]\big)+\BP\h D\dbE[u^{**}]+\BP\si+\z,\q\as\ea\ee
The stationarity condition \rf{stationarity condition-BJ} then becomes
$$\ba{ll}
\ds 0=\BJ^\top\big\{\BB^\top \BY^*+\bar \BB^\top\dbE[\BY^*]+\BD^\top \BZ^*+\bar \BD^\top\dbE[\BZ^*]+\BS\BI_nX^*+\bar \BS\BI_n\dbE[X^*]+\BR\BI_mu^{**}+\bar \BR\BI_m\dbE[u^{**}]+\rho\big\}\\
\ns\ds=\BJ^\top\big\{\BB^\top\big(\BY^*-\dbE[\BY^*]\big)+\h\BB^\top\dbE[\BY^*]+\BD^\top\big(\BZ^*
 -\dbE[\BZ^*]\big)+\h\BD^\top\dbE[\BZ^*]\\
\ns\ds\qq+\BS\BI_n\big(X^*-\dbE[X^*]\big)+\h\BS\BI_n\dbE[X^*]+\BR\BI_m\big(u^{**}-\dbE[u^{**}]
\big)+\h\BR\BI_m\dbE[u^{**}]+\rho\big\}\\
\ns\ds=\BJ^\top\big\{(\BB^\top\BP+\BD^\top\BP C+\BS\BI_n)\big(X^*-\dbE[X^*]\big)+(\h\BB^\top\h\BP+\h\BD^\top\BP\h C
  +\h\BS\BI_n)\dbE[X^*]\\
\ns\ds\qq+(\BR\BI_m+\BD^\top\BP D)\big(u^{**}-\dbE[u^{**}]\big)+(\h\BR\BI_m+\h\BD^\top\BP\h D)\dbE[u^{**}]+\BB^\top\big(\eta-\dbE[\eta]\big)\\
\ns\ds\qq+\h\BB^\top\bar\eta+\BD^\top\big(\z-\dbE[\z]\big)+\h\BD^\top\dbE[\z]+\BD^\top
\BP\big(\si-\dbE[\si]\big)+\h\BD^\top\BP\dbE[\si]+\rho\big\},\q\as\ea$$
Applying $\dbE[\,\cd\,]$ to the above, and assuming that
$$\BSi\deq\BJ^\top\big(\BR\BI_m+\BD^\top\BP D\big),\q \bar\BSi\deq\BJ^\top\big(\h\BR\BI_m+\h\BD^\top\BP\h D\big)\in\dbR^{m\times m},$$
are invertible, we obtain
\bel{Eu*}
\dbE[u^{**}]=-\bar\BSi^{-1}\BJ^\top(\h\BB^\top\h\BP+\h\BD^\top\BP\h C+\h\BS\BI_n)\dbE[X^*]-\bar\BSi^{-1}\BJ^\top\big\{\h\BB^\top\bar\eta+\h\BD^\top\dbE[\z]
         +\h\BD^\top\BP\dbE[\si]+\dbE[\rho]\big\},\ee
\bel{u*-Eu*}\ba{ll}
\ds u^{**}-\dbE[u^{**}]=-\BSi^{-1}\BJ^\top\big(\BB^\top\BP+\BD^\top\BP C+\BS\BI_n\big)(X^*-\dbE[X^*])\\
\ns\ds\qq\qq\qq-\BSi^{-1}\BJ^\top\big\{\BB^\top(\eta-\dbE[\eta])+\BD^\top(\z-\dbE[\z])
             +\BD^\top\BP(\si-\dbE[\si])+\rho-\dbE[\rho]\big\},\ea\ee
and thus
\bel{u*}\ba{ll}
\ds u^{**}=-\mathbf{\Si}^{-1}\BJ^\top\big(\BB^\top\BP+\BD^\top\BP C+\BS\BI_n\big)(X^*(t)-\dbE[X^*])
 -\bar{\mathbf{\Si}}^{-1}\BJ^\top(\h\BB^\top\h\BP+\h\BD^\top\BP\h C+\h\BS\BI_n)\dbE[X^*]\\
\ns\ds\qq\ -\mathbf{\Si}^{-1}\BJ^\top\big\{\BB^\top(\eta-\dbE[\eta])+\BD^\top(\z-\dbE[\z])+\BD^\top\BP(\si-\dbE[\si])+\rho-\dbE[\rho]\big\}\\
\ns\ds\qq\ -\bar{\mathbf{\Si}}^{-1}\BJ^\top\big\{\h\BB^\top\bar\eta+\h\BD^\top\dbE[\z]+\h\BD^\top\BP\dbE[\si]+\dbE[\rho]\big\},\q\as\ea\ee
Now, comparing the drift terms in \rf{Ito-1}, one gets
$$\ba{ll}
\ds0=-(\a-\dbE[\a])+\big(\BP A+\BA^\top\BP+\BC^\top\BP C+\BQ\BI_n\big)\big(X^*-\dbE[X^*]\big)+\big(\h\BP\h A+\h\BA^\top\h\BP+\h\BC^\top\BP\h C+\h\BQ\BI_n\big)\dbE[X^*]\\
\ns\ds\q+\big(\BP B+\BC^\top\BP D+\BS^\top\BI_m\big)(u^{**}-\dbE[u^{**}]\big)+\big(\h\BP\h B+\h\BC^\top\BP\h D+\h\BS^\top\BI_m\big)\dbE[u^{**}]+\BP\big(b-\dbE[b]\big)+\h\BP\dbE[b]+\dot{\bar\eta}\\
\ns\ds\q+\BA^\top\big(\eta-\dbE[\eta]\big)+\h\BA^\top\bar\eta+\BC^\top\big\{\BP(\si
-\dbE[\si]\big)+\z-\dbE[\z]\big\}+\h\BC^\top\big(\BP\dbE[\si]+\dbE[\z]\big)+q\\
\ns\ds=\big\{\BP A+\BA^\top\BP+\BC^\top\BP C+\BQ\BI_n-\big(\BP B+\BC^\top\BP D+\BS^\top\BI_m\big)\BSi^{-1}\BJ^\top\big(\BB^\top\BP+\BD^\top\BP C+\BS\BI_n\big)\big\}(X^*-\dbE[X^*])\\
\ns\ds\q+\Big\{\h\BP\h A+\h\BA^\top\h\BP+\h\BC^\top\BP\h C+\h\BQ\BI_n-\big(\h\BP\h B+\h\BC^\top\BP\h D+\h\BS^\top\BI_m\big)\bar\BSi^{-1}\BJ^\top\big(\h\BB^\top\h\BP
  +\h\BD^\top\BP\h C+\h\BS\BI_n\big)\Big\}\dbE[X^*]\\
\ns\ds\q-\big(\BP B+\BC^\top\BP D+\BS^\top\BI_m\big)\BSi^{-1}\BJ^\top\big\{\BB^\top\big(\eta-\dbE[\eta]\big)+\BD^\top\big(\z
-\dbE[\z]\big)
+\BD^\top\BP\big(\si-\dbE[\si]\big)+\rho-\dbE[\rho]\big\}\\
\ns\ds\q-(\h\BP\h B+\h\BC^\top\BP\h D+\h\BS^\top\BI_m)\bar\BSi^{-1}\BJ^\top\big\{\h\BB^\top\bar\eta
+\h\BD^\top\dbE[\z]+\h\BD^\top\BP\dbE[\si]+\dbE[\rho]\big\}+\BP(b-\dbE[b])+\h\BP
 \dbE[b]\\
\ns\ds\q+\dot{\bar\eta}+\BA^\top\big(\eta-\dbE[\eta]\big)+\h\BA^\top\bar\eta
+\BC^\top\big[\z-\dbE[\z]+\BP\big(\si-\dbE[\si]\big)\big]
 +\h\BC^\top\big(\dbE[\z]+\BP\dbE[\si]\big)+q-(\a-\dbE[\a]).\ea$$
This suggests that $(\BP,\h\BP)$ should be the solution to \rf{ARE-LQ-game-BP} and $(\eta(\cd),\z(\cd),\bar\eta(\cd))$ should satisfy
\bel{alpha}\ba{ll}
\ds 0\1n=\1n\dbE[\a]-\a-\big(\BP B+\BC^\top\BP D+\BS^\top\BI_m\big)\BSi^{-1}\BJ^\top\big\{\BB^\top\big(\eta-\dbE[\eta]\big)
+\BD^\top\big(\z-\dbE[\z]\big)+\BD^\top\BP\big(\si-\dbE[\si]\big)\\
\ns\ds\q +\rho-\dbE[\rho]\big\}\1n-\1n(\h\BP\h B\1n+\1n\h\BC^\top\BP\h D\1n+\1n\h\BS^\top\BI_m)\bar\BSi^{-1}\BJ^\top\big\{\h\BB^\top\bar\eta+\h\BD^\top\dbE[\z]
+\h\BD^\top\BP\dbE[\si]+\dbE[\rho]\big\}+\BP\big(b-\dbE[b]\big)\\
\ns\ds\q +\h\BP\dbE[b]\1n+\1n\dot{\bar\eta}\1n+\1n\BA^\top\big(\eta\1n-\1n\dbE[\eta]\big)\1n+\1n\h\BA^\top\bar\eta\1n
  +\1n\BC^\top\big(\BP(\si\1n-\1n\dbE[\si])\1n+\1n\z\1n-\1n\dbE[\z]\big)\1n+\1n\h\BC^\top\big(\dbE[\z]\1n+\1n\BP\dbE[\si]\big)\1n+\1n q,\q \as\ea\ee
Applying $\dbE[\,\cd\,]$ to the above, we obtain \rf{BODE-bar-eta-BJ}. Putting \rf{BODE-bar-eta-BJ} into \rf{alpha}, from \rf{BSDE-eta-game-form}, we obtain \rf{BSDE-eta-game-BJ}.

Moreover, we obtain the following closed-loop MF-SDE on $[0,\i)$ (recall \rf{A-Th}):
\bel{X*}\left\{\2n\ba{ll}
\ds dX^*(t)=\big\{A_{\Th^{**}}X^*+\bar A_{\BTh^{**}}\dbE[X^*]+Bv^{**}+\bar B\dbE[v^{**}]+b\big\}dt\\
\ns\ds\qq\qq\q+\big\{C_{\Th^{**}}X^*(t)+\bar C_{\BTh^{**}}\dbE[X^*]\1n+\1n Dv^{**}\2n+\1n\bar D\dbE[v^{**}]\1n+\1n\si\big\}dW(t),\q t\ges0,\\
\ns\ds X^*(0)=x,\ea\right.\ee
where $(\Th^{**},\bar\Th^{**},v^{**}(\cd))$ is given by \rf{Th-Th*-v*} and \rf{v*}.

The above procedure implies that if \rf{ARE-LQ-game-BP} admits a solution pair $(\BP,\h\BP)\in\dbR^{2n\times n}\times\dbR^{2n\times n}$ such that $\BTh^{**}\equiv(\Th^{**},\bar\Th^{**})$ defined by \rf{Th-Th*-v*} stabilizes the system $[A,\bar A,C,\bar C;B,\bar B,D,\bar D]$, then ODE \rf{BODE-bar-eta-BJ} admits a solution $\bar\eta(\cd)\in L^2(\dbR^{2n})$, BSDE \rf{BSDE-eta-game-BJ} admits a solution $(\eta(\cd),\zeta(\cd))\in\sX[0,\i)^2\times L^2_\dbF(\dbR^{2n})$ and the triple $(X^*(\cd),\BY^*(\cd),\BZ^*(\cd))\in\sX[0,\i)\times\sX[0,\i)^2\times L^2_\dbF(\dbR^{2n})$, defined through \rf{X*}, \rf{BY} and \rf{BZ}, is an adapted solution to MF-FBSDE \rf{MF-FBSDE}, with respect to the control $u^{**}(\cd)$ defined by \rf{u*}, and the stationarity condition \rf{stationarity condition-BJ} holds. Hence, if, in addition, the convexity condition \rf{convexity condition} holds for $i=1,2$, then by Theorem \ref{T-3.3}, Problem (MF-SDG) admits an open-loop Nash equilibrium for any initial state $x\in\dbR^n$, which has the closed-loop representation \rf{closed-loop ui**}. The proof is complete.
\endpf

\ms

Note that by the definition of $\Th^{**}$ and $\bar\Th^{**}$ in \rf{Th-Th*-v*}, we have
\bel{Th-Th*}
\BSi\Th^{**}+\BJ^\top\big(\BB^\top\BP+\BD^\top\BP C+\BS\BI_n\big)=0,\qq
\bar\BSi\bar\Th^{**}+\BJ^\top\big(\h\BB^\top\h\BP+\h\BD^\top\BP\h C+\h\BS\BI_n\big)=0;
\ee
and we can rewrite \rf{ARE-LQ-game-BP} as
\bel{AE-LQ-game-BP}\left\{\2n\ba{ll}
\ds\BP A+\BA^\top\BP+\BC^\top\BP C+\BQ\BI_n+\big(\BP B+\BC^\top\BP D+\BS^\top\BI_m\big)\Th^{**}=0,\\
\ns\ds\h\BP\h A+\h\BA^\top\h\BP+\h\BC^\top\BP\h C+\h\BQ\BI_n+\big(\h\BP\h B+\h\BC^\top\BP\h D+\h\BS^\top\BI_m\big)\bar\Th^{**}=0.
\ea\right.\ee
We may further write \rf{Th-Th*} and \rf{AE-LQ-game-BP} in the component forms:
\bel{AE-LQ-game-component}\left\{\2n\ba{ll}
\ds P_i A+A^\top P_i+C^\top P_i C+Q_i+\big(P_i B+C^\top P_iD+S_i^\top\big)\Th^{**}=0,\\
\ns\ds\h P_i\h A+\h A^{\,\top}\h P_i+\h C^\top P_i\h C+\h Q_i+\big(\h P_i\h B+\h C^\top P_i\h D+\h S_i^{\,\top}\big)\bar\Th^{**}=0,
\ea\right.\qq i=1,2,\ee
\bel{Th-Th*-component}\left\{\2n\ba{ll}
\ds\begin{pmatrix}R_{111}+D_1^\top P_1D_1&R_{112}+D_1^\top P_1D_2\\R_{221}+D_2^\top P_2D_1&R_{222}+D_2^\top P_2D_2\end{pmatrix}\Th^{**}
 +\begin{pmatrix}B_1^\top P_1+D_1^\top P_1 C+S_{11}\\B_2^\top P_2+D_2^\top P_2 C+S_{22}\end{pmatrix}=0,\\ [4mm]
\ns\ds\begin{pmatrix}\h R_{111}+\h D_1^\top P_1\h D_1&\h R_{112}+\h D_1^\top P_1\h D_2\\
 \h R_{221}+\h D_2^\top P_2\h D_1&\h R_{222}+\h D_2^\top P_2\h D_2\end{pmatrix}\bar\Th^{**}
 +\begin{pmatrix}\h B_1^\top\h P_1+\h D_1^\top P_1\h C+\h S_{11}\\ \h B_2^\top\h P_2+\h D_2^\top P_2\h C+\h S_{22}\end{pmatrix}=0.
\ea\right.\ee
In the above, the coefficient matrices of the equations for $\Th^{**}$ and $\bar\Th^{**}$ are not symmetric in general (even if $P_1,P_2,\h P_1,\h P_2$ are all symmetric). Hence, the equations for $P_i,\h P_i,i=1,2$ and for $P_i^\top,\h P_i^\top, i=1,2$ are different. Consequently, we do not expect $P_i,\h P_i,i=1,2$ to be symmetric in general.

\ms

\subsection{Closed-loop Nash equalibria and symmetric algebraic Riccati equations}
We now look at closed-loop Nash equilibria for Problem (MF-SDG). First, we present the following result, which is a consequence of Theorem \ref{T-3.3}.
\bp{P-3.5} \sl Let {\rm(H2)} hold. If $(\BTh^*,v^*(\cd))$ is a closed-loop Nash equilibrium of Problem {\rm(MF-SDG)}, then $(\BTh^*,0)$ is a closed-loop Nash equilibrium of Problem {\rm(MF-SDG)$^0$}.
\ep

\proof By the observation we made at the end of subsection 3.1, we see that $(\BTh^*,v^*(\cd))$ is a closed-loop Nash equilibrium of Problem (MF-SDG) if and only if $v^*(\cd)$ is an open-loop Nash equilibrium of the problem for any initial state $x\in\dbR^n$, with the stabilized state equation \rf{closed-loop-system-game-1} and cost functionals \rf{cost-LQ-game-closed} which we rewrite here in details for convenience (we suppress some $t$):
\bel{cost-LQ-game-closed-1}\ba{ll}
\ds J^{\BTh^*}_i(x;v(\cd))\equiv J_i\big(x;\Th^*\big\{X(\cd)-\dbE[X(\cd)]\big\}+\bar\Th^*\dbE[X(\cd)]+v(\cd)\big)\\
\ns\ds=\1n\dbE\2n\int_0^\i\2n\Bigg[\Bigg\langle\2n
\begin{pmatrix}\cQ_i^*&\3n(\cS_i^*)^\top\\ \cS_i^*&\3n R_i\end{pmatrix}\2n
\begin{pmatrix}X\\ v\end{pmatrix}\1n,\1n
\begin{pmatrix}X\\ v\end{pmatrix}\2n\Bigg\rangle
\1n+\1n2\Bigg\langle\2n\begin{pmatrix}q^*_i(t)\\ \rho_i(t)\end{pmatrix}\1n,\1n\begin{pmatrix}X\\ v\end{pmatrix}\2n\Bigg\rangle
\1n+\1n\Bigg\langle\2n
\begin{pmatrix}\bar\cQ_i^*&\3n(\bar\cS_i^*)^\top\\ \bar\cS_i^*&\3n\bar R_i\\\end{pmatrix}
\begin{pmatrix}\dbE[X]\\ \dbE[v]\end{pmatrix},\begin{pmatrix}\dbE[X]\\ \dbE[v]\end{pmatrix}\2n\Bigg\rangle\Bigg]dt,\ea\ee
where (with $\h S_i=S_i+\bar S_i$ and $\h R_i=R_i+\bar R_i$, comparing with \rf{Q-Th})
$$\left\{\2n\ba{ll}
\ds\cQ^*_i=Q_i+S_i^\top\Th^*+(\Th^*)^\top S_i+(\Th^*)^\top R_i\Th^*,\\
\ns\ds\bar\cQ^*_i=\bar Q_i+\h S_i^{\,\top}\bar\Th^*+(\bar\Th^*)^\top\h S_i+(\bar\Th^*)^\top\h R_i\bar\Th^*-S_i\Th^*-(\Th^*)^\top S_i-(\Th^*)^\top R_i\Th^*,\\
\ns\ds\cS^*_i=S_i+R_i\Th^*,\q\bar\cS^*_i=\bar S_i+\h R_i\bar\Th^*-R_i\Th^*,\q q_i^*(\cd)=q_i(\cd)+(\Th^*)^\top\big(\rho_i(\cd)-\dbE[\rho_i(\cd)]\big)
+(\bar\Th^*)^\top\dbE[\rho_i(\cd)].\ea\right.$$
Thus by Theorem \ref{T-3.3}, $(\BTh^*,v^*(\cd))$ is a closed-loop Nash equilibrium of Problem (MF-SDG) if and only if for any $x\in\dbR^n$, the solution $(X^*(\cd),Y_i^*(\cd),Z_i^*(\cd))\in
\sX[0,\i)\times\sX[0,\i)\times L^2_\dbF(\dbR^n)$ to the following MF-FBSDE:
\bel{MF-FBSDE-closed*}\left\{\2n\ba{ll}
\ds dX^*(t)=\big\{A_{\Th^*}X^*+\bar A_{\BTh^*}\dbE[X^*]+Bv^*+\bar B\dbE[v^*]+b\big\}dt\\
\ns\ds\qq\qq\qq+\big\{C_{\Th^*}X^*+\bar C_{\BTh^*}\dbE[X^*]+Dv^*+\bar D\dbE[v^*]+\si\big\}dW(t),\qq t\ges0,\\
\ns\ds-dY_i^*(t)=\big\{A_{\Th^*}^\top Y_i^*+\bar A_{\Th^*}^\top\dbE[Y_i^*]+C_{\Th^*}^\top Z_i^*+\bar C_{\Th^*}^\top\dbE[Z_i^*]+\cQ^*_iX^*+\bar\cQ^*_i\dbE[X^*]\\
\ns\ds\qq\qq\qq+(\cS^*_i)^\top v^*+(\bar\cS^*_i)^\top\dbE[v^*]+q^*_i+(\bar\Th^*-\Th^*)^\top\dbE[\rho_i]\big\}dt-Z_i^*dW(t),\qq t\ges0,\\
\ns\ds X^*(0)=x,\ea\right.\ee
for $i=1,2$, satisfies the following stationarity condition:
\bel{stationarity condition-closed*}\ba{ll}
\ds B_i^\top Y_i^*+\bar B_i^\top\dbE[Y_i^*]+D_i^\top Z_i^*+\bar D_i^\top\dbE[Z_i^*]+\cS^*_{ii}X^*+\bar\cS^*_{ii}\dbE[X^*]+R_{ii}v^*+\bar R_{ii}\dbE[v^*]+\rho_{ii}=0,\\
\ns\ds\qq\qq\qq\qq\qq\qq\qq\qq\qq\qq\qq\qq\qq\ae \ t\in[0,\i),~\as,~i=1,2,\ea\ee
where
$$\cS^*_{ii}=S_{ii}+R_{ii}\Th^*,\q\bar\cS^*_{ii}=\bar S_{ii}+\h R_{ii}\bar\Th^*-R_{ii}\Th^*,$$
and for $i=1,2$, the following convexity condition holds:
\bel{convexity condition-closed}\ba{ll}
\ds\dbE\int_0^\i\Big[\blan\cQ^*_iX_i,X_i\bran+2\blan\cS^*_{ii}X_i,v_i\bran+\blan R_{iii}v_i,v_i\bran+\blan\bar\cQ^*_i\dbE[X_i],\dbE[X_i]\bran\\
\ns\ds\qq\qq+2\blan\bar\cS^*_{ii}\dbE[X_i^v],\dbE[v_i]\bran+\blan\bar R_{iii}\dbE[v_i],\dbE[v_i]\bran\Big]dt\ges0,\qq \forall v(\cd)\in L^2_\dbF(\dbR^m),\ea\ee
where $X_i(\cd)\in\sX[0,\i)$ is the solution to the following controlled homogeneous MF-SDE on $[0,\i)$:
\bel{controlled homogeneous MF-SDE-closed}\left\{\2n\ba{ll}
\ds dX_i(t)=\big\{A_{\Th^*}X_i+\bar A_{\BTh^*}\dbE[X_i]+B_iv_i+\bar B_i\dbE[v_i]\big\}dt\\
\ns\ds\qq\qq\qq+\big\{C_{\Th^*}X_i+\bar C_{\BTh^*}\dbE[X_i]+D_iv_i+\bar D_i\dbE[v_i]\big\}dW(t),\q t\ges0,\\
\ns\ds X_i(0)=0.\ea\right.\ee
Since $(\BTh^*,v^*(\cd))$ is independent of $x$ and \rf{MF-FBSDE-closed}--\rf{controlled homogeneous MF-SDE-closed} hold for all $x\in\dbR^n$, by subtracting equations corresponding to $x$ and 0, the latter from the former, we see that for any $x\in\dbR^n$, the following MF-FBSDE:
$$\left\{\2n\ba{ll}
\ds dX_0^*(t)=\big\{A_{\Th^*}X_0^*+\bar A_{\BTh^*}\dbE[X_0^*]\big\}dt+\big\{C_{\Th^*}X_0^*
+\bar C_{\BTh^*}\dbE[X_0^*]\big\}dW(t),\q t\ges0,\\
\ns\ds-dY_{i0}^*(t)=\big\{A_{\Th^*}^\top Y_{i0}^*+\bar A_{\BTh^*}^\top\dbE[Y_{i0}^*]
                       +C_{\Th^*}^\top Z_{i0}^*+\bar C_{\BTh^*}^\top\dbE[Z_{i0}^*]
                       +\cQ^*_iX_0^*+\bar\cQ^*_i\dbE[X_0^*]\big\}dt-Z_{i0}(t)dW(t),\\
\ns\ds\qq\qq\qq\qq\qq\qq\qq\qq\qq\qq\qq\qq\qq t\ges0,\q i=1,2,\\
\ns\ds X_0^*(0)=x,\ea\right.$$
admits an adapted solution satisfying
$$\ba{ll}
 B_i^\top Y_{i0}^*+\bar B_i^\top\dbE[Y_{i0}^*]+D_i^\top Z_{i0}^*+\bar D_i^\top\dbE[Z_{i0}^*]+\cS^*_{ii}X_0^*+\bar\cS^*_{ii}\dbE[X_0^*]=0,\q\ae \ t\in[0,\i),~\as,\ i=1,2.\ea$$
It follows, again from Theorem \ref{T-3.3}, that $(\BTh^*,0)$ is a closed-loop Nash equilibrium of Problem (MF-SDG)$^0$. The proof is complete.
\endpf

\ms

Now, we give a necessary condition for the existence of closed-loop Nash equilibria of Problem (MF-SDG).
\bp{P-3.6} \sl Let {\rm(H2)} hold, and let $(\BTh^*,v^*(\cd))$ be a closed-loop Nash equilibrium of Problem (MF-SDG). Then for $i=1,2$, the following system of coupled AREs admits a solution pair $(P_i,\h P_i)\in\dbS^n\times\dbS^n$:
\bel{AE-LQ-game-closed}\left\{\2n\ba{ll}
\ds P_iA+A^\top P_i+C^\top P_iC+Q_i+(\Th^*)^\top(R_i+D^\top P_iD)\Th^*\\
\ns\ds\qq+(P_iB+C^\top P_iD+S_i^\top)\Th^*+(\Th^*)^\top(B^\top P_i+D^\top P_iC+S_i)=0,\\
\ns\ds\h P_i\h A+\h A^\top\h P_i+\h C^\top P_i\h C+\h Q_i+(\bar\Th^*)^\top(\h R_i+\h D^\top P_i\h D)\bar\Th^*\\
\ns\ds\qq+(\h P_i\h B+\h C^\top P_i\h D+\h S_i^\top)\bar\Th^*+(\bar \Th^*)^\top(\h B^\top\h P_i+\h D^\top P_i\h C+\h S_i)=0,\ea\right.\ee
and the following conditions are satisfied:
\bel{Th-Th*-closed}\left\{\2n\ba{ll}
\ds B_i^\top P_i+D_i^\top P_iC+S_{ii}+(R_{ii}+D_i^\top P_iD)\Th^*=0,\\
\ns\ds\h B_i^\top\h P_i+\h D_i^\top P_i\h C+\h S_{ii}+(\h R_{ii}+\h D_i^\top P_i\h D)\bar\Th^*=0,\ea\right.\ee
and
\bel{Sigma}\Si_i\deq R_{iii}+D_i^\top P_iD_i\ges0,\qq \bar\Si_i\deq\h R_{iii}+\h D_i^\top P_i\h D_i\ges0.\ee
\ep
\proof Suppose that $(\BTh^*,v^*(\cd))$ is a closed-loop Nash equilibrium of Problem (MF-SDG). Then by Proposition \ref{P-3.5}, $(\BTh^*,0)$ is a closed-loop Nash equilibrium of Problem (MF-SDG)$^0$. Denote
\bel{cA-1}\left\{\2n\ba{ll}
\ds\cA_1\deq A+B_2\Th^*_2,\qq \bar\cA_1\deq\bar A+\bar B_2\bar \Th^*_2+B_2(\bar \Th^*_2-\Th^*_2),\\
\ns\ds\cC_1\deq C+D_2\Th^*_2,\qq \bar\cC_1\deq\bar C+\bar D_2\bar \Th^*_2+D_2(\bar \Th^*_2-\Th^*_2),\\
\ns\ds\cQ_1\deq Q_1+S_{12}^\top\Th_2^*+(\Th^*_2)^\top S_{12}+(\Th^*_2)^\top R_{122}\Th^*_2,\\
\ns\ds\bar\cQ_1\deq\bar Q_1+\h S_{12}^\top\bar\Th^*_2+(\bar\Th^*_2)^\top\h S_{12}+(\bar\Th^*_2)^\top\h R_{122}\bar\Th^*_2-S_{12}^\top\Th^*_2-(\Th^*_2)^\top S_{12}-(\Th^*_2)^\top R_{122}\Th^*_2,\\
\ns\ds\cS_{11}\deq S_{11}+R_{112}\Th^*_2,\qq\bar\cS_{11}\deq\bar S_{11}+\bar R_{112}\bar\Th^*_2+R_{112}(\bar\Th^*_2-\Th^*_2).\ea\right.\ee
Then, for any $u_1(\cd)\in L^2_\dbF(\dbR^{m_1})$, let us consider the state equation:
$$\left\{\2n\ba{ll}
\ds dX_1^0(t)=\big\{\cA_1X_1^0(t)+\bar\cA_1\dbE[X_1^0(t)]+B_1u_1(t)+\bar B_1\dbE[u_1(t)]\big\}dt\\
\ns\ds\qq\qq\qq+\big\{\cC_1X_1^0(t)+\bar\cC_1\dbE[X_1^0(t)]+D_1u_1(t)+\bar D_1\dbE[u_1(t)]\big\}dW(t),\q t\ges0,\\
\ns\ds X_1^0(0)=0,\ea\right.$$
and cost functional:
$$\ba{ll}
\ds\bar J_1^0(x;u_1(\cd))\equiv J_1^0\big(x;u_1(\cd);\Th^*_2\big\{X_1^0(\cd)-\dbE[X_1^0(\cd)]\big\}+\bar\Th^*_2\dbE[X_1^0(\cd)]\big)\\
\ns\ds=\dbE\int_0^\i\Bigg[\Bigg\langle\2n
\begin{pmatrix}\cQ_1&\cS_{11}^\top\\ \cS_{11}& R_{111}\end{pmatrix}
\begin{pmatrix}X_1^0\\ u_1\end{pmatrix},
\begin{pmatrix}X_1^0\\ u_1\end{pmatrix}\2n\Bigg\rangle+\Bigg\langle\2n
\begin{pmatrix}\bar\cQ_1&\bar\cS_{11}^\top\\ \bar\cS_{11}&\bar R_{111}\\\end{pmatrix}
\begin{pmatrix}\dbE[X_1^0]\\ \dbE[u_1]\end{pmatrix},\begin{pmatrix}\dbE[X_1^0]\\ \dbE[u_1]\end{pmatrix}\2n\Bigg\rangle\Bigg]dt.\ea$$
It is easy to see that $(\BTh^*_1,0)$ is a closed-loop optimal strategy for the above mean-field LQ stochastic optimal control problem. Thanks to Theorem \ref{T-2.9}, the following system of coupled AREs:
\bel{ARE-1}
\left\{\2n\ba{ll}
\ds P_1\cA_1+\cA_1^\top P_1+\cC_1^\top P_1\cC_1+\cQ_1-\big(P_1B_1+\cC_1^\top P_1D_1+\cS_{11}^\top\big)\Si_1^\dag\big(B_1^\top P_1+D_1^\top P_1\cC_1+\cS_{11}\big)=0,\\
\ns\ds\h P_1\h\cA_1+\h\cA_1^{\,\top}\h P_1+\h\cC_1^{\,\top}P_1\h\cC_1+\h\cQ_1-\big(\h P_1\h B_1+\h\cC_1^{\,\top}P_1\h D_1+\h\cS_{11}^{\,\top}\big)\bar\Si_1^\dag\big(\h B_1^\top\h P_1+\h D_1^\top P_1\h\cC_1+\h\cS_{11}\big)=0,\ea\right.\ee
admit a static stabilizing solution pair $(P_1,\h P_1)\in\dbS^n\times\dbS^n$, satisfying
\bel{Th1-Th*1-closed}\left\{\2n\ba{ll}
\ds B_1^\top P_1+D_1^\top P_1\cC_1+\cS_{11}+\Si_1\Th_1^*=0, \qq \Si_1\equiv R_{111}+D_1^\top P_1D_1\ges0,\\
\ns\ds\h B_1^\top\h P_1+\h D_1^\top P_1\h\cC_1+\h\cS_{11}+\bar\Si_1\bar\Th_1^*=0, \qq \bar\Si_1\equiv\h R_{111}+\h D_1^\top P_1\h D_1\ges0,\ea\right.\ee
where $\h\cA_1=\cA+\bar\cA_1$, $\h\cC_1=\cC_1+\bar\cC_1$, $\h\cQ_1=\cQ_1+\bar\cQ_1$, and
$\h\cS_{11}=\cS_{11}+\bar\cS_{11}$. Similarly, for any $u_2(\cd)\in L^2_\dbF(\dbR^{m_2})$, we can consider the state equation:
$$\left\{\2n\ba{ll}
\ds dX_2^0(t)=\big\{\cA_2X_2^0(t)+\bar\cA_2\dbE[X_2^0(t)]+B_2u_2(t)+\bar B_2\dbE[u_2(t)]\big\}dt\\
\ns\ds\qq\qq\q+\big\{\cC_2X_2^0(t)+\bar\cC_2\dbE[X_2^0(t)]+D_2u_2(t)+\bar D_2\dbE[u_2(t)]\big\}dW(t),\q t\ges0,\\
\ns\ds X_2^0(0)=0,\ea\right.$$
and cost functional:
$$\ba{ll}
\ds\bar J_2^0(x;u_2(\cd))\equiv J_2^0\big(x;\Th^*_1\big(X_2^0(\cd)-\dbE[X_2^0(\cd)]\big)+\bar\Th^*_1\dbE[X_2^0(\cd)];u_2(\cd)\big)\\
\ns\ds=\dbE\int_0^\i\Bigg[\Bigg\langle\2n
\begin{pmatrix}\cQ_2&\cS_{22}^\top\\ \cS_{22}&R_{222}\end{pmatrix}
\begin{pmatrix}X_2^0\\ u_2\end{pmatrix},
\begin{pmatrix}X_2^0\\ u_2\end{pmatrix}\2n\Bigg\rangle+\Bigg\langle\2n
\begin{pmatrix}\bar\cQ_2&\bar\cS_{22}^\top\\ \bar\cS_{22}&\bar R_{222}\\\end{pmatrix}
\begin{pmatrix}\dbE[X_2^0]\\ \dbE[u_2]\end{pmatrix},\begin{pmatrix}\dbE[X_2^0]\\ \dbE[u_2]\end{pmatrix}\2n\Bigg\rangle\Bigg]dt,\ea$$
where $\cA_2,\bar\cA_2,\cC_2,\bar\cC_2$ and $\cQ_2,\bar\cQ_2,\cS_{22},\bar\cS_{22}$ are defined similar to \rf{cA-1}. In the same spirit as above, we can see that $(\BTh^*_2,0)\equiv(\Th^*_2,\bar\Th^*_2,0)$ is a closed-loop optimal strategy for the above mean-field LQ stochastic optimal control problem. Making use of Theorem \ref{T-2.9} again, the following system of coupled AREs
\bel{ARE-2}\left\{\2n\ba{ll}
\ds P_2\cA_2+\cA_2^\top P_2+\cC_2^\top P_2\cC_2+\cQ_2-\big(P_2B_2+\cC_2^\top P_2D_2+\cS_{22}^\top\big)\Si_2^\dag\big(B_2^\top P_2+D_2^\top P_2\cC_2+\cS_{22}\big)=0,\\
\ns\ds\h P_2\h\cA_2+\h\cA_2^{\,\top}\h P_2+\h\cC_2^{\,\top}P_2\h\cC_2+\h\cQ_2-\big(\h P_2\h B_2+\h\cC_2^{\,\top}P_2\h D_2+\h\cS_{22}^{\,\top}\big)\bar\Si_2^\dag
\big(\h B_2^\top\h P_2+\h D_2^\top P_2\h\cC_2+\h\cS_{22}\big)=0,\ea\right.\ee
admit a static stabilizing solution pair $(P_2,\h P_2)\in\dbS^n\times\dbS^n$, satisfying
\bel{Th2-Th*2-closed}\left\{\2n\ba{ll}
\ds0=B_2^\top P_2+D_2^\top P_2\cC_2+\cS_{22}+\Si_2\Th_2^*, \qq\Si_2\equiv R_{222}+D_2^\top P_2D_2\ges0,\\
\ns\ds0=\h B_2^\top\h P_2+\h D_2^\top P_2\h\cC_2+\h\cS_{22}+\bar\Si_2\bar\Th_2^*, \qq \bar\Si_2\equiv\h R_{222}+\h D_2^\top P_2\h D_2\ges0,\ea\right.\ee
where $\h\cA_2=\cA_2+\bar\cA_2$, $\h\cC_2=\cC_2+\bar\cC_2$, $\h\cQ_2=\cQ_2+\bar\cQ_2$, and
$\h\cS_{22}=\cS_{22}+\bar\cS_{22}$. By \rf{Th1-Th*1-closed} and \rf{Th2-Th*2-closed} and putting $\BTh^*\equiv(\Th^*,\bar\Th^*)=\bigg(\begin{pmatrix}\Th_1^*\\\Th_2^*\end{pmatrix},\begin{pmatrix}\bar\Th_1^*\\\bar\Th_2^*\end{pmatrix}\bigg)\in\dbR^{2m\times 2n}$, we get \rf{Sigma} and
$$\left\{\2n\ba{ll}
\ds0=B_1^\top P_1+D_1^\top P_1\cC_1+\cS_{11}+\Si_1\Th_1^*\\
\ns\ds\q=B_1^\top P_1+D_1^\top P_1C+S_{11}+R_{111}\Th_1^*+R_{112}\Th^*_2+D_1^\top P_1(D_1\Th_1^*+D_2\Th^*_2)\\
\ns\ds\q=B_1^\top P_1+D_1^\top P_1C+S_{11}+(R_{11}+D_1^\top P_1D)\Th^*,\\
\ns\ds0=\h B_1^\top\h P_1+\h D_1^\top P_1\h\cC_1+\h\cS_{11}+\bar\Si_1\bar\Th_1^*\\
\ns\ds\q=\h B_1^\top\Pi_1+\h D_1^\top P_1\h C+\h S_{11}+\h R_{111}\Th_1^*+\h R_{112}
\bar\Th^*_2+\h D_1^\top P_1\h D_1\Th_1^*+\h D_1^\top P_1\h D_2\bar\Th^*_2\\
\ns\ds\q=\h B_1^\top\h P_1+\h D_1^\top P_1\h C+\h S_{11}+\big(\h R_{11}+\h D_1^\top P_1\h D\big)\bar\Th^*.\ea\right.$$
Similarly,
$$\left\{\2n\ba{ll}
\ds0=B_2^\top P_2+D_2^\top P_2\cC_2+\cS_{22}+\Si_2\Th_2^*
  =B_2^\top P_2+D_2^\top P_2C+S_{22}+(R_{22}+D_2^\top P_2D)\Th^*,\\
\ns\ds0=\h B_2^\top\h P_2+\h D_2^\top P_2\h\cC_2+\h\cS_{22}+\bar\Si_2\bar\Th_2^*=\h B_2^\top\h P_2+\h D_2^\top P_2\h C+\h S_{22}+\big(\h R_{22}+\h D_2^\top P_2\h D\big)\bar\Th^*,\ea\right.$$
which implies \rf{Th-Th*-closed}. By \rf{Th-Th*-closed}, \rf{ARE-1} and \rf{ARE-2}, we have (by a straightforward calculation)
$$\left\{\ba{ll}
\ns\ds0=P_1\cA_1+\cA_1^\top P_1+\cC_1^\top P_1\cC_1+\cQ_1-\big(P_1B_1+\cC_1^\top P_1D_1+\cS_{11}^\top\big)\Si_1^\dag\big(B_1^\top P_1+D_1^\top P_1\cC_1+\cS_{11}\big)\\
\ns\ds\q=P_1A+A^\top P_1+C^\top P_1C+Q_1+(\Th^*)^\top(R_1+D^\top P_1D)\Th^*\\
\ns\ds\qq+\big(P_1B+C^\top P_1D+S_1^\top\big)\Th^*+(\Th^*)^\top\big(B^\top P_1+D^\top P_1C+S_1\big),\\
\ns\ds0=\h P_1\h\cA_1+\h\cA_1^{\,\top}\h P_1+\h\cC_1^{\,\top}P_1\h\cC_1
  +\h\cQ_1-\big(\h P_1\h B_1+\h\cC_1^{\,\top}P_1\h D_1+\h\cS_{11}^{\,\top}\big)\bar\Si_1^\dag\big(\h B_1^\top\h P_1
  +\h D_1^\top P_1\h\cC_1+\h\cS_{11}\big)\\
\ns\ds\q=\h P_1\h A+\h A^{\,\top}\h P_1+\h C^\top P_1\h C+\h Q_1+(\bar\Th^*)^\top\big(\h R_1
+\h D^\top P_1\h D\big)\bar\Th^*\\
\ns\ds\qq+\big(\h P_1\h B+\h C^\top P_1\h D+\h S_1^\top\big)\bar\Th^*+(\bar\Th^*)^\top\big(\h B^\top\h P_1+\h D^\top P_1\h C+\h S_1\big).\ea\right.$$
In the same way, we have
$$\left\{\2n\ba{ll}
\ds0=P_2\cA_2+\cA_2^\top P_2+\cC_2^\top P_2\cC_2+\cQ_2-\big(P_2B_2+\cC_2^\top P_2D_2+\cS_{22}^\top\big)\Si_2^\dag
  \big(B_2^\top P_2+D_2^\top P_2\cC_2+\cS_{22}\big)\\
\ns\ds\q=P_2A+A^\top P_2+C^\top P_2C+Q_2+(\Th^*)^\top(R_2+D^\top P_2D)\Th^*\\
\ns\ds\qq+\big(P_2B+C^\top P_2D+S_2^\top\big)\Th^*+(\Th^*)^\top\big(B^\top P_2+D^\top P_2C+S_2\big),\\
\ns\ds0=\h P_2\h\cA_2+\h\cA_2^{\,\top}\h P_2+\h\cC_2^{\,\top}P_2\h\cC_2+\h\cQ_2-\big(\h P_2\h B_2+\h\cC_2^{\,\top}P_2\h D_2+\h\cS_{22}^{\,\top}\big)\bar\Si_2^\dag\big(\h B_2^\top\h P_2+\h D_2^\top P_2\h\cC_2+\h\cS_{22}\big)\\
\ns\ds\q=\h P_2\h A+\h A^{\,\top}\h P_2+\h C^\top P_2\h C+\h Q_2+(\bar\Th^*)^\top\big(\h R^2+\h D^\top P_2\h D\big)\bar\Th^*
\\
\ns\ds\qq+\big(\h P_2\h B+\h C^\top P_2\h D+\h S_2^\top\big)\bar\Th^*+(\bar\Th^*)^\top\big(\h B^\top\h P_2+\h D^\top P_1\h C+\h S_2\big).\ea\right.$$
This yields \rf{AE-LQ-game-closed}. The proof is complete.
\endpf

\ms
Now, we are ready to present the main result of this subsection, which characterizes the closed-loop Nash equilibrium of Problem (MF-SDG).
\bt{T-3.6} \sl
Let {\rm(H2)} hold. Problem (MF-SDG) admits a closed-loop Nash equilibrium $(\BTh^*,v^*(\cd))$ if and only if the following statements hold:

\ms

{\rm(i)}\ The system of coupled algebraic equations \rf{AE-LQ-game-closed} admits a solution pair $(P_i,\h P_i)\in\dbS^n\times\dbS^n$ and \rf{Th-Th*-closed}--\rf{Sigma} are satisfied.

\ms

{\rm(ii)}\ For $i=1,2$, the following BSDE:
\bel{BSDE-eta_i-game}\ba{ll}
d\eta_i(t)=-\big\{A_{\Th^*}^\top\eta_i(t)\1n+\1n C_{\Th^*}^\top\big[\z_i(t)\1n+\1n P_i\si(t)\big]\1n+\1n P_ib(t)+\1n q_i(t)+\1n(\Th^*)^\top\rho_i(t)\1n\big\}dt\1n+\1n\z_i(t)dW(t),\q t\ges0,
\ea\ee
admits a solution $(\eta_i(\cd),\z_i(\cd))\in\sX[0,\i)\times L^2_\dbF(\dbR^n)$ such that
\bel{v*-Ev*}\ba{ll}
\ds0=(R_{ii}+D_i^\top P_iD)\big(v^*(t)-\dbE[v^*(t)]\big)
+B_i^\top\big(\eta_i(t)-\dbE[\eta_i(t)]\big)+D_i^\top\big(\z_i(t)-\dbE[\z_i(t)]\big)\\
\ns\ds\qq\q+D_i^\top P_i\big(\si(t)-\dbE[\si(t)]\big)+\rho_{ii}(t)-\dbE[\rho_{ii}(t)],\qq\ae \ t\in[0,\i),~\as,\qq i=1,2,\ea\ee
and the following ODE admits a solution $\bar\eta_i(\cd)\in L^2(\dbR^n)$:
\bel{BODE-bar-eta_i-game}\ba{ll}
\ds\dot{\bar\eta}_i(t)+\h A_{\bar\Th^*}^{\,\top}\bar\eta_i(t)+\h C_{\bar \Th^*}^{\,\top}\big(\dbE[\z_i(t)]+P_i\dbE[\si(t)]\big)+\h P_i\dbE[b(t)]+\dbE[q_i(t)]+(\bar\Th^*)^\top\dbE[\rho_i(t)]=0,\q t\ges0,\ea\ee
satisfying
\bel{Ev*}
(\h R_{ii}+\h D_i^\top P_i\h D)\dbE[v^*(t)]+\h B_i^\top\bar\eta_i(t)
+\h D_i^\top\dbE[\z_i(t)]+\h D_i^\top P_i\dbE[\si(t)]+\dbE[\rho_{ii}(t)]=0,\q\ae\ t\in[0,\i).
\ee
\et

\proof {\it Necessity}. Suppose that $(\BTh^*,v^*(\cd))$ is a closed-loop Nash equilibrium of Problem (MF-SDG). Then by Proposition \ref{P-3.6}, (i) holds.

For (ii), we first note that system \rf{AE-LQ-game-closed} is equivalent to
\bel{AE-LQ-game-closed-1}\left\{\2n\ba{ll}
\ds P_iA_{\Th^*}+A_{\Th^*}^\top P_i+C_{\Th^*}^\top P_iC_{\Th^*}+Q_i
+S_i^\top\Th^*+(\Th^*)^\top S_i+(\Th^*)^\top R_i\Th^*=0,\\
\ns\ds\h P_i\h A_{\bar\Th^*}+\h A_{\bar\Th^*}^{\,\top}\h P_i+\h C_{\bar\Th^*}^{\,\top} P_i\h C_{\Th^*}+\h Q_i+\h S_i^\top\bar\Th^*+(\bar\Th^*)^\top\h S_i+(\bar\Th^*)^\top\h R_i\bar\Th^*=0.\ea\right.\ee
Let $(X^*(\cd),Y_i^*(\cd),Z_i^*(\cd))\in\sX[0,\i)\times\sX[0,\i)\times L^2_\dbF(\dbR^n)$ be the solution to MF-FBSDE
\bel{MF-FBSDE-closed}\left\{\2n\ba{ll}
\ds dX^*(t)=\big\{A_{\Th^*}X^*+\bar A_{\BTh^*}\dbE[X^*]+Bv^*+\bar B\dbE[v^*]+b\big\}dt\\
\ns\ds\qq\qq\qq+\big\{C_{\Th^*}X^*+\bar C_{\BTh^*}\dbE[X^*]+Dv^*+\bar D\dbE[v^*]+\si\big\}dW(t),\\
\ns\ds-dY_i^*(t)=\big\{A_{\Th^*}^\top Y_i^*+\bar A_{\Th^*}^\top\dbE[Y_i^*]+C_{\Th^*}^\top Z_i^*+\bar C_{\Th^*}^\top\dbE[Z_i^*]+\cQ^*_iX^*+\bar\cQ^*_i\dbE[X^*]\\
\ns\ds\qq\qq\qq+(\cS^*_i)^\top v^*+(\bar\cS^*_i)^\top\dbE[v^*]+q^*_i+(\bar\Th^*-\Th^*)^\top\dbE[\rho_i]\big\}dt-Z_i^*dW(t),\qq t\ges0,\\
\ns\ds X^*(0)=x,\ea\right.\ee
for $i=1,2$. Proceeding as in the proof of Proposition \ref{P-3.5}, we see that $(X^*(\cd),Y_i^*(\cd),Z_i^*(\cd))$ satisfies
\bel{stationarity condition-closed}\ba{ll}
\ds B_i^\top Y_i^*+\bar B_i^\top\dbE[Y_i^*]+D_i^\top Z_i^*+\bar D_i^\top\dbE[Z_i^*]+\cS^*_{ii}X^*+\bar\cS^*_{ii}\dbE[X^*]+R_{ii}v^*+\bar R_{ii}\dbE[v^*]+\rho_{ii}=0,\\
\ns\ds\qq\qq\qq\qq\qq\qq\qq\qq\qq\qq\qq\qq\qq\ae \ t\in[0,\i),~\as,
\ea\ee
where $\cS^*_{ii}\deq S_{ii}+R_{ii}\Th^*,\,\bar\cS^*_{ii}\deq \bar S_{ii}+\h R_{ii}\bar\Th^*-R_{ii}\Th^*,\,i=1,2$. Now, we define
\bel{definition}\left\{\2n\ba{ll}
\ds\a_i\deq Y_i^*-\dbE[Y_i^*]-P_i\big(X^*-\dbE[X^*]\big), \qq \bar\eta_i\deq\dbE[Y_i^*]-\h P_i\dbE[X^*],\\
\ns\ds\z_i\deq Z_i^*-P_iC_{\Th^*}(X^*-\dbE[X^*])-P_i\h C_{\bar\Th^*}\dbE[X^*]-P_iD(v^*-\dbE[v^*])-P_i\si-P_i\h D\dbE[v^*],\q i=1,2.\ea\right.\ee
We want to show that
\bel{alpha_i}
\a_i(\cd)=\eta_i(\cd)-\dbE[\eta_i(\cd)],\ee
with $\eta_i(\cd),\z_i(\cd),\bar{\eta}(\cd)$ satisfy \rf{BSDE-eta_i-game}--\rf{Ev*}, for $i=1,2$. For this target, applying It\^{o}'s formula to $\a_i(\cd)$ yields
\bel{Ito-2}\ba{ll}
\ds d\a_i(t)=-\big\{A_{\Th^*}^\top \big(Y_i^*\2n-\1n\dbE[Y_i^*]\big)\1n+\1n C_{\Th^*}^\top\1n
\big(Z_i^*\2n-\1n\dbE[Z_i^*]\big)\1n+\1n\cQ_i\big(X^*\2n-\1n\dbE[X^*]\big)\1n+\1n\cS_i
\big(v^*\2n-\1n\dbE[v^*]\big)\1n+\1n q_i^*-\dbE[q_i^*]\big\}dt\\
\ns\ds\qq\qq+Z_idW(t)\1n-\1n P_i\big\{\cA_{\Th^*}^\top\big(X^*\2n-\1n\dbE[X^*]\big)\1n+\1nB\big(v^*\2n-\1n\dbE[v^*]\big)+b-\dbE[b]\big\}dt-P_i\big\{C_{\Th^*}^\top\big(X^*
\2n-\1n\dbE[X^*]\big)\\
\ns\ds\qq\qq+\big(\h C+\h D\bar\Th^*\big)\dbE[\mathrm{X}^*]+D(v^*-\dbE[v^*])+\h D\dbE[v^*]+\si\big\}dW(t)\\
\ns\ds\qq\q=-\big\{A_{\Th^*}^\top\a_i+C_{\Th^*}^\top\big(\z_i-\dbE[\z_i]\big)
+(P_iA_{\Th^*}+A_{\Th^*}^\top P_i+C_{\Th^*}^\top P_iC_{\Th^*}+\cQ_i)\big(X^*-\dbE[X^*]\big)\\
\ns\ds\qq\qq+\big[P_iB+C^\top P_iD+S_i^\top+(\Th^*)^\top(R_i+D^\top P_iD)\big]\big(v^*-\dbE[v^*]\big)+\cC_{\Th^*}^\top P_i(\si-\dbE[\si])\\
\ns\ds\qq\qq+P_i\big(b-\dbE[b]\big)+q_i-\dbE[q_i]+(\Th^*)^\top\big(\rho_i-\dbE[\rho_i]\big)\big\}dt+\z_idW(t)\\
\ns\ds\qq\q=-\big\{A_{\Th^*}^\top\a_i+C_{\Th^*}^\top\big(\z_i-\dbE[\z_i]\big)
+C_{\Th^*}^\top P_i\big(\si-\dbE[\si]\big)\\
\ns\ds\qq\qq+(\Th^*)^\top\big(\rho_i-\dbE[\rho_i]\big)+P_i(b-\dbE[b])+q_i-\dbE[q_i]\big\}dt+\z_idW(t).
\ea\ee
Since the solution $(\eta_i(\cd),\z_i(\cd))$ to \rf{BSDE-eta_i-game} satisfies
$$\ba{ll}
\ds d\big(\eta_i(t)-\dbE[\eta_i(t)]\big)=-\big\{A_{\Th^*}^\top\big(\eta_i-\dbE[\eta_i]\big)
+C_{\Th^*}^\top\big(\z_i-\dbE[\z_i]\big)+C_{\Th^*}^\top P_i\big(\si-\dbE[\si]\big)\\
\ns\ds\qq\qq\qq\qq\qq+(\Th^*)^\top(\rho_i-\dbE[\rho_i])+P_i(b-\dbE[b])+q_i-\dbE[q_i]\big\}dt+\z_idW(t),\q t\ges0,\ea$$
by the uniqueness of solutions, we obtain \rf{alpha_i}.
Moreover, we get
\bel{Ito-3}\ba{ll}
\ds-\dot{\bar\eta}_i(t)=\big[\h A+\h B\bar\Th^*\big]^\top\dbE[Y_i^*]+\big[\h C+\h D\bar\Th^*\big]^\top\dbE[Z_i^*]+\h\cQ_i\dbE[X^*]\\
\ns\ds\q+\h\cS_i\dbE[v^*]+\dbE[q_i^*]+\h P_i\big[\h A+\h B\bar\Th^*\big]\dbE[X^*]
 +\h P_i\h B\dbE[v^*]+\h P_i\dbE[b]\\
\ns\ds=\big[\h A+\h B\bar\Th^*\big]^\top\bar{\eta}_i+\big\{\h P_i\big(\h A+\h B\bar\Th^*\big)
 +\big(\h A+\h B\bar\Th^*\big)^\top\h P_i+\big[\h C+\h D\Th^*\big]^\top P_i\big[\h C+\h D\Th^*\big]+\h Q_i+\h S_i^\top\bar\Th^*\\
\ns\ds\q+(\bar\Th^*)^\top\h S_i+(\bar\Th^*)^\top\h R_i\bar\Th^*\big\}\dbE[X^*]+\big\{\h P_i\h B+\h C^\top P_i\h D+\h S_i^\top+(\bar\Th^*)^\top\big[\h R_i+\h D^\top P_i\h D\big]\big\}\dbE[v^*]\\
\ns\ds\q+(\bar\Th^*)^\top\big[\h D^\top(P_i\dbE[\si]+\dbE[\z_i])+\dbE[\rho_i]\big]+\h C^\top(P_i\dbE[\si]+\dbE[\z_i])+\h P_i\dbE[b]+\dbE[q_i]\\
\ns\ds=\1n\big(\h A\1n+\1n\h B\bar\Th^*\big)^\top\bar{\eta}_i\1n+\1n(\bar \Th^*)^\top\big[\h D^\top(P_i\dbE[\si]+\dbE[\z_i])+\dbE[\rho_i]\big]+\h C^\top(P_i\dbE[\si]+\dbE[\z_i])+\h P_i\dbE[b]
+\dbE[q_i],\ea\ee
which is \rf{BODE-bar-eta_i-game}. Further, from \rf{stationarity condition-closed} we have
\bel{stationarity condition-closed-1}
\h R_{ii}\dbE[v^*]+\h B_i^\top\dbE[Y_i^*]+\h D_i^\top\dbE[Z_i^*]+\big[\h S_{ii}+\h R_{ii}\bar\Th^*\big]\dbE[X^*]+\dbE[\rho_{ii}]=0,\qq\ae t\in[0,\i),\ i=1,2,\ee
and
\bel{stationarity condition-closed-2}\ba{ll}
\ds R_{ii}(v^*\1n-\1n\dbE[v^*])\1n+\1n B_i^\top(Y_i^*\1n-\1n\dbE[Y_i^*])\1n+\1n D_i^\top(Z_i^*\1n-\1n\dbE[Z_i^*])\1n+\1n(S_{ii}\1n+\1n R_{ii}
\Th^*)(X^*\1n-\1n\dbE[X^*])\1n+\1n\rho_{ii}\1n-\1n\dbE[\rho_{ii}]\1n=\1n0,\\
\ns\ds\qq\qq\qq\qq\qq\qq\qq\qq\qq\qq \ae \ t\in[0,\i),~\as,\ i=1,2.\ea\ee
Now, \rf{Th-Th*-closed}, \rf{definition} and \rf{stationarity condition-closed-1} yields
$$\ba{ll}
\ds 0=\h R_{ii}\dbE[v^*]\1n+\1n\h B_i^\top\1n\big(\bar\eta_i\1n+\1n\h P_i\dbE[X^*]\big)
  \1n+\1n\h D_i^\top\1n\big\{\dbE[\z_i]\1n+\1n P_i\dbE[\si]\1n+\1n P_i\h C\dbE[X^*]\1n+\1n P_i\h D\dbE[v^*]\big\}\1n+\1n\big[\h S_{ii}\1n+\1n\h R_{ii}\bar\Th^*\big]\dbE[X^*]\1n+\1n\dbE[\rho_{ii}]\\
\ns\ds=\big(\h R_{ii}+\h D_i^\top P_i\h D\big)\dbE[v^*]+\big\{\h B_i^\top\h P_i
+\h D_i^\top P_i\h C+\h S_{ii}+\big(\h R_{ii}+\h D_i^\top P_i\h D\big)\Th^*\big\}\dbE[X^*]\\
\ns\ds\qq+\h B_i^\top\bar\eta_i+\h D_i^\top\big(\dbE[\z_i]+P_i\dbE[\si]\big)+\dbE[\rho_{ii}]\\
\ns\ds=\big(\h R_{ii}+\h D_i^\top P_i\h D\big)\dbE[v^*]+\h B_i^\top\bar\eta_i+\h D_i^\top\big(\dbE[\z_i]+P_i\dbE[\si]\big)+\dbE[\rho_{ii}],
 \q\ae \ t\in[0,\i),\ i=1,2.\ea$$
Thus \rf{Ev*} holds. Furthermore, \rf{Th-Th*-closed}, \rf{definition} and \rf{stationarity condition-closed-2} yields
\bel{3.68}\ba{ll}
\ds 0=R_{ii}(v^*-\dbE[v^*])+B_i^\top\big[\a_i+P_i\big(X^*-\dbE[X^*]\big)\big]+D_i^\top\big\{
\z_i-\dbE[\z_i]+P_iC_{\Th^*}\big(X^*-\dbE[X^*]\big)\\
\ns\ds\q+P_iD\big(v^*-\dbE[v^*]\big)+P_i\big(\si-\dbE[\si]\big)\big\}+(S_{ii}+R_{ii}\Th^*)
\big(X^*-\dbE[X^*]\big)+\rho_{ii}-\dbE[\rho_{ii}]\\
\ns\ds=(R_{ii}+D_i^\top P_iD)(v^*-\dbE[v^*])+B_i^\top\a_i+D_i^\top(\z_i-\dbE[\z_i])+D_i^\top P_i\big(\si-\dbE[\si]\big)\\
\ns\ds\q+\big[B_i^\top P_i+D_i^\top P_iC+S_{ii}+(R_{ii}+D_i^\top P_iD)\Th^*\big]\big(X^*-\dbE[X^*]\big)+\rho_{ii}-\dbE[\rho_{ii}]\\
\ns\ds=(R_{ii}+D_i^\top P_iD)\big(v^*-\dbE[v^*]\big)+B_i^\top\a_i+D_i^\top(\z_i-\dbE[\z_i])+D_i^\top P_i\big(\si-\dbE[\si]\big)+\rho_{ii}-\dbE[\rho_{ii}],\\
\ns\ds\qq\qq\qq\qq\qq\qq\qq\qq\qq\qq\q\ae \ t\in[0,\i), \ \as,\ i=1,2.\ea\ee
Thus, by \rf{3.68}, we get \rf{v*-Ev*}.

\ms

{\it Sufficiency}. We choose any $x\in\dbR^n$ and $v(\cd)\equiv\begin{pmatrix}v_1(\cd)\\ v_2(\cd)\end{pmatrix}\in L^2_\dbF(\dbR^m)$. Denote $w(\cd)\equiv\begin{pmatrix}v_1(\cd)\\ v_2^*(\cd)\end{pmatrix}$ and let $X^w(\cd)\equiv X\big(\cd;x,\BTh_1^*,v_1(\cd);\BTh_2^*,v_2^*(\cd)\big)$ be the solution to the state equation
\bel{Xw}\left\{\2n\ba{ll}
\ds dX^w(t)=\big\{A_{\Th^*}X^w+\bar A_{\BTh^*}\dbE[X^w]+Bw+\bar B\dbE[w]+b\big\}dt\\
\ns\ds\qq\qq\qq+\big\{C_{\Th^*}X^w+\bar C_{\BTh^*}\dbE[X^w]+Dw+\bar D\dbE[w]+\si\big\}dW(t),\qq t\ges0,\\
\ns\ds X^w(0)=x,\ea\right.\ee
corresponding to $x$ and $(\BTh_1^*,v_1(\cd),\BTh_2^*,v_2^*(\cd))$. Similarly as \rf{cost-LQ-game-closed-1}, we have
$$\ba{ll}
\ds J_1\big(x;\Th^*\big(X^w(\cd)-\dbE[X^w(\cd)]\big)+\bar\Th^*\dbE[X^w(\cd)]+w(\cd)\big)\\
\ns\ds=\dbE\int_0^\i\Big[\blan(Q_1+S_1^\top\Th^*+(\Th^*)^\top S_1+(\Th^*)^\top R_1\Th^*)(X^w-\dbE[X^w]),X^w-\dbE[X^w]\bran\\
\ns\ds\qq\qq+2\blan(S_1+R_1\Th^*)^\top(w-\dbE[w]),X^w-\dbE[X^w]\bran+\blan R_1(w-\dbE[w]),w-\dbE[w]\bran\\
\ns\ds\qq\qq+2\blan q_1+(\Th^*)^\top\rho_1,X^w-\dbE[X^w]\bran+2\blan\rho_1,w-\dbE[w]\bran\\
\ns\ds\qq\qq+\blan\h Q_1+\h S_1^{\,\top}\bar\Th^*+(\bar\Th^*)^\top\h S_1+(\bar\Th^*)^\top\h R_1\bar\Th^*\big)\dbE[X^w],\dbE[X^w]\bran+\blan\h S_1+\h R_1\bar\Th^*\big)^\top\dbE[w],\dbE[X^w]\bran\\
\ns\ds\qq\qq+\blan\h R_1\dbE[w],\dbE[w]\bran+2\blan\dbE[q_1]+(\bar\Th^*)^\top\dbE[\rho_1],\dbE[X^w]\bran
+2\blan\dbE[\rho_1],\dbE[w]\bran\Big]dt.\ea$$
Applying It\^{o}'s formula to
$$\blan P_1\big(X^w(\cd)-\dbE[X^w(\cd)]\big)+2\eta_1(\cd),X^w(\cd)-\dbE[X^w(\cd)]\bran+\blan \h P_1\dbE[X^w(\cd)]+2\bar\eta_1(\cd),\dbE[X^w(\cd)]\bran,$$
we have
$$\ba{ll}
\ds\q-\blan\h P_1x+2\bar\eta_1(0),x\bran=\dbE\int_0^\i\Big[\blan\big(P_1A_{\Th^*}+A_{\Th^*}^\top P_1+C_{\Th^*}^\top P_1C_{\Th^*}\big)(X^w-\dbE[X^w]),X^w-\dbE[X^w]\bran\\
\ns\ds\qq+2\blan\big(P_1B+C_{\Th^*}^\top P_1D\big)(w-\dbE[w]),X^w-\dbE[X^w]\bran+\blan D^\top P_1D(w-\dbE[w]),w-\dbE[w]\bran\\
\ns\ds\qq-2\blan(\Th^*)^\top\rho_1+q_1,X^w-\dbE[X^w]\bran+2\blan B^\top\eta_1+D^\top\z_1+D^\top P_1\si,w-\dbE[w]\bran\\
\ns\ds\qq+\blan\big(\h P_1\h A_{\bar\Th^*}+\h A_{\bar\Th^*}^{\,\top}\h P_1
  +\h C_{\bar\Th^*}^{\,\top}P_1\h C_{\bar\Th^*}\big)\dbE[X^w],\dbE[X^w]\bran+2\blan\big(\h P_1\h B+\h C_{\bar\Th^*}^\top P_1\h D\big)\dbE[w],\dbE[X^w]\bran\\
\ns\ds\qq+\blan\h D^\top P_1\h D\dbE[w],\dbE[w]\bran-2\blan(\bar \Th^*)^\top\dbE[\rho_1]+\dbE[q_1],\dbE[X^w]\bran
  +2\blan\h B^\top\bar\eta_1+\h D^\top\big(P_1\dbE[\si]+\dbE[\z_1]),\dbE[w]\bran\\
\ns\ds\qq+\blan P_1\si,\si\bran+2\blan\eta_1,b-\dbE[b]\bran+2\blan\z_1,\si\bran
  +2\blan\bar\eta_1,\dbE[b]\bran\Big]dt.
\ea$$
Combining the above two equalities, together with conditions \rf{v*-Ev*} and \rf{Ev*}, we obtain
$$\ba{ll}
\ds J_1\big(x;\Th^*\big(X^w(\cd)-\dbE[X^w(\cd)]\big)+\bar\Th^*\dbE[X^w(\cd)]+w(\cd)\big)
-\blan\h P_1x+2\bar\eta_1(0),x\bran\\
\ns\ds=\dbE\int_0^\i\Big[\blan(R_1+D^\top P_1D)(w-\dbE[w]),w-\dbE[w]\bran+2\blan B^\top\eta_1+D^\top\z_1+D^\top P_1\si+\rho_1,w-\dbE[w]\bran\\
\ns\ds\qq\qq+\blan\big[\h R_1+\h D^\top P_1\h D\big]\dbE[w],\dbE[w]\bran
  +2\blan\h B^\top\bar\eta_1+\h D^\top(P_1\dbE[\si]+\dbE[\z_1])+\dbE[\rho_1],\dbE[w]\bran\\
\ns\ds\qq\qq+\blan P_1\si,\si\bran+2\blan\eta_1,b-\dbE[b]\bran
+2\blan\z_1,\si\bran+2\blan\bar\eta_1,\dbE[b]\bran\Big]dt\\
\ns\ds=\dbE\int_0^\i\Big[\blan(R_{111}+D_1^\top P_1D_1)(v_1-\dbE[v_1]),
v_1-\dbE[v_1]\bran+2\blan\big(R_{112}+D_1^\top P_1D_2\big)(v_2^*-\dbE[v_2^*]),v_1-\dbE[v_1]\bran\\
\ns\ds\qq\qq+\blan\big(R_{122}+D_2^\top P_1D_2\big)(v_2^*-\dbE[v_2^*]),v_2^*-\dbE[v_2^*]\bran\\
\ns\ds\qq\qq+2\blan B_1^\top(\eta_1-\dbE[\eta_1])+D_1^\top(\z_1-\dbE[\z_1])+D_1^\top P_1(\si-\dbE[\si])+\rho_{11}-\dbE[\rho_{11}],v_1-\dbE[v_1]\bran\\
\ns\ds\qq\qq+2\blan B_2^\top(\eta_1-\dbE[\eta_1])+D_2^\top(\z_1-\dbE[\z_1])+D_2^\top P_1(\si-\dbE[\si])+\rho_{12}-\dbE[\rho_{12}],v_2^*-\dbE[v_2^*]\bran\\
\ns\ds\qq\qq+\blan(\h R_{111}+\h D_1^\top P_1\h D_1)\dbE[v_1],\dbE[v_1]
\bran+2\blan(\h R_{112}+\h D_1^\top P_1\h D_2)\dbE[v_2^*],\dbE[v_1]\bran\\
\ns\ds\qq\qq+\blan(\h R_{122}+\h D_2^\top P_1\h D_2)\dbE[v_2^*],
\dbE[v_2^*]\bran+2\blan\h B_1^\top\bar\eta_1+\h D_1^\top \big(P_1\dbE[\si]+\dbE[\z_1]\big)+\dbE[\rho_{11}],\dbE[v_1]\bran\\
\ns\ds\qq\qq+2\blan\h B_2^\top\bar\eta_1+\h D_2^\top(P_1\dbE[\si]+\dbE[\z_1])+\dbE[\rho_{12}],\dbE[v_2^*]\bran\\
\ns\ds\qq\qq+\blan P_1\si,\si\bran+2\blan\eta_1,b-\dbE[b]\bran+2\blan\z_1,
\si\bran+2\blan\bar\eta_1,\dbE[b]\bran\Big]dt\\
\ns\ds=\dbE\int_0^\i\Big[\blan\big(R_{111}+D_1^\top P_1D_1\big)\big[v_1-\dbE[v_1]-(v_1^*-\dbE[v_1^*])\big],v_1-\dbE[v_1]-(v_1^*-\dbE[v_1^*])\bran\\
\ns\ds\qq\qq-\blan\big(R_{111}+D_1^\top P_1D_1\big)(v_1^*-\dbE[v_1^*]),v_1^*-\dbE[v_1^*]\bran
  +\blan\big(R_{122}+D_2^\top P_1D_2\big)(v_2^*-\dbE[v_2^*]),v_2^*-\dbE[v_2^*]\bran\\
\ns\ds\qq\qq+2\blan B_2^\top(\eta_1-\dbE[\eta_1])+D_2^\top(\z_1-\dbE[\z_1])+D_2^\top\mathrm{P}_1(\si-\dbE[\si])
+\rho_{12}-\dbE[\rho_{12}],v_2^*-\dbE[v_2^*]\bran\\
\ns\ds\qq\qq+\blan\big[\h R_{111}+\h D_1^\top P_1\h D_1\big](\dbE[v_1]-\dbE[v_1^*]),\dbE[v_1]-\dbE[v_1^*]\bran-\blan\big[\h R_{111}+\h D_1^\top P_1\h D_1\big]\dbE[v_1^*],\dbE[v_1^*]\bran\\
\ns\ds\qq\qq+\blan\big[\h R_{122}+\h D_2^\top P_1\h D_2\big]\dbE[v_2^*],\dbE[v_2^*]\bran+2\blan\h B_2^\top\bar\eta_1+\h D_2^\top P_1\dbE[\si]+\dbE[\z_1])+\dbE[\rho_{12}],\dbE[v_2^*]\bran\\
\ns\ds\qq\qq+\blan P_1\si,\si\bran+2\blan\eta_1,b-\dbE[b]\bran+2\blan\z_1,\si\bran+2\blan\bar\eta_1,\dbE[b]\bran\Big]dt.
\ea$$
Consequently, one gets
$$\ba{ll}
\ds J_1\big(x;\Th^*\big\{X^w(\cd)-\dbE[X^w(\cd)]\big\}+\bar\Th^*\dbE[X^w(\cd)]+w(\cd)\big)
-J_1\big(x;\Th^*\big\{X^*(\cd)-\dbE[X^*(\cd)]\big\}+\bar\Th^*\dbE[X^*(\cd)]+v^*(\cd)\big)\\
\ns\ds=\dbE\int_0^\i\Big[\blan\big(R_{111}+D_1^\top P_1D_1\big)\big[v_1-\dbE[v_1]-(v_1^*-\dbE[v_1^*])\big],v_1-\dbE[v_1]-(v_1^*-\dbE[v_1^*])\bran\\
\ns\ds\qq\qq+\blan\big[\h R_{111}+\h D_1^\top P_1\h D_1\big](\dbE[v_1]-\dbE[v_1^*]),\dbE[v_1]-\dbE[v_1^*]\bran\Big]dt\ges0,\ea$$
since \rf{Sigma} holds with $i=1$.

\ms

Similarly, by \rf{Sigma} with $i=2$, for any $\bar w(\cd)\equiv(v_1^*(\cd),v_2(\cd))$ and $X^{\bar w}(\cd)\equiv X\big(\cd\,;x,\BTh_1^*,v_1^*(\cd);\BTh_2^*,v_2(\cd)\big)$, we can prove that the following holds:
$$\ba{ll}
\ds J_2\big(x;\Th^*\big\{X^{\bar w}(\cd)-\dbE[X^{\bar w}(\cd)]\big\}+\bar\Th^*\dbE[X^{\bar w}(\cd)]+\bar w(\cd)\big)-J_2\big(x;\Th^*\big\{X^*(\cd)-\dbE[X^*(\cd)]\big\}+\bar\Th^*\dbE[X^*(\cd)]
+v^*(\cd)\big)\\
\ns\ds=\dbE\int_0^\i\Big[\blan\big(R_{222}+D_2^\top\mathrm{P}_2D_2\big)\big[v_2-\dbE[v_2]-(v_2^*-\dbE[v_2^*])\big],v_2-\dbE[v_2]-(v_2^*-\dbE[v_2^*])\bran\\
\ns\ds\qq\qq+\blan\big[\h R_{222}+\h D_2^\top P_2\h D_2\big](\dbE[v_2]-\dbE[v_2^*]),\dbE[v_2]-\dbE[v_2^*]\bran\Big]dt\ges0.\ea$$
By Definition \ref{D-3.2}, this proves the sufficiency. The proof is complete.
\endpf

\ms
To conclude this section, let us rewrite system  \rf{AE-LQ-game-closed} in a more compact form so that one can see an interesting feature of it. We define
$$\BP\deq\begin{pmatrix}P_1&0\\0&P_2\end{pmatrix},\q \h\BP\deq\begin{pmatrix}\h P_1&0\\0&\h P_2\end{pmatrix}.$$
Note that \rf{Th-Th*-closed} is equivalent to (recalling the notation introduced in \rf{BA} and \rf{IJ})
\bel{Th-Th*-closed-BJ}
0=\begin{pmatrix}B_1^\top P_1+D_1^\top P_1C+S_{11}\\B_2^\top P_2+D_2^\top P_2C+S_{22}\end{pmatrix}
+\begin{pmatrix}R_{11}+D_1^\top P_1D\\R_{22}+D_2^\top P_2D\end{pmatrix}\Th^*\equiv\BJ^\top\big(\BB^\top\BP+\BD^\top\BP\BC+\BS\big)\BI_n+\BSi\Th^*,\ee
and
\bel{Th-Th*-closed-BBJJ}
0=\begin{pmatrix}\h B_1^\top\h P_1+\h D_1^\top P_1\h C+\h S_{11}\\
  \h B_2^\top\h P_2+\h D_2^\top P_2\h C+\h S_{22}\end{pmatrix}
  \1n+\1n\begin{pmatrix}\h R_{11}+\h D_1^\top P_1\h D\\ \h R_{22}+\h D_2^\top P_2
  \h D\end{pmatrix}\bar\Th^*\equiv\BJ^\top\big(\h\BB^\top\h\BP+\h\BD^\top\BP\h \BC+\h\BS\big)\BI_n+\bar\BSi\bar\Th^*,\ee
where
$$\BSi\deq\BJ^\top\1n(\BR\1n+\1n\BD^\top\1n\BP\BD)\BI_m\1n\equiv\1n\begin{pmatrix}R_{11}\1n+\1n
D_1^\top\1n P_1D\\R_{22}\1n+\1n D_2^\top\1n P_2D\end{pmatrix},\q \bar\BSi\deq\BJ^\top\1n\big(\h\BR\1n+\1n\h\BD^\top\1n\BP\h\BD\big)\BI_m\1n
\equiv\1n\begin{pmatrix}\h R_{11}\1n+\1n\h D_1^\top\1n P_1\h D\\ \h R_{22}\1n+\1n\h D_2^\top\1n P_2\h D\end{pmatrix}\1n\in\1n\dbR^{m\times m}.$$
If we assume both $\BSi$ and $\bar\BSi$ are invertible, then we have
\bel{Th-Th*-BJ-closed}
\Th^*=-\BSi^{-1}\BJ^\top\big(\BB^\top\BP+\BD^\top\BP\BC+\BS\big)\BI_n,\q
\bar\Th^*=-\bar\BSi^{-1}\BJ^\top\big(\h\BB^\top\h\BP
+\h\BD^\top\BP\h\BC+\h\BS\big)\BI_n\in\dbR^{m\times n}.\ee
On the other hand, \rf{AE-LQ-game-closed} can be written as
\bel{AE-LQ-game-closed--}\left\{\2n\ba{ll}
\ns\ds\BP\BA+\BA^\top\BP+\BC^\top\BP\BC+\BQ+(\BPhi^*)^\top(\BR+\BD^\top\BP\BD)\BPhi^*\\
\ns\ds\qq+(\BP\BB+\BC^\top\BP\BD+\BS^\top)\BPhi^*+(\BPhi^*)^\top(\BB^\top\BP+\BD^\top\BP\BC+\BS)=0,\\
\ns\ds\h\BP\h\BA+\h\BA^\top\h\BP+\h\BC^\top\BP\h\BC+\h\BQ+(\bar\BPhi^*)^\top\big(\h\BR+\h\BD^\top\BP\h\BD\big)\bar\BPhi^*\\
\ns\ds\qq+\big(\h\BP\h\BB+\h\BC^\top\BP\h\BD+\h\BS^\top\big)\bar\BPhi^*
+(\bar\BPhi^*)^\top\big(\h\BB^\top\h\BP+\h\BD^\top\BP\h\BC+\h\BS\big)=0,\ea\right.\ee
with
$$\BPhi^*\deq\begin{pmatrix}\Th^*&0\\0&\Th^*\end{pmatrix},\q\bar\BPhi^*
\deq\begin{pmatrix}\bar\Th^*&0\\0&\bar\Th^*\end{pmatrix}.$$
Clearly, both equations in \rf{AE-LQ-game-closed--} are symmetric with solutions $\BP,\h\BP\in\dbS^{2n}$.
Recall that system \rf{ARE-LQ-game-BP} (or, equivalently, \rf{AE-LQ-game-BP}) for closed-loop representation of an open-loop Nash equilibrium are coupled, and are not symmetric. Therefore, the closed-loop representation of open-loop Nash equilibria is different from the outcome of closed-loop Nash equilibria, for Problem (MF-SDG), in general.

\section{Mean-Field LQ Zero-Sum Stochastic Differential Games}
In this section, we will look at the situation for mean-field LQ zero-sum stochastic differential games. According to \rf{zero-zero}, let us simplify some notation:
\bel{coeff-zero-sum}\left\{\2n\ba{ll}
\ds Q_1=-Q_2\equiv Q,\q \bar Q_1=-\bar Q_2\equiv\bar Q,\q q_1(\cd)=-q_2(\cd)\equiv q(\cd),\\
\ns\ds S_1\equiv\begin{pmatrix}S_{11}\\ S_{12}\end{pmatrix}=-\begin{pmatrix}S_{21}\\ S_{22}\end{pmatrix}\equiv-S_2\equiv\begin{pmatrix}S^1\\ S^2\end{pmatrix}\equiv S,\q
 \bar S_1\equiv\begin{pmatrix}\bar S_{11}\\\bar S_{12}\end{pmatrix}=-\begin{pmatrix}\bar S_{21}\\\bar S_{22}\end{pmatrix}\equiv-\bar S_2\equiv\begin{pmatrix}\bar S^1\\\bar S^2\end{pmatrix}\equiv\bar S,\\
\ns\ds R_1\equiv\begin{pmatrix}R_{111}&R_{112}\\ R_{121}& R_{122}\end{pmatrix}=-\begin{pmatrix}R_{211}&R_{212}\\ R_{221}& R_{222}\end{pmatrix}
\equiv-R_2\equiv R\equiv\begin{pmatrix}R_{11}&R_{12}\\ R_{21}& R_{22}\end{pmatrix},\\
\ns\ds\bar R_1\equiv\begin{pmatrix}\bar R_{111}&\bar R_{112}\\\bar R_{121}&\bar R_{122}\end{pmatrix}=-\begin{pmatrix}\bar R_{211}&\bar R_{212}\\\bar R_{221}&\bar R_{222}\end{pmatrix}
\equiv-\bar R_2\equiv\bar R\equiv\begin{pmatrix}\bar R_{11}&\bar R_{12}\\\bar R_{21}&\bar R_{22}\end{pmatrix},\\
\ns\ds\rho_1(\cd)\equiv\begin{pmatrix}\rho_{11}(\cd)\\ \rho_{12}(\cd)\end{pmatrix}=-\begin{pmatrix}\rho_{21}(\cd)\\ \rho_{22}(\cd)\end{pmatrix}
\equiv-\rho_2(\cd)\equiv\begin{pmatrix}\rho^1(\cd)\\ \rho^2(\cd)\end{pmatrix}\equiv\rho(\cd).
\ea\right.\ee
Then the cost functional will be the following:
\bel{cost-zero-sum}\ba{ll}
\ds J_1(x;u_1(\cd),u_2(\cd))=-J_2(x;u_1(\cd),u_2(\cd))\equiv J(x;u_1(\cd),u_2(\cd))\equiv J(x;u(\cd))\\
\ns\ds\equiv\dbE\int_0^\i\bigg[\bigg\langle\2n
\begin{pmatrix}Q&S^\top\\ S&R\\\end{pmatrix}
\begin{pmatrix}X(t)\\ u(t)\end{pmatrix},\begin{pmatrix}X(t)\\ u(t)\end{pmatrix}\2n\bigg\rangle
+2\bigg\langle\2n\begin{pmatrix}q(t)\\ \rho(t)\end{pmatrix},
\begin{pmatrix}X(t)\\ u(t)\end{pmatrix}\2n\bigg\rangle\\
\ns\ds\qq\qq\qq+\bigg\langle\2n
\begin{pmatrix}\bar Q&\bar S^\top\\ \bar S&\bar R\end{pmatrix}
\begin{pmatrix}\dbE[X(t)]\\ \dbE[u(t)]\end{pmatrix},\begin{pmatrix}\dbE[X(t)]\\ \dbE[u(t)]\end{pmatrix}
\2n\bigg\rangle\bigg]dt.\ea\ee
Similar to Problem (MF-SDG), we assume that (H2) holds. Then, for any $x\in\dbR^n$, we define $\sU_{ad}(x)$ and $\sU_{ad}^i(x)$ ($i=1,2$) similar to \rf{Uad-LQ-game-u}--\rf{Uad-LQ-game-2}.
Let us state the following zero-sum problem.

\ms

{\bf Problem (MF-SDG)$_0$.}\q For any initial state $x\in\dbR^n$, Player 1 wants to find a control $u_1^*(\cd)$ to minimize the cost functional $J(x;u_1(\cd),u_2(\cd))$, and Player 2 wants to find a control $u_2^*(\cd)$ to maximize $J(x;u_1(\cd),u_2(\cd))$ respectively, subject to \rf{state-LQ-game} (or equivalently, \rf{state-LQ-game-u}) such that $(u_1^*(\cd),u_2^*(\cd))\in\mathcal{U}_{ad}(x)$.

\ms

Next, we introduce the following definitions.
\bde{D-4.1} \rm
(i) A pair $(u_1^*(\cd),u_2^*(\cd))\in\mathcal{U}_{ad}(x)$ is called an {\it open-loop saddle point} of Problem (MF-SDG)$_0$ for the initial state $x\in\dbR^n$ if
\bel{open-loop saddle point}
J(x;u_1^*(\cd),u_2(\cd))\les J(x;u_1^*(\cd),u_2^*(\cd))\les J(x;u_1(\cd),u_2^*(\cd)),\qq\forall(u_1(\cd),u_2(\cd))\in\sU_{ad}(x).\ee

\ms

(ii) A closed-strategy $(\BTh^*,v^*(\cd))\in\sS[A,\bar A,C,\bar C;B,\bar B,D,\bar D]\times L^2_\dbF(\dbR^m)$ is called a {\it closed-loop saddle point} of Problem (MF-SDG)$_0$ if
for any $x\in\dbR^n$, $\BTh_1\in\sS^1(\BTh_2^*)$, $\BTh_2\in\sS^2(\BTh_1^*)$, $v_1(\cd)\in L^2_\dbF(\dbR^{m_1})$ and $v_2(\cd)\in L^2_\dbF(\dbR^{m_2})$,
\bel{closed-loop saddle point}
J\big(x;\BTh_1^*,v_1^*(\cd);\BTh_2,v_2(\cd)\big)\1n\les\1n J\big(x;\BTh^*,v^*(\cd)\big)\1n\les\1n J\big(x;\BTh_1,v_1(\cd);\BTh_2^*,v_2^*(\cd)\big).\ee
\ede

Similar as in Section 3, it is easy to see that $(\BTh_1^*,v_1^*(\cd);\BTh_2^*,v_2^*(\cd))$ is a closed-loop saddle point of Problem (MF-SDG)$_0$ if and only if one of the following holds:

\ms

(i)\ For any $v_1(\cd)\in L^2_\dbF(\dbR^{m_1})$ and $v_2(\cd)\in L^2_\dbF(\dbR^{m_2})$,
\bel{closed-loop saddle point-1}
J\big(x;\BTh_1^*,v_1^*(\cd);\BTh_2^*,v_2(\cd)\big)\1n\les\1n J\big(x;\BTh^*,v^*(\cd)\big)\1n\les\1n J\big(x;\BTh_1^*,v_1(\cd);\BTh_2^*,v_2^*(\cd)\big).\ee

(ii)\ For any $u_1(\cd)\in L^2_\dbF(\dbR^{m_1})$ and $u_2(\cd)\in L^2_\dbF(\dbR^{m_2})$,
\bel{closed-loop saddle point-2}
J\big(x;\BTh_1^*,v_1^*(\cd);u_2(\cd)\big)\les J\big(x;\BTh^*,v^*(\cd)\big)\les J\big(x;u_1(\cd);\BTh_2^*,v_2^*(\cd)\big).\ee

\bde{D-4.2} \rm
The following maps
$$\left\{\2n\ba{ll}
\ds V^+(x)\deq\inf_{u_1(\cd)\in\sU_{ad}^1(x)}\sup_{u_2(\cd)\in \sU_{ad}^2(x)}J(x;u_1(\cd),u_2(\cd)),\\
\ns\ds V^-(x)\deq\sup_{u_2(\cd)\in\sU_{ad}^2(x)}\inf_{u_1(\cd)\in \sU_{ad}^1(x)}J(x;u_1(\cd),u_2(\cd)),
\ea\right.$$
are called the {\it upper value function} and the {\it lower value function} of Problem (MF-SDG)$_0$, respectively. In the case that
\begin{equation*}
V^+(x)=V^-(x)\equiv V(x),
\end{equation*}
we call the map $x\mapsto V(x)$ the value function of Problem (MF-SDG)$_0$.
\ede

Now, let $(\BTh^*,v^*(\cd))\in\dbR^{2m\times n}\times L^2_\dbF(\dbR^m)$, and assume the open-loop saddle points of Problem (MF-SDG)$_0$ admit the closed-loop representation.
Then by \rf{coeff-zero-sum}, we see that $(P_1,\h P_1)=(-P_2,-\h P_2)\equiv(P,\h P)$ satisfy the same equations (see \rf{AE-LQ-game-component})
\bel{AE-zero-sum-open-loop}\left\{\2n\ba{ll}
\ds PA+A^\top P+C^\top PC+Q+\big(PB+C^\top PD+S^\top\big)\Th^*=0,\\
\ns\ds\h P\h A+\h A^{\,\top}\h P+\h C^\top P\h C+\h Q+\big(\h P\h B+\h C^\top P\h D+\h S^{\,\top}\big)\bar\Th^*=0,\ea\right.\ee
and \rf{Th-Th*-component} is equivalent to
\bel{Th-Th*-zero-sum-open-loop}\left\{\2n\ba{ll}
\ds\Si_o\Th^*+B^\top P+D^\top PC+S=0,\qq\Si_o\equiv R+D^\top PD,\\
\ns\ds\bar\Si_o\bar\Th^*+\h B^\top\h P+\h D^\top P\h C+\h S=0,\qq\bar\Si_o\equiv\h R+\h D^\top P\h D.\ea\right.\ee
These equations (for $\Th^*$ and $\bar\Th^*$) are solvable if and only if
\bel{solvable condition for AE-zero-sum-open-loop}\left\{\2n\ba{ll}
\ds\sR\big(B^\top P+D^\top PC+S\big)\subseteq\sR(\Si_o),\qq\Th^*=-\Si_o^\dag(B^\top P+D^\top PC+S)+\big(I-\Si_o^\dag\Si_o\big)\th,\\
\ns\ds\sR\big(\h B^\top\h P+\h D^\top P\h C+\h S\,\big)\subseteq\sR(\bar\Si_o),\qq
\bar\Th^*=-\bar\Si_o^\dag\big(\h B^\top\h P+\h D^\top P\h C+\h S\,\big)+(I-\bar\Si_o^\dag\bar\Si_o)\bar\th,\ea\right.\ee
for some $\th,\bar\th\in\dbR^{m\times n}$ are chosen such that $(\Th^*,\bar\Th^*)\in\sS[A,\bar A,C,\bar C;B,\bar B,D,\bar D]$.
Putting \rf{solvable condition for AE-zero-sum-open-loop} into \rf{AE-zero-sum-open-loop}, it yields
\bel{ARE-zero-sum-open-loop}\left\{\2n\ba{ll}
\ds PA+A^\top P+C^\top P C+Q-\big(PB+C^\top P D+S^\top\big)\Si_o^\dag(B^\top P+D^\top PC+S)=0,\\
\ns\ds\h P\h A+\h A^{\,\top}\h P+\h C^\top P\h C+\h Q-(\h P\h B+\h C^\top P\h D+\h S^\top)\bar\Si_o^\dag(\h B^\top\h P+\h D^\top P\h C+\h S)=0,
\ea\right.\ee
which is a system of coupled AREs, and both of them are symmetric. Thus, we have $P,\h P\in\dbS^n$.
Next, by \rf{coeff-zero-sum} again, from the componentwise form of \rf{BSDE-eta-game-BJ}, it is easy to see that
$(\eta_1(\cd),\z_1(\cd))=(-\eta_2(\cd),-\z_2(\cd))\\\equiv(\eta_o(\cd),\z_o(\cd))\in\sX[0,\i)\times L^2_\dbF(\dbR^n)$ satisfies
\bel{BSDE-eta-zero-sum-open-loop}\ba{ll}
\ds-d\eta_o(t)=\Big\{A^\top\eta_o(t)-\big(PB+C^\top PD+S^\top\big)\Si_o^\dag\big\{B^\top\eta_o(t)+ D^\top[\z_o(t)+P\si(t)]+\rho(t)\big\}\\
\ns\ds\qq\qq\q+C^\top\big[\z_o(t)+P\si(t)\big]+Pb(t)+q(t)\Big\}dt-\z_o(t)dW(t),\q t\ges0,\ea\ee
with constraint
\bel{range-eta-zero-sum-open-loop}\ba{ll}
\ds B^\top(\eta_o-\dbE[\eta_o])+D^\top(\z_o-\dbE[\z_o])+D^\top P(\si-\dbE[\si])+\rho-\dbE[\rho]\in\mathcal{R}(\Si_o),\q\ae,~\as\ea\ee

Likewise, by \rf{coeff-zero-sum}, from the componentwise form of \rf{BODE-bar-eta-BJ}, we see that $\bar\eta_1(\cd)=-\bar\eta_2(\cd)\equiv\bar\eta_o(\cd)\in L^2(\dbR^n)$ satisfy the following:
\bel{BODE-bar-eta-zero-sum-open-loop}\ba{ll}
\ds\dot{\bar\eta}_o(t)+\h A^{\,\top}\bar\eta_o(t)-(\h P\h B+\h C^\top P\h D+\h S^\top)\bar\Si_o^\dag\big(\h B^\top\bar\eta_o(t)+\dbE[\bar\z_o(t)]+P\dbE[\si(t)]+\dbE[\rho(t)]\big)\\
\ns\ds\qq+\h C^\top\big(\dbE[\bar\z_o(t)]+P\dbE[\si(t)]\big)+\dbE[q(t)]+\h P\dbE[b(t)]=0,\qq t\ges0,\ea\ee
with constraint
\bel{range-bar-eta-zero-sum-open-loop}\ba{ll}
\ds\h B^\top\bar\eta_o+\h D^\top\dbE[\z_o]+\h D^\top P\dbE[\si]+\dbE[\rho]\in\sR(\bar\Si_o),\q\ae\ea\ee
Finally, from \rf{closed-loop ui}, it is direct that the closed-loop representation of an open-loop saddle point is given by
\bel{closes-loop representation of an open-loop saddle point}\ba{ll}
\ds u^*(\cd)=\big\{-\Si_o^\dag(B^\top P+D^\top PC+S)+\big(I-\Si_o^\dag\Si_o\big)\th\big\}\big(X^*(\cd)-\dbE[X^*(\cd)]\big)\\
\ns\ds\qq\qq+\big\{-\bar\Si_o^\dag\big(\h B^\top\h P+\h D^\top P\h C+\h S\,\big)+\big(I-\bar\Si_o^\dag\bar\Si_o\big)\bar\th\big\}\dbE[X^*(\cd)]\\
\ns\ds\qq\qq-\Si_o^\dag\big\{B^\top(\eta_o(\cd)-\dbE[\eta_o(\cd)])+D^\top\big[\z_o(\cd)
-\dbE[\z_o(\cd)]+P\big(\si(\cd)-\dbE[\si(\cd)]\big)\big]+\rho(\cd)-\dbE[\rho(\cd)]\big\}\\
\ns\ds\qq\qq-\bar\Si_o^\dag\big\{\h B^\top\bar\eta_o(\cd)+\h D^\top\big(P\dbE[\si(\cd)]+\dbE[\z_o(\cd)]\big)+\dbE[\rho(\cd)]+\bar\rho(\cd)\big\}\\
\ns\ds\qq\qq+\big(I-\Si_o^\dag\Si_o\big)(\nu(\cd)-\dbE[\nu(\cd)])+\big(I-\bar\Si_o^\dag\bar\Si_o\big)\bar\nu(\cd),
\ea\ee
where $\n(\cd),\bar\n(\cd)\in L^2(\dbR^m)$, and $X^*(\cd)\in\sX[0,\i)$ is the solution to \rf{closed-loop-system-game} corresponding to $u^*(\cd)$.

\ms

To summarize, we have the following result.
\bt{T-4.3} \sl
Let {\rm(H2)} hold and the initial state $x\in\dbR^n$ be given. Then any open-loop saddle point $u^*(\cd)$ of Problem (MF-SDG)$_0$ admits a closed-loop representation if and only if the following hold:

\ms

{\rm(i)}\ The following convexity-concavity condition holds: For $i=1,2$,
\bel{convexity-concavity condition}\ba{ll}
\ds(-1)^{i-1}\dbE\int_0^\i\Big[\blan QX_i(t),X_i(t)\bran+2\blan S^iX_i(t),u_i(t)\bran+\blan R_{ii}u_i(t),u_i(t)\bran+\blan \bar Q\dbE[X_i(t)],\dbE[X_i(t)]\bran\\
\ns\ds\qq+2\blan \bar S^i\dbE[X_i(t)],\dbE[u_i(t)]\bran+\blan\bar R_{ii}\dbE[u_i(t)],\dbE[u_i(t)]\bran\Big]dt\ges0,\
 \forall u(\cd)\equiv\begin{pmatrix}u_1(\cd)\\ u_2(\cd)\end{pmatrix}\in\sU_{ad}(x),\ea\ee
where $X_i(\cd)\in\sX[0,\i)$ is the solution to MF-SDE \rf{homogeneous controlled MF-SDE}.

\ms

{\rm(ii)}\ System \rf{ARE-zero-sum-open-loop} admits a static stabilizing solution $(P,\h P)\in\dbS^n\times\dbS^n$, such that the solution $(\eta_o(\cd),\z_o(\cd))\in\sX[0,\i)\times L^2_\dbF(\dbR^n)$ to \rf{BSDE-eta-zero-sum-open-loop} satisfies \rf{range-eta-zero-sum-open-loop}, and the solution $\bar\eta_o(\cd)\in L^2(\dbR^n)$ to \rf{BODE-bar-eta-zero-sum-open-loop} satisfies \rf{range-bar-eta-zero-sum-open-loop}.

\ms

In the above case, any open-loop saddle point admits the closed-loop representation \rf{closes-loop representation of an open-loop saddle point}.
\et

The following result characterizes the closed-loop saddle points of Problem (MF-SDG)$_0$.
\bt{T-4.4} \sl Problem (MF-SDG)$_0$ admits a closed-loop saddle point $(\Th^*,\bar \Th^*,v^*(\cd))$ if and only if the following statements hold:

\ms

{\rm(i)}\ The following system:
\bel{ARE-zero-sum-closed-loop}\left\{\2n\ba{ll}
\ds P_cA+A^\top P_c+C^\top P_cC+Q-(P_cB+C^\top P_cD+S^\top)\Si_c^\dag(B^\top P_c+D^\top P_cC+S)=0,\\
\ns\ds\h P_c\h A+\h A^{\,\top}\h P_c+\h C^\top P_c\h C+\h Q-(\h P_c\h B+\h C^\top P_c\h D+\h S^{\,\top})\bar\Si_c^\dag(\h B^\top\h P_c+\h D^\top P_c\h C+\h S\,)=0,\\
\ns\ds\sR(B^\top P_c+D^\top P_cC+S)\subseteq\sR(\Si_c),\qq \sR(\h B^\top\h P_c+\h D^\top P_c\h C+\h S\,)\subseteq\sR(\bar\Si_c),
\ea\right.\ee
admits a static stabilizing solution $(P_c,\h P_c)\in\dbS^n\times\dbS^n$ such that
\bel{R11-R22}\left\{\2n\ba{ll}
\ds R_{11}+D_1^\top P_cD_1\ges0,\qq\h R_{11}+\h D_1^\top P_c\h D_1\ges0,\\
\ns\ds R_{22}+D_2^\top P_cD_2\les0,\qq\h R_{22}+\h D_2^\top P_c\h D_2\les0.
\ea\right.\ee

{\rm(ii)}\ The following BSDE on $[0,\i)$:
\bel{BSDE-eta-zero-sum-closed-loop}\left\{\2n\ba{ll}
\ds-d\eta_c(t)=\big\{A^\top\eta_c(t)-\big(P_cB+C^\top P_cD+S^\top\big)\Si_c^\dag \big[B^\top\eta_c(t)+D^\top\big(\z_c(t)+P_c\si(t)\big)+\rho(t)\big]\\
\ns\ds\qq\qq\qq+C^\top\big[\z_c(t)+P_c\si(t)\big]+P_cb(t)+q(t)\big\}dt-\z_c(t)dW(t),\q t\ges0,\\
\ns\ds B^\top\big[\eta_c(t)-\dbE[\eta_c(t)]\big]+D^\top\big[\z_c(t)-\dbE[\z_c(t)]\big]+D^\top P_c\big[\si(t)-\dbE[\si(t)]\big]+\rho(t)-\dbE[\rho(t)]\in\sR\big(\Si_c\big),\\
\ns\ds\qq\qq\qq\qq\qq\qq\qq\qq\qq\qq\qq\qq\qq\qq\qq\qq\ae\ t\in[0,\i),\ \as,\ea\right.\ee
admits a solution $(\eta_c(\cd),\z_c(\cd))\in\sX[0,\i)\times L^2_\dbF(\dbR^n)$, and the following ODE:
\bel{BODE-bar-eta-zero-sum-closed-loop}\left\{\2n\ba{ll}
\ds\dot{\bar\eta}_c(t)+\h A^{\,\top}\bar\eta_c(t)-\big(\h P\h B+\h C^\top P_c\h D
  +\h S^\top\big)\bar\Si_c^\dag\big\{\h B^\top\bar\eta_c(t)+\h D^\top\big(\dbE[\bar\z_o(t)]+P_c\dbE[\si(t)]\big)+\dbE[\bar\rho(t)]\big\}\\
\ns\ds\qq+\h C^\top\big(\dbE[\bar\z_o(t)]+P_c\dbE[\si(t)]\big)+\dbE[q(t)]+\h P_c\dbE[b(t)],\qq t\ges0,\\
\ns\ds\h B^\top\bar\eta_c(t)+\h D^\top\dbE[\z_c(t)]+\h D^\top P_c\dbE[\si(t)]+\dbE[\rho(t)]\in\sR(\bar\Si_c),\q \ae \ t\in[0,\i),\ea\right.\ee
admits a solution $\bar\eta_c(\cd)\in L^2(\dbR^n)$. In the above case, the closed-loop saddle point is given by
\bel{closed-loop saddle point-3}\left\{\ba{ll}
    \Th^*=-\Si_c^\dag(B^\top P_c+D^\top P_cC+S)+\big(I-\Si_c^\dag\Si_c\big)\th,\\ [1mm]
\bar\Th^*=-\bar\Si_c^\dag\big[\h B^\top\Pi_c+\h D^\top P_c\h C+\h S\big]+\big(I-\bar\Si_c^\dag\bar\Si_c\big)\bar\th,
\ea\right.\ee
where $\th,\bar\th\in\dbR^{m\times n}$ are chosen such that $(\Th^*,\bar\Th^*)\in\sS[A,\bar A,C,\bar C;B,\bar B,D,\bar D]$, and
\bel{closed-loop saddle point-4}\ba{ll}
\ds v^*(\cd)=-\Si_c^\dag\big\{B^\top\big(\eta_c(\cd)-\dbE[\eta_c(\cd)]\big)
+D^\top\big(\z_c(\cd)-\dbE[\z_c(\cd)]\big)+D^\top P\big(\si(\cd)-\dbE[\si(\cd)]\big)+\rho(\cd)-\dbE[\rho(\cd)]\big\}\\
\ns\ds\qq\q\,\,-\bar\Si_c^\dag\big\{\h B^\top\bar\eta_c(\cd)+\h D^\top\big(P_c\dbE[\si(\cd)]+\dbE[\z_c(\cd)]\big)+\dbE[\rho(\cd)]\big\}\\
\ns\ds\qq\q\,\,+\big(I-\Si_c^\dag\Si_c\big)\big(\n(\cd)-\dbE[\n(\cd)]\big)+\big(I-\bar\Si_c^\dag\bar\Si_c\big)\bar\n(\cd),\ea\ee
for some $\n(\cd)\in L^2_\dbF(\dbR^m),\bar\n(\cd)\in L^2(\dbR^m)$.
\et

\proof Let $(\BTh^*,v^*(\cd))$ be a closed-loop saddle point of Problem (MF-SDG)$_0$. System \rf{AE-LQ-game-closed} for $(P_1,\h P_1)\equiv(-P_2,-\h P_2)=(P_c,\h P_c)$ becomes
\bel{AE-zero-sum-closed-loop}\left\{\2n\ba{ll}
\ds P_cA+A^\top P_c+C^\top P_cC+Q+(\Th^*)^\top(R+D^\top P_cD)\Th^*\\
\ns\ds\q+(P_cB+C^\top P_cD+S^\top)\Th^*+(\Th^*)^\top(B^\top P_c+D^\top P_cC+S)=0,\\
\ns\ds\h P_c\h A+\h A^{\,\top}\h P_c+\h C^\top P_c\h C+\h Q+(\bar\Th^*)^\top(\h R+\h D^\top P_c\h D)\bar\Th^*\\
\ns\ds\q+\big(\h P_c\h B+\h C^\top P_c\h D+\h S^\top\big)\bar\Th^*+(\bar\Th^*)^\top\big(\h B^\top\h P_c+\h D^\top P_c\h C+\h S\,\big)=0.
\ea\right.\ee
Thus, \rf{Th-Th*-closed} becomes
$$\left\{\2n\ba{ll}
\ds B^\top P_c+D^\top P_cC+S+\Si_c\Th^*=0,\\
\ns\ds\h B^\top\h P_c+\h D^\top P_c\h C+\h S+\bar\Si_c\bar\Th^*=0.\ea\right.$$
This system is solvable if and only if
\bel{solvable condition for AE-zero-sum-closed-loop}\left\{\2n\ba{ll}
\ds\sR\big(B^\top P_c+D^\top P_cC+S\big)\subseteq\sR(\Si_c),\qq\Th^*=-\Si_c^\dag(B^\top P_c+D^\top P_cC+S)+\big(I-\Si_c^\dag\Si_c\big)\th,\\
\ns\ds\sR\big(\h B^\top\h P_c+\h D^\top P_c\h C+\h S\,\big)\subseteq\sR(\bar\Si_c),
\qq\bar\Th^*=-\bar\Si_c^\dag\big(\h B^\top\h P_c+\h D^\top P_c\h C+\h S\,\big)+\big(I-\bar\Si_c^\dag\bar\Si_c\big)\bar\th,\ea\right.\ee
where $\th,\bar\th\in\dbR^{m\times n}$ are chosen such that $(\Th^*,\bar\Th^*)\in\sS[A,\bar A,C,\bar C;B,\bar B,D,\bar D]$.
This proves \rf{closed-loop saddle point-3}.
Putting \rf{solvable condition for AE-zero-sum-closed-loop} into \rf{AE-zero-sum-closed-loop} yields that \rf{ARE-zero-sum-closed-loop} admits a static stabilizing solution $(P_c,\h P_c)\in\dbS^n\times\dbS^n$.
\rf{R11-R22} can be easily obtained from \rf{Sigma} and \rf{coeff-zero-sum}.
\rf{BSDE-eta-zero-sum-closed-loop} and \rf{BODE-bar-eta-zero-sum-closed-loop} can be similarly proved from \rm(ii) of Theorem \ref{T-3.6}.

\ms

Finally, from \rf{v*-Ev*} we have
\bel{v*-Ev*-zero-sum}\ba{ll}
\ds v^*(\cd)-\dbE[v^*(\cd)]=-\Si_c^\dag\big\{B^\top\big(\eta_c(\cd)-\dbE[\eta_c(\cd)]\big)
+D^\top\big(\z_c(\cd)-\dbE[\z_c(\cd)]\big)+D^\top P_c\big(\si(\cd)-\dbE[\si(\cd)]\big)\\
\ns\ds\qq\qq\qq\qq\qq+\rho(\cd)-\dbE[\rho(\cd)]\big\}+\big(I-\Si_c^\dag\Si_c\big)\big(\n(\cd)-\dbE[\nu(\cd)]\big),\ea\ee
for some $\n(\cd)\in L^2_\dbF(\dbR^m)$. It follows from \rf{Ev*} that we get
\bel{Ev*-zero-sum}
\dbE[v^*(\cd)]=-\bar\Si_c^\dag\big\{\h B^\top\bar\eta_c(\cd)+\h D^\top\big(P_c\dbE[\si(\cd)]+\dbE[\z_c(\cd)]\big)+\dbE[\rho(\cd)]\big\}+\big(I-\bar\Si_c^\dag\bar\Si_c\big)\bar\nu(\cd),
\ee
for some $\bar\nu(\cd)\in L^2(\dbR^m)$. Combining the above two expressions leads to \rf{closed-loop saddle point-4}. The proof is complete.
\endpf

\ms

Comparing Theorems \ref{T-4.3} and \ref{T-4.4}, we have the following result which might not be true for general non-zero sum differential games (see \cite{Sun-Yong2019}).
\bt{T-4.5} \sl If both the closed-loop representation of open-loop saddle points and the closed-loop saddle points of Problem (MF-SDG)$_0$ exist, then the closed-loop representation coincides with the outcome of the closed-loop saddle points. In the above case, the value function admits the following representation:
\bel{V(x)-zero-sum}\ba{ll}
\ds V(x)=\blan\h Px,x\bran+2\blan\bar\eta_c(0),x\bran\\
\ns\ds\qq+\dbE\int_0^\i\Big[\blan P\si(t),\si(t)\bran+2\blan\eta_c(t),b(t)-\dbE[b(t)]\bran+2\blan\z_c(t),\si(t)\bran
+2\blan\bar\eta_c(t),\dbE[b(t)]\bran\\
\ns\ds\qq-\big|\big(\Si^\dag\big)^{\frac{1}{2}}\big\{B^\top\1n(\eta_c(t)
-\dbE[\eta_c(t)])+D^\top\big[\z_c(t)-\dbE[\z_c(t)]+P\big(\si(t)-\dbE[\si(t)]\big)\big]+\rho(t)-\dbE[\rho(t)]\big\}\big|^2\\
\ns\ds\qq-\big|\big(\bar\Si^\dag\big)^{\frac{1}{2}}\big\{\h B^\top\bar\eta_c(t)+\h D^\top\big(\dbE[\z_c(t)]
+P\dbE[\si(t)]\big)+\dbE[\rho(t)]\big\}\big|^2\Big]dt,
\ea\ee
where $\Si\equiv\Si_o\equiv\Si_c\triangleq R+D^\top PD\equiv R+D^\top P_cD$, $\bar\Si\equiv\bar\Si_o\equiv\bar\Si_c\triangleq \h R+\h D^\top P\h D\equiv\h R+\h D^\top P_c\h D$, $\h P\equiv\h P_c$.
\et
\proof
Let $(\BTh^*,v^*(\cd))$ be a closed-loop saddle point of Problem (MF-SDG)$_0$. By \rf{closed-loop saddle point-2}, the outcome
$$ u^*(\cd)=\Th^*\big(X^*(\cd)-\dbE[X^*(\cd)]\big)+\bar\Th^*\dbE[X^*(\cd)]+v^*(\cd),$$
of $(\Th^*,\bar \Th^*,v^*(\cd))$ is an open-loop saddle point of Problem (MF-SDG)$_0$, where $X^*(\cd)\in\sX[0,\i)$ is the solution to
$$\left\{\2n\ba{ll}
\ds dX^*(t)=\big\{A_{\Th^*}X^*+\bar A_{\BTh^*}\dbE[X^*]+Bv^*+\bar B\dbE[v^*]+b\big\}dt\\
\ns\ds\qq\qq\qq+\big\{C_{\Th^*}X^*+\bar C_{\BTh^*}\dbE[X^*]+Dv^*+\bar D\dbE[v^*]+\si\big\}dW(t),\q t\ges0,\\
\ns\ds X^*(0)=x.\ea\right.$$
Thus, similar to the proof of the sufficiency on Theorem \ref{T-3.3}, noting \rf{v*-Ev*-zero-sum} and \rf{Ev*-zero-sum}, we have
\bel{completion of square-closed}\ba{ll}
\ds J\big(x;\Th^*\big(X^*(\cd)-\dbE[X^*(\cd)]\big)+\bar\Th^*\dbE[X^*(\cd)]+v^*(\cd)\big)-\blan\h P_cx+2\bar\eta_c(0),x\bran\\
\ns\ds=\dbE\int_0^\i\Big[\blan\Si_c\big(v^*-\dbE[v^*]\big),v^*-\dbE[v^*]\bran+2\blan B^\top\eta_c+D^\top(\z_c+P_c\si)+\rho,v^*-\dbE[v^*]\bran\\
\ns\ds\qq\qq+\blan\bar\Si_c\dbE[v^*],\dbE[v^*]\bran+2\blan\h B^\top\bar\eta_c+\h D^\top\big(\dbE[\z_c]+P_c\dbE[\si]\big)+\dbE[\rho],\dbE[v^*]\bran\\
\ns\ds\qq\qq+\blan P_c\si,\si\bran+2\blan\eta_c,b-\dbE[b]\bran+2\blan\z_c,\si\bran+2\blan\bar\eta_c,\dbE[b]\bran\Big]dt\\
\ns\ds=\dbE\int_0^\i\Big[\blan P_c\si,\si\bran+2\blan\eta_c,b-\dbE[b]\bran+2\blan\z_c,\si\bran+2\blan\bar\eta_c,\dbE[b]\bran\\
\ns\ds\qq\qq-\big|(\Si_c^\dag)^{\frac{1}{2}}\big\{B^\top(\eta_c-\dbE[\eta_c])
+D^\top(\z_c-\dbE[\z_c])+D^\top P_c(\si-\dbE[\si])+\rho-\dbE[\rho]\big\}\big|^2\\
\ns\ds\qq\qq-\big|(\bar\Si_c^\dag)^{\frac{1}{2}}\big\{\h B^\top\bar\eta_c+\h D^\top\big(P_c\dbE[\si]+\dbE[\z_c]\big)+\dbE[\rho]\big\}\big|^2\Big]dt,\ea\ee
where $(P_c,\h P_c)$ satisfies \rf{ARE-zero-sum-closed-loop}, $(\eta_c(\cd),\z_c(\cd))$ satisfies \rf{BSDE-eta-zero-sum-closed-loop} and $\bar\eta_c(\cd)$ satisfies \rf{BODE-bar-eta-zero-sum-closed-loop}. By Theorem \ref{T-4.4}, $(\Th^*,\bar \Th^*,v^*(\cd))$ is given by \rf{closed-loop saddle point-3}--\rf{closed-loop saddle point-4}.

\ms

On the other hand, let $(P,\h P)$ be the solution pair to \rf{ARE-zero-sum-open-loop}, $(\eta_o(\cd),\z_o(\cd))$ be the solution to \rf{BSDE-eta-zero-sum-open-loop}, $\bar\eta_o(\cd)$ be the solution to \rf{BODE-bar-eta-zero-sum-open-loop}, and choose $\th,\bar\th\in\dbR^{m\times n}$, $\n(\cd),\bar\n(\cd)\in L^2(\dbR^m)$, such that the following
$$\left\{\2n\ba{ll}
\ds\Th^{**}=-\Si_o^\dag(B^\top P+D^\top PC+S)+\big(I-\Si_o^\dag\Si_o\big)\th,\q
\ds\bar\Th^{**}=-\bar\Si_o^\dag\big(\h B^\top\h P+\h D^\top P\h C+\h S\,\big)+\big(I-\bar\Si_o^\dag\bar\Si_o\big)\bar\th,\\
\ns\ds v^{**}(\cd)=-\Si_o^\dag\big\{B^\top\big(\eta_o(\cd)-\dbE[\eta_o(\cd)]\big)
+D^\top\big(\z_o(\cd)-\dbE[\z_o(\cd)]\big)+D^\top P\big(\si(\cd)-\dbE[\si(\cd)]\big)+\rho(\cd)-\dbE[\rho(\cd)]\big\}\\
\ns\ds\qq\qq-\bar\Si_o^\dag\big\{\h B^\top\1n\bar\eta_o(\cd)\1n+\1n\h D^\top\2n\big(P\dbE[\si(\cd)]\1n+\1n\dbE[\z_o(\cd)]\big)\1n+\1n\dbE[\rho(\cd)]\big\}
\1n+\1n\big(I\1n-\1n\Si_o^\dag\Si_o\big)\big(\n(\cd)\1n-\1n\dbE[\nu(\cd)]\big)
\1n+\1n\big(I\1n-\1n\bar\Si_o^\dag\bar\Si_o\big)\bar\n(\cd),\ea\right.$$
satisfies $(\Th^{**},\bar\Th^{**})\in\sS[A,\bar A,C,\bar C;B,\bar B,D,\bar D]$. For any initial state $x$, define $u^{**}(\cd)\in\sU_{ad}(x)$ as follows:
$$u^{**}(\cd)=\Th^{**}\big\{X^{**}(\cd)-\dbE[X^{**}(\cd)]\big\}+\bar\Th^{**}\dbE[X^{**}(\cd)]+v^{**}(\cd),$$
where $X^{**}(\cd)\in\sX[0,\i)$ is the solution to
$$\left\{\2n\ba{ll}
\ds dX^{**}(t)=\big\{A_{\Th^{**}}X^{**}+\bar A_{\BTh^{**}}\dbE[X^{**}]+Bv^{**}+\bar B\dbE[v^{**}]+b\big\}dt\\
\ns\ds\qq\qq\qq+\big\{C_{\Th^{**}}X^{**}+\bar C_{\BTh^{**}}\dbE[X^{**}]+Dv^{**}+\bar D\dbE[v^{**}]+\si\big\}dW(t),\q t\ges0,\\
\ns\ds X^{**}(0)=x.\ea\right.$$
By Theorem \ref{T-4.3}, $u^{**}(\cd)$ is an open-loop saddle point of Problem (MF-SDG)$_0$ for $x$. By the same argument as \rf{completion of square-closed},
noting \rf{range-eta-zero-sum-open-loop} and \rf{range-bar-eta-zero-sum-open-loop}, we obtain
\bel{completion of square-open}\ba{ll}
\ds J\big(x;\Th^{**}\big(X^{**}(\cd)-\dbE[X^{**}(\cd)]\big)+\bar\Th^{**}\dbE[X^{**}(\cd)]+v^{**}(\cd)\big)-\blan\h Px+2\bar\eta_o(0),x\bran\\
\ns\ds=\dbE\int_0^\i\Big[\blan\Si_o(v^{**}-\dbE[v^{**}]),v^{**}-\dbE[v^{**}]\bran+2\blan B^\top\eta_o+D^\top\z_o+D^\top P\si+\rho,v^{**}-\dbE[v^{**}]\bran\\
\ns\ds\qq\qq+\blan\bar\Si_o\dbE[v^{**}],\dbE[v^{**}]\bran+2\blan\h B^\top\bar\eta_o+\h D^\top P\dbE[\si]+\dbE[\z_o])+\dbE[\rho],\dbE[v^{**}]\bran\\
\ns\ds\qq\qq+\blan P\si,\si\bran+2\blan\eta_o,b-\dbE[b]\bran+2\blan\z_o,\si\bran+2\blan\bar\eta_o,\dbE[b]\bran\Big]dt\\
\ns\ds=\dbE\int_0^\i\Big[\blan P\si,\si\bran+2\blan\eta_o,b-\dbE[b]\bran+2\blan\z_o,\si\bran+2\blan\bar\eta_o,\dbE[b]\bran\\
\ns\ds\qq\qq-\big|(\Si_o^\dag)^{\frac{1}{2}}\big\{B^\top(\eta_o-\dbE[\eta_o])
+D^\top(\z_o-\dbE[\z_o])+D^\top P(\si-\dbE[\si])+\rho-\dbE[\rho]\big\}\big|^2\\
\ns\ds\qq\qq-\big|(\bar\Si_o^\dag)^{\frac{1}{2}}\big\{\h B^\top\bar\eta_o+\h D^\top\big(P\dbE[\si]+\dbE[\z_o]\big)+\dbE[\rho]\big\}\big|^2\Big]dt.\ea\ee
Since both $u^{**}(\cd)\equiv(u_1^{**}(\cd),u_2^{**}(\cd))$ and $u^*(\cd)\equiv(u_1^*(\cd),u_2^*(\cd))$ are open-loop saddle points of Problem (MF-SDG)$_0$ for $x$, we have
$$J(x;u_1^*(\cd),u_2^*(\cd))\les J(x;u_1^{**}(\cd),u_2^*(\cd))\les J(x;u_1^{**}(\cd),u_2^{**}(\cd))\les J(x;u_1^*(\cd),u_2^{**}(\cd))\les J(x;u_1^*(\cd),u_2^*(\cd)).$$
Therefore, $J(x;u^*(\cd))=J(x;u^{**}(\cd))$ for all $x$, which, together with \rf{completion of square-closed} and \rf{completion of square-open}, yields
$$P_c=P,\q\h P_c=\h P,\q\eta_c(\cd)=\eta_o(\cd),\q \z_c(\cd)=\z_o(\cd),\q \bar\eta_c(\cd)=\bar\eta_o(\cd).$$
Thus, the value function is given by \rf{V(x)-zero-sum}. The proof is complete.
\endpf

\ms

Finally, we have the following corollary for Problem (MF-SLQ), since it is a special case of Problem (MF-SDG)$_0$ when $m_2=0$.

\bc{C-4.6} \sl
For Problem {\rm(MF-SLQ)}, if the open-loop optimal control admits a closed-loop representation, then each open-loop optimal control must be an outcome of a closed-loop optimal strategy.
\ec

\section{Examples}

In this section, we present some examples illustrating the results in the previous sections.

\ms

The following example shows that for the closed-loop saddle points of the mean-field LQ zero-sum stochastic differential game, the system of generalized AREs may only admit non-static stabilizing solutions even if the system $[A,\bar A,C,\bar C;B,\bar B,D,\bar D]$ is MF-$L^2$-stabilizable, which leads to the non-existence of a closed-loop saddle point.

\bex{Example 5.1} \rm We consider one example for Problem (MF-SDG)$_0$. Consider the following one-dimensional state equation
\bel{e5.1-1}\left\{\2n\ba{ll}
\ds dX(t)=-\frac{1}{2}X(t)dt+\big[u_1(t)+u_2(t)\big]dW(t),\q t\ges0,\\
\ns\ds X(0)=x,\ea\right.\ee
with the cost functional
\bel{e5.1-2}\ba{ll}
\ds J(x;u_1(\cd),u_2(\cd))=\dbE\int_0^\i\Bigg[\Bigg\langle\2n
\begin{pmatrix}1&1&-1\\1&1&0\\ -1&0&-1\\\end{pmatrix}
\begin{pmatrix}X(t)\\u_1(t)\\u_2(t)\end{pmatrix},\begin{pmatrix}X(t)\\u_1(t)\\u_2(t)\end{pmatrix}
\2n\Bigg\rangle\\
\ns\ds\qq\qq\qq\qq\qq\qq+\Bigg\langle\2n
\begin{pmatrix}1&-1&1\\-1&0&0\\1&0&0\end{pmatrix}
\begin{pmatrix}\dbE[X(t)]\\ \dbE[u_1(t)]\\ \dbE[u_2(t)]\end{pmatrix},\begin{pmatrix}\dbE[X(t)]\\ \dbE[u_1(t)]\\ \dbE[u_2(t)]\end{pmatrix}
\2n\Bigg\rangle\Bigg]dt.\ea\ee
In this example,
$$\left\{\2n\ba{ll}
\ds A=-\frac{1}{2},\q \bar A=0,\q B=\bar B=(0,0), \q C=\bar C=0, \q D=(1,1),\q \bar D=(0,0),\\
\ns\ds Q=\bar Q=1,\q S=\begin{pmatrix}1\\ -1\end{pmatrix},\q \bar S=\begin{pmatrix}-1\\-1\end{pmatrix},\q R=\begin{pmatrix}1&0\\0&-1\end{pmatrix},\q
\bar R=\begin{pmatrix}0&0\\0&0\end{pmatrix},\ea\right.$$
and from \rf{hat}, we have
$$\h A=-\frac{1}{2}, \q \h B=(0,0),\q \h C=0,\q \h D=(1,1), \q  \h Q=2, \q \h S=\begin{pmatrix}0\\0\end{pmatrix},\q \h R=\begin{pmatrix}1&0\\0&-1\end{pmatrix}.$$
According to Lemma 2.3 in \cite{Sun-Yong-Zhang2016}, the system $[A,C;B,D]=[A,0,C,0;B,0,D,0]\equiv[A,\bar A,C,\bar C;B,\bar B,D,\bar D]$ above is (MF-)$L^2$-stabilizable, and
$\Th=(\Th_1,\Th_2)^\top\in\sS[A,C;B,D]$ if and only if there exists a $P_0>0$, such that
$$2(A+B\Th)P_0+(C+D\Th)^2P_0=-P_0+(\Th_1+\Th_2)^2P_0<0.$$
Due to the fact that $P_0>0$, hence, we obtain that
\bel{e5.1-3}
-1<\Th_1+\Th_2<1.
\ee
Note that $\Si_c=\bar \Si_c$ is invertible for all $(P_c,\h P_c)\in\dbS^1\times\dbS^1$ with
$$
\Sigma^{-1}_c=\bar{\Sigma}^{-1}_c=\begin{pmatrix}P_c+1&P_c\\P_c&P_c-1\end{pmatrix}^{-1}
=\begin{pmatrix}-P_c+1&P_c\\P_c&-P_c-1\end{pmatrix}.$$
Then the corresponding system of generalized AREs \rf{ARE-zero-sum-closed-loop} reads
\bel{e5.1-4}
\left\{\2n\ba{ll}
\ds 0=P_cA+A^\top P_c+C^\top P_cC+Q-(P_cB+C^\top P_cD+S^\top)\Sigma^\dag_c(B^\top P_c+D^\top P_cC+S)\\
\ns\ds\q =3P_c+1,\\
\ns\ds 0=\h P _c\h A +\h A^\top\h P_c+\h C^\top P_c\h C+\h Q-(\h P_c\h B+\h C^\top P_c\h D+\h S^\top)\bar\Si^\dag_c(\h B^\top\h P_c
 +\h D^\top P_c\h C+\h S)\\
\ns\ds\q =-\h P_c+2. \ea\right.\ee
Thus, $P_c=-\frac{1}{3},\h P_c=2$ and
\bel{e5.1-5}
\left\{\2n\ba{ll}
\ds R_{11}+D^\top_1P_c D_1=\h R_{11}+\h D^\top_1P_c\h D_1=\frac{2}{3}\ges0,\\
\ns\ds R_{22}+D^\top_2P_c D_2=\h R_{22}+\h D^\top_2P_c\h D_2=-\frac{4}{3}\les0.\ea\right.\ee
Also, the range condition
$$
\left\{\2n\ba{ll}
\ds\sR(B^\top P_c+D^\top P_cC+S)\subseteq\sR(\Si_c),\\
\ns\ds \sR(\h B^\top\h P_c+\h D^\top P_c\h C+\h S)\subseteq\sR(\bar \Si_c)\ea\right.
$$
hold automatically since $\Si_c$ and $\bar \Si_c$ are invertible. However, we have
\bel{e5.1-6}
(\Th_1,\Th_2)^\top=-\Si^\dag_c(B^\top P_c+D^\top P_cC+S)+(I-\Si^\dag_c\Si_c)\th=\left(-\frac{5}{3},-\frac{1}{3}\right)^\top,\qq \forall \th\in\dbR,\ee
which is not a stabilizer of the system $[A,C;B,D]$. Hence, by Theorem \ref{T-4.4}, Problem (MF-SDG)$_{0}$ does not admit a closed-loop saddle point. From this example, we see that generalized AREs \rf{ARE-zero-sum-closed-loop} may only admit non-static stabilizing solutions.

\ms

The following example tells us that for the closed-loop saddle points of the mean-field LQ zero-sum stochastic differential game, it may admit uncountably many closed-loop saddle points.

\ex

\bex{Example 5.2} \rm We consider one example for Problem (MF-SDG)$_0$. Consider the following one-dimensional state equation
\bel{e5.2-1}
\left\{\2n\ba{ll}
\ds dX(t)=-\Big[\frac{1}{4}X(t)+\frac{1}{2}u_2(t)\Big]dt+\big[-X(t)+u_1(t)\big]dW(t), \q t\ges 0,\\
\ns\ds X(0)=x,\ea\right.\ee
with the cost functional
\bel{e5.2-2}\ba{ll}
\ds J(x;u_1(\cd),u_2(\cd))=\dbE\int_0^\i\Bigg[\Bigg\langle\2n
\begin{pmatrix}\frac{1}{2}&-1&-\frac{1}{2}\\-1&1&0\\-\frac{1}{2}&0&0\\\end{pmatrix}
\begin{pmatrix}X(t)\\u_1(t)\\u_2(t)\end{pmatrix},\begin{pmatrix}X(t)\\u_1(t)\\u_2(t)\end{pmatrix}
\2n\Bigg\rangle\\
\ns\ds\qq\qq\qq\qq\qq\qq+\Bigg\langle\2n
\begin{pmatrix}\frac{1}{2}&0&\frac{1}{2}\\0&0&0\\\frac{1}{2}&0&-1\end{pmatrix}
\begin{pmatrix}\dbE[X(t)]\\ \dbE[u_1(t)]\\ \dbE[u_2(t)]\end{pmatrix},\begin{pmatrix}\dbE[X(t)]\\ \dbE[u_1(t)]\\ \dbE[u_2(t)]\end{pmatrix}
\2n\Bigg\rangle\Bigg]dt.\ea\ee
In this example,
$$\left\{\2n\ba{ll}
\ds A=-\frac{1}{4}, \q \bar A=0, \q B=\Big(0,-\frac{1}{2}\Big), \q \bar B=(0,0), \q C=-1, \q \bar C=0, \q D=(1,0),\q \bar D=(0,0),\\
\ns\ds Q=\bar Q=\frac{1}{2}, \q S=\begin{pmatrix}-1\\-\frac{1}{2}\end{pmatrix},\q \bar S=\begin{pmatrix}0\\\frac{1}{2}\end{pmatrix},\q
R=\begin{pmatrix}1&0\\0&0\end{pmatrix},\q \bar R=\begin{pmatrix}0&0\\0&-1\end{pmatrix}\ea\right.$$
and from \rf{hat}, we have
$$
\h A=-\frac{1}{4}, \q \h B=\Big(0,-\frac{1}{2}\Big),\q \h C=-1,\q \h D=(1,0), \q \h Q=1, \q \h S=\begin{pmatrix}-1\\0\end{pmatrix},\q \h R=\begin{pmatrix}1&0\\0&-1\end{pmatrix}.
$$
According to Lemma 2.3 in \cite{Sun-Yong-Zhang2016}, the system $[A,C;B,D]=[A,0,C,0;B,0,D,0]\equiv[A,\bar A,C,\bar C;B,\bar B,D,\bar D]$ above is (MF-)$L^2$-stabilizable, and
$\Th=(\Th_1,\Th_2)^\top\in \sS[A,C;B,D]$ if and only if there exists a $P_0>0$, such that
$$2(A+B\Th)P_0+(C+D\Th)^2P_0=2\Big(-\frac{1}{4}-\frac{1}{2}\Th_2\Big)P_0+(-1+\Th_1)^2P_0<0,$$
that is,
\bel{e5.2-3}
\Th^2_1-2\Th_1+\frac{1}{2}<\Th_2.\ee
Then the corresponding system of generalized AREs \rf{ARE-zero-sum-closed-loop} reads
\bel{e5.2-4}\ba{ll}
\ds 0=P_cA+A^\top P_c+C^\top P_cC+Q-(P_cB+C^\top P_cD+S^\top)\Si^\dag_c(B^\top P_c+D^\top P_cC+S)\\
\ns\ds =\frac{1}{2}P_c+\frac{1}{2}-\left(-P_c-1,-\frac{1}{2}P_c-\frac{1}{2}\right)
\begin{pmatrix}P_c+1&0\\0&0\end{pmatrix}^\dag\begin{pmatrix}-P_c-1\\-\frac{1}{2}P_c-\frac{1}{2}\end{pmatrix}\\
\ns\ds =\frac{1}{2}P_c+\frac{1}{2}-\frac{1}{4}(P_c+1)^2\begin{pmatrix}2,1\end{pmatrix}\begin{pmatrix}P_c+1&0\\0&0\end{pmatrix}^\dag\begin{pmatrix}2\\1\end{pmatrix},\ea\ee
which admits a unique solution $P_c=-1$; and
\bel{e5.2-5}\ba{ll}
\ds0=\h P_c\h A +\h A^\top\h P_c+\h C^\top P_c\h C+\h Q-(\h P_c\h B+\h C^\top P_c\h D+\h S^\top)\bar \Si^\dag_c(\h B^\top\h P_c+\h D^\top P_c\h C+\h S)\\
\ns\ds =-\frac{1}{2}\hat{P}_c+\frac{1}{4}\hat{P}^2_c,\ea\ee
which admits a solution $\h P_{c}=0$. Thus,
\bel{e5.2-6}\ba{ll}
\Si_c=\begin{pmatrix}0&0\\0&0\end{pmatrix}, \q
\bar \Si_c=\begin{pmatrix}0&0\\0&-1\end{pmatrix}, \q B^\top P_c+D^\top P_cC+S=\h B^\top\h P_c+\h D^\top P_c\h C+\h S=\begin{pmatrix}0\\0\end{pmatrix}.\ea\ee
Hence,
\bel{e5.2-7}\left\{\2n\ba{ll}
\ds R_{11}+D^\top_1P_c D_1=0\geq0, \qq \h R_{11}+\h D^\top_1P_c\h D_1=1\ges0,\\
\ns\ds R_{22}+D^\top_2P_c D_2=0\leq0, \qq \h R_{22}+\h D^\top_2P_c\h D_2=-1\les0,\ea\right.\ee
and the range condition
$$\left\{\2n\ba{ll}
\ds \sR(B^\top P_c+D^\top P_cC+S)\subseteq\sR(\Si_c),\\
\ns\ds \sR(\h B^\top\h P_c+\h D^\top P_c\h C+\h S)\subseteq\sR(\bar \Si_c)\ea\right.$$
hold. By Theorem \ref{T-4.4}, we see that Problem (MF-SDG)$_{0}$ admits a closed-loop saddle point if the condition \rf{e5.2-3} holds. But,
\bel{e5.2-8}
\Si^\dag_c(B^\top P_c+D^\top P_cC+S)=\bar \Si^\dag_c(\h B^\top\h P_c+\h D^\top P_c\h C+\h S)=(0,0)^\top\notin \sS[A,C;B,D].\ee
However, we can choose $\th,\bar \th\neq 0$ in \rf{closed-loop saddle point-3}, such that $\big((I-\Si^\dag_c\Si_c)\th,(I-\bar \Si^\dag_c\bar \Si_c)\bar \th\big)$ is a stabilizer of the system $[A,C;B,D]$. Thus, Problem (MF-SDG)$_{0}$ may still admit uncountably many closed-loop saddle points.

\ms

The following example shows that for the closed-loop saddle points of the mean-field LQ zero-sum stochastic differential game, it may happen that the system $[A,\bar A,C,\bar C;B,\bar B,D,\bar D]$ has more than one (uncountably many) MF-$L^2$-stabilizer, while the closed-loop saddle point is unique.

\ex

\bex{Example 5.3.} \rm We consider one example for Problem (MF-SDG)$_0$. Consider the following one-dimensional state equation
\bel{e5.3-1}
\left\{\2n\ba{ll}
\ds dX(t)=\big[-8X(t)+u_1(t)-u_2(t)\big]dt+\big[u_1(t)+u_2(t)\big]dW(t), \q  t\ges 0,\\
\ns\ds X(0)=x,
\ea\right.
\ee
with the cost functional
\bel{e5.3-2}\ba{ll}
\ds J(x;u_1(\cd),u_2(\cd))=\dbE\int_0^\i\Big[12X^2(t)+u^2_1(t)-u^2_2(t)-12|\dbE[X(t)]|^2+|\dbE[u_1(t)]|^2-|\dbE[u_2(t)]|^2]\Big]dt.\ea\ee
In this example,
$$\left\{\2n\ba{ll}
\ns\ds A=-8, \q \bar A=0,\q  B=(1,-1),\q \bar B=(0,0),\q C=\bar C=0,\q D=(1,1),\q \bar D=(0,0),\\
\ns\ds Q=12,\q \bar Q=52, \q S=\bar S=\begin{pmatrix}0\\0\end{pmatrix},\q R=\bar R=\begin{pmatrix}1&0\\0&-1\end{pmatrix},\ea\right.$$
and from \rf{hat}, we have
$$
\h A=-8, \q \h B=(1,-1),\q \h C=0,\q \h D=(1,1), \q  \h Q=64, \q \h S=\begin{pmatrix}0\\0\end{pmatrix},\q \h R=\begin{pmatrix}2&0\\0&-2\end{pmatrix}.
$$
According to Lemma 2.3 in \cite{Sun-Yong-Zhang2016}, the system $[A,C;B,D]=[A,0,C,0;B,0,D,0]\equiv[A,\bar A,C,\bar C;B,\bar B,D,\bar D]$ above is (MF-)$L^2$-stabilizable, and
$\Th=(\Th_1,\Th_2)^\top\in\sS[A,C;B,D]$ if and only if there exists a $P_0>0$, such that
$$
2(A+B\Th)P_0+(C+D\Th)^2P_0=\big[-16+2(\Th_1-\Th_2)+(\Th_1+\Th_2)^2\big]P_0<0,
$$
that is
\bel{e5.3-3}
-16+2(\Th_1-\Th_2)+(\Th_1+\Th_2)^2<0,
\ee
Then the corresponding generalized ARE \rf{ARE-zero-sum-closed-loop} for $P_c$ reads
\bel{e5.3-4}\ba{ll}
\ds 0=P_cA+A^\top P_c+C^\top P_cC+Q-(P_cB+C^\top P_cD+S^\top)\Si^\dag_c(B^\top P_c+D^\top P_cC+S)\\
\ns\ds\q =-16P_c+12-P^2_c(1,-1)\begin{pmatrix}P_c+1&P_c\\P_c&P_c-1\end{pmatrix}^\dag\begin{pmatrix}1\\-1\end{pmatrix}\\
\ns\ds\q =-16P_c+12-P^2_c(1,-1)\begin{pmatrix}-P_c+1&P_c\\P_c&-P_c-1\end{pmatrix}\begin{pmatrix}1\\-1\end{pmatrix}\\
\ns\ds\q =4P^3_c-16P_c+12,
\ea\ee
which has three solutions:
$$
P_c=1,\q P_c=\frac{-1+\sqrt{13}}{2},\q P_c=\frac{-1-\sqrt{13}}{2}.
$$
All of them satisfy the range condition
$$
\sR(B^\top P_c+D^\top P_cC+S)\subseteq\sR(\Si_c)
$$
since $\Si_c=R+D^\top P_cD$ is invertible for any $P_c\in \dbR$. However, only $P_c=1$ satisfies
\bel{e5.3-5}
R_{11}+D^\top_1P_c D_1=2\ges0,\q R_{22}+D^\top_2P_c D_2=0\les0.
\ee
And in this case, $\Si_c=R+D^\top P_cD=\begin{pmatrix}0&1\\1&-2\end{pmatrix}$.

Now putting $P_{c}=1$ into the following generalized ARE \rf{ARE-zero-sum-closed-loop} to resolve $\h P_c$:
\bel{e5.3-6}\ba{ll}
\ds 0=\h P_c\h A +\h A^\top\h P_c+\h C^\top P_c\h C+\h Q-(\h P_c\h B+\h C^\top P_c\h D+\h S^\top)\bar \Si^\dag_c(\h B^\top\h P_c+\h D^\top P_c\h C+\h S)\\
\ns\ds =-16\h P_c+64-\frac{1}{4}\h P_c^2(1,-1)\begin{pmatrix}1&1\\1&-3\end{pmatrix}\begin{pmatrix}1\\-1\end{pmatrix}\\
\ns\ds =-16\h P_c+64+\h P_c^2,
\ea\ee
which has a unique solution $\h P_c=8$. Due to the fact that $\bar \Si_c=\h R+\h D^\top P_c\h D=\begin{pmatrix}3&1\\1&-1\end{pmatrix}$ is invertible,
it is easy to verify that the range condition
$$\sR(\h B^\top\h P_c+\h D^\top P_c\h C+\h S)\subseteq\sR(\bar \Si_c)$$
holds, and
\bel{e5.3-7}
\h R_{11}+\h D^\top_1P_c\h D_1=3\ges0,\qq \h R_{22}+\h D^\top_2P_c\h D_2=-1\les0.
\ee
For any $\th \in \dbR$, from \rf{closed-loop saddle point-3}, we have
\bel{e5.3-8}\left\{\ba{ll}
\ds \Th^*=-\Si_c^\dag(B^\top P_c+D^\top P_cC+S)+\big(I-\Si_c^\dag\Si_c\big)\th=(1,-3)^\top,\\
\ns\ds \bar \Th^*=-\bar\Si_c^\dag\big[\h B^\top\Pi_c+\h D^\top P_c\h C+\h S\big]+\big(I-\bar\Si_c^\dag\bar\Si_c\big)\bar\th=(0,-1)^\top.\ea\right.
\ee
We can see that $(\Th^*_1,\Th^*_2)=(1,-3)$ satisfies \rf{e5.3-3} and hence is a stabilizer of the system $[A,C;B,D]$. By Theorem \ref{T-4.4}, the above problem admits a unique closed-loop saddle point.

On the other hand, by verifying \rf{e5.3-3}, we see that
$$
(\Th^*_1,\Th^*_2)=(1,-3), \q (\Th^*_1,\tilde{\Th}_2)=(1,2), \q (\tilde{\Th}_1,\tilde{\Th}_2)=(0,2), \q (\tilde{\Th}_1,\Th^*_2)=(0,-3)
$$
are stabilizers of $[A,B;C,D]$, but only $(\Th^{\ast}_{1},\Th^{\ast}_{2})$ is the closed-loop saddle point of the problem.

\ms

The following example shows that for the closed-loop saddle points of the mean-field LQ zero-sum stochastic differential game, it may happen that the system of generalized AREs may only admit non-static stabilizing solutions and the system $[A,\bar A,C,\bar C;B,\bar B,D,\bar D]$ is not MF-$L^2$-stabilizable. Thus no closed-loop saddle points exist.

\ex
\bex{Example 5.4} \rm We give one example of Problem (MF-SDG)$_0$. Consider the following two-dimensional state equation
\bel{e5.4-1}\left\{\2n\ba{ll}
\ds dX(t)=\big\{AX(t)+\bar A\dbE[X(t)]+Bu(t)+\bar B\dbE[u(t)]\big\}dt\\
\ns\ds\qq\qq\q+\big\{CX(t)+\bar C\dbE[X(t)]+Du(t)+\bar D\dbE[u(t)]\big\}dW(t),\q t\ges0,\\
\ns\ds X(0)=x,\ea\right.\ee
with the cost functional
\bel{e5.4-2}\ba{ll}
\ds J(x;u_1(\cd),u_2(\cd))=\dbE\int_0^\i\Bigg[\Bigg\langle\2n
\begin{pmatrix}Q&S^\top\\ S&R\end{pmatrix}
\begin{pmatrix}X(t)\\ u(t)\end{pmatrix},\begin{pmatrix}X(t)\\ u(t)\end{pmatrix}\2n\Bigg\rangle \\
\ns\ds\qq\qq\qq\qq\qq+\Bigg\langle\2n
\begin{pmatrix}\bar Q&\bar S^\top\\ \bar S&\bar R\\\end{pmatrix}
\begin{pmatrix}\dbE[X(t)]\\ \dbE[u(t)]\end{pmatrix},\begin{pmatrix}\dbE[X(t)]\\ \dbE[u(t)]\end{pmatrix}
\2n\Bigg\rangle\Bigg]dt.\ea\ee
In this example,
$$\left\{\2n\ba{ll}
\ds A=\begin{pmatrix}1&0.5\\0.5&1\\\end{pmatrix}, \q \bar A=\begin{pmatrix}1&0\\0&1\\\end{pmatrix},\q B=\begin{pmatrix}1&0\\0&1\\\end{pmatrix}, \q \bar B=\begin{pmatrix}2.5&0\\0&0.5\\\end{pmatrix}, \\
\ns\ds C=\begin{pmatrix}1&0.5\\0.5&1\\\end{pmatrix},\q \bar C=\begin{pmatrix}0&-0.5\\-0.5&0\\\end{pmatrix},\q D=\begin{pmatrix}1&0\\0&1\\\end{pmatrix}, \q \bar D=\begin{pmatrix}0.5&0\\0&0.5\\\end{pmatrix}, \\
\ns\ds Q=\begin{pmatrix}7.125&2.5\\2.5&1.5\\\end{pmatrix},\q \bar Q=\begin{pmatrix}5.875&-2.5\\-2.5&1.1\\\end{pmatrix},\q S=\begin{pmatrix}0&0\\0&0\end{pmatrix},\q \bar S=\begin{pmatrix}0&0\\0&0\end{pmatrix}, \\
\ns\ds R=\begin{pmatrix}2&0\\0&-1\end{pmatrix},\q \bar R=\begin{pmatrix}0.5&0\\0&-0.5\end{pmatrix},\ea\right.$$
and from \rf{hat}, we have
$$\left\{\2n\ba{ll}
\ds\h A=\begin{pmatrix}2&0.5\\0.5&2\\\end{pmatrix}, \q \h B=\begin{pmatrix}3.5&0\\0&1.5\\\end{pmatrix}, \q \h C=\begin{pmatrix}1&0\\0&1\\\end{pmatrix},\q
\h D=\begin{pmatrix}1.5&0\\0&1.5\\\end{pmatrix}, \\
\ns\ds \h Q=\begin{pmatrix}13&0\\0&2.6\\\end{pmatrix}, \q \h S=\begin{pmatrix}0&0\\0&0\end{pmatrix},\q \h R=\begin{pmatrix}2.5&0\\0&-1.5\end{pmatrix}.
\ea\right.$$
According to Proposition A.5 in \cite{Huang-Li-Yong2015}, the system $[A,\bar A,C,\bar C;B,\bar B,D,\bar D]$ above is MF-$L^2$-stabilizable, and
$(\Th,\bar\Th)\in\sS[A,\bar A,C,\bar C;B,\bar B,D,\bar D] $ if and only if there exist $P_0>0$ and $\bar P_0>0$, such that
\bel{e5.4-3}\left\{\2n\ba{ll}
\ds (A+B\Th)P_0+P_0(A+B\Th)^\top+(C+D\Th)P_0(C+D\Th)^\top+(\h C+\h D\bar\Th)\bar P_0(\h C+\h D\bar\Th)^\top\\
\ns\ds =\left[\begin{pmatrix}1&0.5\\0.5&1\\\end{pmatrix}+\Th\right]P_0+P_0\left[\begin{pmatrix}1&0.5\\0.5&1\\\end{pmatrix}+\Th\right]^\top
+\left[\begin{pmatrix}1&0.5\\0.5&1\\\end{pmatrix}+\Th\right]P_0\left[\begin{pmatrix}1&0.5\\0.5&1\\\end{pmatrix}+\Th\right]^\top\\
\ns\ds \q+\left[\begin{pmatrix}1&0\\0&1\\\end{pmatrix}+\begin{pmatrix}1.5&0\\0&1.5\\\end{pmatrix}\bar\Th\right]\bar P_0
\left[\begin{pmatrix}1&0\\0&1\\\end{pmatrix}+\begin{pmatrix}1.5&0\\0&1.5\\\end{pmatrix}\bar\Th\right]^\top<0, \\
\ns\ds (\h A+\h B\bar\Th)\bar P_0+\bar P_0(\h A+\h B\bar\Th)^\top\\
\ns\ds =\left[\begin{pmatrix}2&0.5\\0.5&2\\\end{pmatrix}+\begin{pmatrix}3.5&0\\0&1.5\\\end{pmatrix}\bar\Th\right]\bar P_0
+\bar P_0\left[\begin{pmatrix}2&0.5\\0.5&2\\\end{pmatrix}+\begin{pmatrix}3.5&0\\0&1.5\\\end{pmatrix}\bar\Th\right]^\top<0.
\ea\right.\ee
Note that the corresponding system of generalized AREs \rf{ARE-zero-sum-closed-loop} reads
\bel{e5.4-4}\left\{\2n\ba{ll}
\ds 0=P_cA+A^\top P_c+C^\top P_cC+Q-(P_cB+C^\top P_cD+S^\top)\Sigma^\dag_c(B^\top P_c+D^\top P_cC+S), \\
\ns\ds 0=\h P _c\h A +\h A^\top\h P_c+\h C^\top P_c\h C+\h Q-(\h P_c\h B+\h C^\top P_c\h D+\h S^\top)\bar\Si^\dag_c(\h B^\top\h P_c+\h D^\top P_c\h C+\h S),
\ea\right.\ee
with
$$\left\{\ba{ll}
\ds \Si_c=R+D^\top P_cD=\begin{pmatrix}2&0\\0&-1\\\end{pmatrix}+P_c,\\
\ns\ds \bar\Si_c=\h R+\h D^\top P_c\h D=\begin{pmatrix}2.5&0\\0&-1.5\\\end{pmatrix}+2.25P_c.
\ea\right.$$
Then, solving (\ref{e5.4-4}) yields
$$
P_c=\begin{pmatrix}-1&0\\0&-1\\\end{pmatrix}, \qq \h P_c=\begin{pmatrix}1&0\\0&-1\\\end{pmatrix}.
$$
Thus,
\bel{e5.4-5}\left\{\2n\ba{ll}
\ds R_{11}+D^\top_1P_c D_1=1\ges0, & \h R_{11}+\h D^\top_1P_c\h D_1=0.25\ges0,\\
\ns\ds R_{22}+D^\top_2P_c D_2=-2\les0, & \h R_{22}+\h D^\top_2P_c\h D_2=-3.75\les0.\ea\right.\ee
Also, the range condition
$$
\left\{\2n\ba{ll}
\ds\sR(B^\top P_c+D^\top P_cC+S)\subseteq\sR(\Si_c),\\
\ns\ds \sR(\h B^\top\h P_c+\h D^\top P_c\h C+\h S)\subseteq\sR(\bar \Si_c)\ea\right.
$$
hold automatically since $\Si_c$ and $\bar \Si_c$ are invertible. However, we have
\bel{e5.4-6}\left\{\2n\ba{ll}
\ds \Th^*=-\Si^\dag_c(B^\top P_c+D^\top P_cC+S)+(I-\Si^\dag_c\Si_c)\th=\begin{pmatrix}2&0.5\\-0.25&-1\end{pmatrix}, \qq \forall\ \th\in\dbR, \\
\ns\ds \bar\Th^*=-\bar\Si^\dag_c(\h B^\top\h P_c+\h D^\top P_c\h C+\h S)+(I-\bar\Si^\dag_c\bar\Si_c)\th=\begin{pmatrix}-8&0\\0&-0.8\end{pmatrix}, \qq \forall\ \th\in\dbR.
\ea\right.\ee
However, substituting the above $\Th^*$ and $\bar\Th^*$ of (\ref{e5.4-6}) into (\ref{e5.4-3}), we no longer obtain the result of $P_0>0$ and $\bar P_0>0$, which conflicts with (\ref{e5.4-3}). Therefore, $(\Th^*,\bar\Th^*)$ is not an MF-$L^2$-stabilizer of the system $[A,\bar A,C,\bar C;B,\bar B,D,\bar D]$. Hence, by Theorem \ref{T-4.4}, Problem (MF-SDG)$_{0}$ does not admit a closed-loop saddle point. From this example, we see that the system of generalized AREs \rf{ARE-zero-sum-closed-loop} may only admit non-static stabilizing solutions.
%
\ms

The following example shows that for the mean-field LQ zero-sum stochastic differential game, it may happen that the closed-loop saddle point uniquely exists, which coincides with the closed-loop representation of the open-loop saddle point. Moreover, the system of generalized AREs may only admit static stabilizing solutions and the system $[A,\bar A,C,\bar C;B,\bar B,D,\bar D]$ is MF-$L^2$-stabilizable.

\ex

\ms

The following example shows that for the mean-field LQ zero-sum stochastic differential game, it may happen that the closed-loop saddle point uniquely exists, which coincides with the closed-loop representation of the open-loop saddle point. Moreover, the system of generalized AREs admits static stabilizing solutions and the system $[A,\bar A,C,\bar C;B,\bar B,D,\bar D]$ is MF-$L^2$-stabilizable.

\bex{Example 5.5} \rm Consider the following two-dimensional state equation
\bel{e5.5-1}\left\{\2n\ba{ll}
\ds dX(t)=\big\{AX(t)+\bar A\dbE[X(t)]+Bu(t)+\bar B\dbE[u(t)]\big\}dt\\
\ns\ds\qq\qq\q+\big\{CX(t)+\bar C\dbE[X(t)]+Du(t)+\bar D\dbE[u(t)]\big\}dW(t),\q t\ges0,\\
\ns\ds X(0)=x,\ea\right.\ee
with the cost functional
\bel{e5.5-2}\ba{ll}
\ds J(x;u_1(\cd),u_2(\cd))\1n=\1n\dbE\2n\int_0^\i\2n\bigg[\bigg\langle\2n
\begin{pmatrix}Q&S^\top\\ S&R\end{pmatrix}
\begin{pmatrix}X(t)\\ u(t)\end{pmatrix},\begin{pmatrix}X(t)\\ u(t)\end{pmatrix}\2n\bigg\rangle+\bigg\langle\2n
\begin{pmatrix}\bar Q&\bar S^\top\\ \bar S&\bar R\\\end{pmatrix}
\begin{pmatrix}\dbE[X(t)]\\ \dbE[u(t)]\end{pmatrix},\begin{pmatrix}\dbE[X(t)]\\ \dbE[u(t)]\end{pmatrix}
\2n\bigg\rangle\bigg]dt.\ea\ee
Let
$$\left\{\2n\ba{ll}
\ds A=\begin{pmatrix}-1&-1\\0&-3\\\end{pmatrix}, \q \bar A=\begin{pmatrix}-1&0\\-1&-20\\\end{pmatrix},\q B=\begin{pmatrix}1&0\\0&1\\\end{pmatrix}, \q \bar B=\begin{pmatrix}2.5&0\\0&0.5\\\end{pmatrix}, \\
\ns\ds C=\begin{pmatrix}1&0.5\\0.5&1\\\end{pmatrix},\q \bar C=\begin{pmatrix}0&-0.5\\-0.5&0\\\end{pmatrix},\q D=\begin{pmatrix}1&0\\0&1\\\end{pmatrix}, \q \bar D=\begin{pmatrix}0.5&0\\0&0.5\\\end{pmatrix}, \\
\ns\ds Q=\begin{pmatrix}2.3070&0.7781\\0.7781&0.3123\\\end{pmatrix},\q \bar Q=\begin{pmatrix}5.6930&0.7219\\0.7219&22.2317\\\end{pmatrix},\q S=\begin{pmatrix}0&0\\0&0\end{pmatrix},\q \bar S=\begin{pmatrix}0&0\\0&0\end{pmatrix}, \\
\ns\ds R=\begin{pmatrix}2&0\\0&-2\end{pmatrix},\q \bar R=\begin{pmatrix}0.75&0\\0&-0.5\end{pmatrix},\ea\right.$$
and from \rf{hat}, we have
$$\left\{\2n\ba{ll}
\ds \h A=\begin{pmatrix}-2&-1\\-1&-23\\\end{pmatrix}, \q \h B=\begin{pmatrix}3.5&0\\0&1.5\\\end{pmatrix}, \q \h C=\begin{pmatrix}1&0\\0&1\\\end{pmatrix},\q \h D=\begin{pmatrix}1.5&0\\0&1.5\\\end{pmatrix}, \\
\ns\ds \h Q=\begin{pmatrix}8&1.5\\1.5&22.5440\\\end{pmatrix}, \q \h S=\begin{pmatrix}0&0\\0&0\end{pmatrix},\q \h R=\begin{pmatrix}2.75&0\\0&-2.5\end{pmatrix}.
\ea\right.$$

First, we look at closed-loop representation of open-loop saddle point. According to \rf{ARE-zero-sum-open-loop}, we have the following coupled AREs:
\bel{e5.5-ARE-zero-sum-open-loop}\left\{\2n\ba{ll}
\ds PA+A^\top P+C^\top P C+Q-\big(PB+C^\top P D+S^\top\big)\Si_o^\dag(B^\top P+D^\top PC+S)=0,\\
\ns\ds\h P\h A+\h A^{\,\top}\h P+\h C^\top P\h C+\h Q-(\h P\h B+\h C^\top P\h D+\h S^\top)\bar\Si_o^\dag(\h B^\top\h P+\h D^\top P\h C+\h S)=0,
\ea\right.\ee
where
$$\Si_o\equiv R+D^\top PD=\begin{pmatrix}2&0\\0&-2\\\end{pmatrix}+P,\qq
\bar\Si_o\equiv\h R+\h D^\top P\h D=\begin{pmatrix}2.75&0\\0&-2.5\\\end{pmatrix}+2.25P.$$
Then, solving \rf{e5.5-ARE-zero-sum-open-loop} yields
$$
P=\begin{pmatrix}1&0\\0&0.1\\\end{pmatrix}, \qq \h P=\begin{pmatrix}1&0\\0&0.5\\\end{pmatrix}.
$$
The above solutions satisfy \rf{AE-zero-sum-open-loop} and \rf{Th-Th*-zero-sum-open-loop}. Further, since $\Si_o$ and $\bar \Si_o$ are invertible, we have
\bel{e5.5-solvable condition for AE-zero-sum-open-loop}\left\{\2n\ba{ll}
\ds\sR\big(B^\top P+D^\top PC+S\big)\subseteq\sR(\Si_o), \qq \sR\big(\h B^\top\h P+\h D^\top P\h C+\h S\,\big)\subseteq\sR(\bar\Si_o),\\
\ns\ds\Th^*=-\Si_o^{-1}(B^\top P+D^\top PC+S)=\begin{pmatrix}-0.6667&-0.1667\\ \ns 0.0263&0.1053\end{pmatrix},\\
\ns\ds\bar\Th^*=-\bar\Si_o^{-1}\big(\h B^\top\h P+\h D^\top P\h C+\h S\,\big)=\begin{pmatrix}-1&0\\0&0.3956\end{pmatrix}.
\ea\right.\ee
According to Proposition A.5 in \cite{Huang-Li-Yong2015}, the system $[A,\bar A,C,\bar C;$ $B,\bar B,D,\bar D]$ above is MF-$L^2$-stabilizable, and
$(\Th,\bar\Th)\in\sS[A,\bar A,C,\bar C;B,\bar B,D,\bar D] $ if and only if there exist $P_0>0$ and $\bar P_0>0$, such that
\bel{e5.5-static-stabilizing-judgement-condition}\left\{\2n\ba{ll}
\ds (A+B\Th)P_0+P_0(A+B\Th)^\top+(C+D\Th)P_0(C+D\Th)^\top+(\h C+\h D\bar\Th)\bar P_0(\h C+\h D\bar\Th)^\top\\
\ns\ds =\left[\begin{pmatrix}-1&-1\\0&-3\\\end{pmatrix}+\Th\right]P_0+P_0\left[\begin{pmatrix}-1&-1\\0&-3\\\end{pmatrix}+\Th\right]^\top
+\left[\begin{pmatrix}1&0.5\\0.5&1\\\end{pmatrix}+\Th\right]P_0\left[\begin{pmatrix}1&0.5\\0.5&1\\\end{pmatrix}+\Th\right]^\top\\
\ns\ds \q+\left[\begin{pmatrix}1&0\\0&1\\\end{pmatrix}+\begin{pmatrix}1.5&0\\0&1.5\\\end{pmatrix}\bar\Th\right]\bar P_0
\left[\begin{pmatrix}1&0\\0&1\\\end{pmatrix}+\begin{pmatrix}1.5&0\\0&1.5\\\end{pmatrix}\bar\Th\right]^\top<0, \\
\ns\ds (\h A+\h B\bar\Th)\bar P_0+\bar P_0(\h A+\h B\bar\Th)^\top\\
\ns\ds =\left[\begin{pmatrix}-2&-1\\-1&-23\\\end{pmatrix}+\begin{pmatrix}3.5&0\\0&1.5\\\end{pmatrix}\bar\Th\right]\bar P_0
+\bar P_0\left[\begin{pmatrix}-2&-1\\-1&-23\\\end{pmatrix}+\begin{pmatrix}3.5&0\\0&1.5\\\end{pmatrix}\bar\Th\right]^\top<0.
\ea\right.\ee
Substituting the above $\Th^*$ and $\bar\Th^*$ of (\ref{e5.5-solvable condition for AE-zero-sum-open-loop}) into (\ref{e5.5-static-stabilizing-judgement-condition}), we have
$$P_0=\begin{pmatrix}0.5742&-0.0424\\-0.0424&0.4258\end{pmatrix}>0, \qq \bar P_0=\begin{pmatrix}0.8185&-0.0377\\-0.0377&0.1815\end{pmatrix}>0,$$
such that the negative definiteness conditions of (\ref{e5.5-static-stabilizing-judgement-condition}) hold. Then $(\Th^*,\bar\Th^*)$ above is an MF-$L^2$-stabilizer of the system $[A,\bar A,C,\bar C;B,\bar B,D,\bar D]$.
Now, we solve corresponding systems \rf{BSDE-eta-zero-sum-open-loop}--\rf{range-bar-eta-zero-sum-open-loop} to get $(\eta_o(\cd),\z_o(\cd))$ and $\bar\eta_o(\cd)$.
From \rf{closes-loop representation of an open-loop saddle point}, the closed-loop representation of this open-loop saddle point can be given by
\bel{e5.5-closes-loop representation of an open-loop saddle point}\ba{ll}
\ds u^*(\cd)=-\Si_o^{-1}(B^\top P+D^\top PC+S)\big(X^*(\cd)-\dbE[X^*(\cd)]\big)-\bar\Si_o^{-1}\big(\h B^\top\h P+\h D^\top P\h C+\h S\,\big)\dbE[X^*(\cd)]\\
\ns\ds\qq\q\ -\Si_o^{-1}\big\{B^\top(\eta_o(\cd)-\dbE[\eta_o(\cd)])+D^\top\big(\z_o(\cd)-\dbE[\z_o(\cd)]\big)\big\}
-\bar\Si_o^{-1}\big\{\h B^\top\bar\eta_o(\cd)+\h D^\top\dbE[\z_o(\cd)]\big\},\ea\ee
where $X^*(\cd)\in\sX[0,\i)$ is the solution to \rf{e5.1-1} corresponding to $u^*(\cd)$.

\ms

Next, we consider the closed-loop saddle point. Note that the corresponding system of generalized AREs \rf{ARE-zero-sum-closed-loop} reads
\bel{e5.5-ARE-zero-sum-closed-loop}\left\{\2n\ba{ll}
\ds 0=P_cA+A^\top P_c+C^\top P_cC+Q-(P_cB+C^\top P_cD+S^\top)\Sigma^\dag_c(B^\top P_c+D^\top P_cC+S), \\
\ns\ds 0=\h P _c\h A +\h A^\top\h P_c+\h C^\top P_c\h C+\h Q-(\h P_c\h B+\h C^\top P_c\h D+\h S^\top)\bar\Si^\dag_c(\h B^\top\h P_c+\h D^\top P_c\h C+\h S), \ea\right.\ee
with
$$\left\{\ba{ll}
\ds \Si_c\equiv R+D^\top P_cD=\begin{pmatrix}2&0\\0&-2\\\end{pmatrix}+P_c,\\
\ns\ds \bar\Si_c\equiv\h R+\h D^\top P_c\h D=\begin{pmatrix}2.75&0\\0&-2.5\\\end{pmatrix}+2.25P_c.
\ea\right.$$
Then, solving (\ref{e5.5-ARE-zero-sum-closed-loop}) yields
$$
P_c=\begin{pmatrix}1&0\\0&0.1\\\end{pmatrix}, \qq \h P_c=\begin{pmatrix}1&0\\0&0.5\\\end{pmatrix}.
$$
Thus,
\bel{e5.5-R11-R22}\left\{\2n\ba{ll}
\ds R_{11}+D^\top_1P_c D_1=3\ges0, & \h R_{11}+\h D^\top_1P_c\h D_1=5\ges0,\\
\ns\ds R_{22}+D^\top_2P_c D_2=-1.9\les0, & \h R_{22}+\h D^\top_2P_c\h D_2=-2.275\les0.\ea\right.\ee
Also, the range condition
$$\sR(B^\top P_c+D^\top P_cC+S)\subseteq\sR(\Si_c),\qq\sR(\h B^\top\h P_c+\h D^\top P_c\h C+\h S)\subseteq\sR(\bar \Si_c)$$
holds automatically since $\Si_c$ and $\bar \Si_c$ are invertible. We see that
\bel{e5.5-closed-loop saddle point-3}\left\{\2n\ba{ll}
\ds \Th^*=-\Si^{-1}_c(B^\top P_c+D^\top P_cC+S)=\begin{pmatrix}-0.6667&-0.1667\\ \ns 0.0263&0.1053\end{pmatrix},\\
\ns\ds \bar\Th^*=-\bar\Si^{-1}_c(\h B^\top\h P_c+\h D^\top P_c\h C+\h S)=\begin{pmatrix}-1&0\\0&0.3956\end{pmatrix},
\ea\right.\ee
which coincides with (\ref{e5.5-solvable condition for AE-zero-sum-open-loop}). Hence, by Theorem \ref{T-4.4}, the problem admits a closed-loop saddle point.
Solving corresponding backward systems \rf{BSDE-eta-zero-sum-closed-loop}--\rf{BODE-bar-eta-zero-sum-closed-loop}, we can get $(\eta_c(\cd),\z_c(\cd))$ and $\bar\eta_c(\cd)$, the closed-loop saddle point is given by \rf{e5.5-closed-loop saddle point-3} and
\bel{e5.5-closed-loop saddle point-4}\ba{ll}
\ds v^*(\cd)=-\Si_c^{-1}\big\{B^\top\big(\eta_c(\cd)-\dbE[\eta_c(\cd)]\big)
      +D^\top\big(\z_c(\cd)-\dbE[\z_c(\cd)]\big)\big\}-\bar\Si_c^{-1}\big\{\h B^\top\bar\eta_c(\cd)+\h D^\top\dbE[\z_c(\cd)]\big\}.\ea\ee
Since $P=P_c$ and $\h P=\h P_c$, we get $\Si=\Si_c,\bar\Si=\bar\Si_c$ and $\eta_o(\cd)=\eta_c(\cd)$, $\z_o(\cd)=\z_c(\cd)$, $\bar\eta_o(\cd)=\bar\eta_c(\cd)$.
This implies that the closed-loop representation of the open-loop saddle point coincides with the closed-loop saddle point.

\ms

The following example shows that for the mean-field LQ non-zero sum stochastic differential game, it may happen that the closed-loop representations of open-loop Nash equilibria are different from the closed-loop Nash equilibria. However, the solutions to the system of algebraic Riccati equations for the closed-loop representation of open-loop Nash equilibrium are symmetric.

\ex
\bex{Example 5.6} \rm Consider the following two-dimensional state equation
\bel{e5.6-1}\left\{\2n\ba{ll}
\ds dX(t)=\big\{AX(t)+\bar A\dbE[X(t)]+Bu(t)+\bar B\dbE[u(t)]\big\}dt\\
\ns\ds\qq\qq\q+\big\{CX(t)+\bar C\dbE[X(t)]+Du(t)+\bar D\dbE[u(t)]\big\}dW(t),\q t\ges0,\\
\ns\ds X(0)=x,\ea\right.\ee
with the cost functional
\bel{e5.6-2}\left\{\2n\ba{ll}
\ds J_1(x;u_1(\cd),u_2(\cd))=\dbE\int_0^\i\bigg[\bigg\langle\2n
\begin{pmatrix}Q_1&S_1^\top\\ S_1&R_1\end{pmatrix}
\begin{pmatrix}X(t)\\ u(t)\end{pmatrix},\begin{pmatrix}X(t)\\ u(t)\end{pmatrix}\2n\bigg\rangle \\
\ns\ds\qq\qq\qq\qq\qq\qq+\bigg\langle\2n
\begin{pmatrix}\bar Q_1&\bar S_1^\top\\ \bar S_1&\bar R_1\\\end{pmatrix}
\begin{pmatrix}\dbE[X(t)]\\ \dbE[u(t)]\end{pmatrix},\begin{pmatrix}\dbE[X(t)]\\ \dbE[u(t)]\end{pmatrix}
\2n\bigg\rangle\bigg]dt, \\
\ns\ds J_2(x;u_1(\cd),u_2(\cd))=\dbE\int_0^\i\bigg[\bigg\langle\2n
\begin{pmatrix}Q_2&S_2^\top\\ S_2&R_2\end{pmatrix}
\begin{pmatrix}X(t)\\ u(t)\end{pmatrix},\begin{pmatrix}X(t)\\ u(t)\end{pmatrix}\2n\bigg\rangle \\
\ns\ds\qq\qq\qq\qq\qq\qq+\bigg\langle\2n
\begin{pmatrix}\bar Q_2&\bar S_2^\top\\ \bar S_2&\bar R_2\\\end{pmatrix}
\begin{pmatrix}\dbE[X(t)]\\ \dbE[u(t)]\end{pmatrix},\begin{pmatrix}\dbE[X(t)]\\ \dbE[u(t)]\end{pmatrix}
\2n\bigg\rangle\bigg]dt.\ea\right.\ee
We let
$$\left\{\2n\ba{ll}
\ds A=\begin{pmatrix}-1&-1\\0&-1\\\end{pmatrix}, \q \bar A=\begin{pmatrix}-1&0\\-1&-1\\\end{pmatrix},\q B=\begin{pmatrix}1&0\\0&1\\\end{pmatrix}, \q \bar B=\begin{pmatrix}2.5&0\\0&0.5\\\end{pmatrix}, \\
\ns\ds C=\begin{pmatrix}1&0.5\\0.5&1\\\end{pmatrix},\q \bar C=\begin{pmatrix}0&-0.5\\-0.5&0\\\end{pmatrix},\q D=\begin{pmatrix}1&0\\0&1\\\end{pmatrix}, \q \bar D=\begin{pmatrix}0.5&0\\0&0.5\\\end{pmatrix}, \\
\ns\ds Q_1=\begin{pmatrix}2.85&0.9\\0.9&2.475\\\end{pmatrix},\q Q_2=\begin{pmatrix}1.35&0.4\\0.4&2.5375\\\end{pmatrix},\q \bar Q_1=\begin{pmatrix}6.0324&0.6\\0.6&0.325\\\end{pmatrix}, \\
\ns\ds \bar Q_2=\begin{pmatrix}7.15&1.6\\1.6&2.8625\\\end{pmatrix}, \q S_1=S_2=\begin{pmatrix}0&0\\0&0\end{pmatrix},\q \bar S_1=\bar S_2=\begin{pmatrix}0&0\\0&0\end{pmatrix}, \\
\ns\ds R_1=\begin{pmatrix}1&0\\0&1\end{pmatrix}, \q R_2=\begin{pmatrix}1&0\\0&1.5\end{pmatrix},\q \bar R_1=\begin{pmatrix}1&0\\0&1\end{pmatrix}, \q \bar R_2=\begin{pmatrix}0.5&0\\0&0\end{pmatrix},
\ea\right.$$
and from \rf{hat}, we have
$$\left\{\2n\ba{ll}
\ds\h A=\begin{pmatrix}-2&-1\\-1&-2\\\end{pmatrix}, \q \h B=\begin{pmatrix}3.5&0\\0&1.5\\\end{pmatrix}, \q \h C=\begin{pmatrix}1&0\\0&1\\\end{pmatrix},\q
\h D=\begin{pmatrix}1.5&0\\0&1.5\\\end{pmatrix}, \q \h Q_1=\begin{pmatrix}8.8824&1.5\\1.5&2.8\\\end{pmatrix}, \\
\ns\ds \h Q_2=\begin{pmatrix}8.5&2\\2&5.4\\\end{pmatrix}, \q \h S_1 = \h S_2=\begin{pmatrix}0&0\\0&0\end{pmatrix},\q
\h R_1=\begin{pmatrix}2&0\\0&2\end{pmatrix},\q \h R_2=\begin{pmatrix}1.5&0\\0&1.5\end{pmatrix}.
\ea\right.$$
To look at the closed-loop representation of open-loop Nash equilibria, we solve the corresponding (\ref{AE-LQ-game-component})--(\ref{Th-Th*-component}) to get
$$
P_1=\begin{pmatrix}1&0\\0&1\\\end{pmatrix}, \qq P_2=\begin{pmatrix}0.5&0\\0&1\\\end{pmatrix}, \qq
\h P_1=\begin{pmatrix}1&0\\0&0.5\\\end{pmatrix}, \qq \h P_2=\begin{pmatrix}1&0\\0&1\\\end{pmatrix},
$$
which are symmetric. Since
$$\begin{pmatrix}R_{111}+D_1^\top P_1D_1&R_{112}+D_1^\top P_1D_2\\R_{221}+D_2^\top P_2D_1&R_{222}
+D_2^\top P_2D_2\end{pmatrix}=\begin{pmatrix}2&0\\0&2.5\end{pmatrix}$$
and
$$\begin{pmatrix}\h R_{111}+\h D_1^\top P_1\h D_1&\h R_{112}+\h D_1^\top P_1\h D_2\\\h R_{221}+\h D_2^\top P_2\h D_1&\h R_{222}+\h D_2^\top P_2\h D_2\end{pmatrix}
=\begin{pmatrix}4.25&0\\0&3.75\end{pmatrix}$$
are invertible, it follows that 
\bel{e5.6-open-loop-Th-Th*}\ba{ll}
\ds\Th^{**}=\begin{pmatrix}-1&-0.25\\
\ns -0.2&-0.8\end{pmatrix}, \q
\bar\Th^{**}=\begin{pmatrix}-1.1765&0\\0&-0.8\end{pmatrix}.
\ea\ee
Also, it follows from \rf{BSDE-eta-game-BJ}--\rf{BODE-bar-eta-BJ}, we can get $(\eta(\cd),\z(\cd))$ and $\bar\eta(\cd)$. 
Then, from (\ref{Th-Th*-v*}), we have
\bel{e5.6-v*}\ba{ll}
\ds v^{**}(\cd)=-\BSi^{-1}\BJ^\top\1n\big\{\BB^\top\big(\eta(\cd)\1n-\1n\dbE[\eta(\cd)]\big)\1n
+\1n\BD^\top\big(\z(\cd)\1n-\1n\dbE[\z(\cd)]\big)\big\}-\bar\BSi^{-1}\BJ^\top\big\{\h\BB^\top\bar\eta(\cd)+\h\BD^\top\dbE[\z(\cd)]\big\}.\ea
\ee
By Theorem \ref{T-3.3}, making use of \rf{Th-Th*-v*}--\rf{BODE-bar-eta-BJ}, the problem admits an open-loop Nash equilibrium $u^*(\cd)$ for any initial state $x\in\dbR^n$. And it has the closed-loop representation \rf{closed-loop control} with $(\Th^{**},\bar\Th^{**},v^{**}(\cd))$ given by \rf{e5.6-open-loop-Th-Th*} and \rf{e5.6-v*}.

\ms

Next, we consider the closed-loop Nash equilibria.
According to Proposition A.5 in \cite{Huang-Li-Yong2015}, system $[A,\bar A,C,\\\bar C;B,\bar B,D,\bar D]$ above is MF-$L^2$-stabilizable, and
$(\Th,\bar\Th)\in\sS[A,\bar A,C,\bar C;B,\bar B,D,\bar D]$ if and only if there exist $P_0>0$ and $\bar P_0>0$, such that
\bel{e5.6-static-stabilizing-judgement-condition}\left\{\2n\ba{ll}
\ds (A+B\Th)P_0+P_0(A+B\Th)^\top+(C+D\Th)P_0(C+D\Th)^\top+(\h C+\h D\bar\Th)\bar P_0(\h C+\h D\bar\Th)^\top\\
\ns\ds =\left[\begin{pmatrix}-1&-1\\0&-1\\\end{pmatrix}+\Th\right]P_0+P_0\left[\begin{pmatrix}-1&-1\\0&-1\\\end{pmatrix}+\Th\right]^\top
+\left[\begin{pmatrix}1&0.5\\0.5&1\\\end{pmatrix}+\Th\right]P_0\left[\begin{pmatrix}1&0.5\\0.5&1\\\end{pmatrix}+\Th\right]^\top\\
\ns\ds \q+\left[\begin{pmatrix}1&0\\0&1\\\end{pmatrix}+\begin{pmatrix}1.5&0\\0&1.5\\\end{pmatrix}\bar\Th\right]\bar P_0
\left[\begin{pmatrix}1&0\\0&1\\\end{pmatrix}+\begin{pmatrix}1.5&0\\0&1.5\\\end{pmatrix}\bar\Th\right]^\top<0, \\
\ns\ds (\h A+\h B\bar\Th)\bar P_0+\bar P_0(\h A+\h B\bar\Th)^\top\\
\ns\ds =\left[\begin{pmatrix}-2&-1\\-1&-2\\\end{pmatrix}+\begin{pmatrix}3.5&0\\0&1.5\\\end{pmatrix}\bar\Th\right]\bar P_0
+\bar P_0\left[\begin{pmatrix}-2&-1\\-1&-2\\\end{pmatrix}+\begin{pmatrix}3.5&0\\0&1.5\\\end{pmatrix}\bar\Th\right]^\top<0.
\ea\right.\ee

For $i=1,2$, we solve the corresponding (\ref{AE-LQ-game-closed})--(\ref{Sigma}) to get
$$\left\{\2n\ba{ll}
\ds P_1=\begin{pmatrix}0.9949&-0.0168\\-0.0168&0.9201\\\end{pmatrix}, \qq P_2=\begin{pmatrix}0.6255&-0.0104\\-0.0104&1.01741\\\end{pmatrix}, \\
\ns\ds\h P_1=\begin{pmatrix}1.0023&-0.0155\\-0.0155&0.6472\\\end{pmatrix}, \qq \h P_2=\begin{pmatrix}0.8919&0.0126\\0.0126&0.9964\\\end{pmatrix},
\ea\right.$$
$$\left\{\2n\ba{ll}
\Si_1\deq R_{111}+D_1^\top P_1D_1=1.9949, \qq \Si_2\deq R_{222}+D_2^\top P_2D_2=2.5174, \\
\bar\Si_1\deq \h R_{111}+\h D_1^\top P_1\h D_1=4.2386, \qq \bar\Si_2\deq \h R_{222}+\h D_2^\top P_2\h D_2=3.7891,
\ea\right.$$
and
\bel{e5.6-closed-loop-Th-Th*}\ba{ll}
\ds \Th^*=\begin{pmatrix}-0.9949&-0.2393\\
\ns -0.1979&-0.8072\end{pmatrix}, \qq 
\bar\Th^*=\begin{pmatrix}-1.1798&0.0117\\-0.0082&-0.7971\end{pmatrix}.
\ea\ee
Substituting the above $\Th^*$ and $\bar\Th^*$ of (\ref{e5.6-closed-loop-Th-Th*}) into (\ref{e5.6-static-stabilizing-judgement-condition}), we can get
$$P_0=\begin{pmatrix}0.3323&-0.0865\\-0.0865&0.3365\end{pmatrix}>0, \qq
\bar P_0=\begin{pmatrix}0.1196&-0.0327\\-0.0327&0.2115\end{pmatrix}>0.$$
such that the negative definiteness conditions of (\ref{e5.6-static-stabilizing-judgement-condition}) hold. Then $(\Th^*,\bar\Th^*)$ above is an MF-$L^2$-stabilizer of the system $[A,\bar A,C,\bar C;B,\bar B,D,\bar D]$.
Hence, by Theorem \ref{T-3.6}, the problem admits a closed-loop Nash equilibrium.

\ms

Now, since
$$B(\Th^*-\Th^{**})=\begin{pmatrix}0.0051&0.0107\\0.0021&-0.0072\end{pmatrix}\ne0,$$
by Remark 2.7, we see that the closed-loop representation of the open-loop Nash equilibrium and the closed-loop Nash equilibrium are intrinsically different.

\ms

From this example, we see that the solutions to system  (\ref{AE-LQ-game-component})--(\ref{Th-Th*-component}), which come from the open-loop Nash equilibrium, may be symmetric, but they are different from the solution to the system of AREs (\ref{AE-LQ-game-closed})--(\ref{Sigma}), which come from the closed-loop Nash equilibrium. It is obvious that the closed-loop representation of the open-loop Nash equilibrium is different from the closed-loop Nash equilibrium.

\ms

The following example shows that for the closed-loop representation of open-loop Nash equilibria of the mean-field LQ non-zero sum stochastic differential game, it may happen that the solutions to the system of AREs may be asymmetric when the system $[A,\bar A,C,\bar C;B,\bar B,D,\bar D]$ is MF-$L^2$-stabilizable. And, the closed-loop representation of the open-loop Nash equilibrium is different from the closed-loop Nash equilibrium.

\ex
\bex{Example 5.7} \rm

Consider the two-dimensional state equation of form \rf{e5.6-1} with cost functionals \rf{e5.6-2}, but we let
$$\left\{\2n\ba{ll}
\ds A=\begin{pmatrix}-1&-1\\0&-1\\\end{pmatrix}, \q \bar A=\begin{pmatrix}-1&0\\-1&-1\\\end{pmatrix},\q B=\begin{pmatrix}0&0\\0&0\\\end{pmatrix}, \q \bar B=\begin{pmatrix}0&0\\0&0\\\end{pmatrix}, \\
\ns\ds C=\begin{pmatrix}0&0\\0&0\\\end{pmatrix},\q \bar C=\begin{pmatrix}0&0\\0&0\\\end{pmatrix},\q D=\begin{pmatrix}1&0\\0&1\\\end{pmatrix}, \q \bar D=\begin{pmatrix}1&0\\0&1\\\end{pmatrix}, \\
\ns\ds Q_1=\begin{pmatrix}2&1\\1&\frac{3}{2}\\\end{pmatrix},\q Q_2=\begin{pmatrix}1&\frac{1}{2}\\\frac{1}{2}&3\\\end{pmatrix},\q \bar Q_1=\begin{pmatrix}\frac{83}{44}&1\\ \ns 1&\frac{105}{44}\\\end{pmatrix}, \q
\bar Q_2=\begin{pmatrix}3&\frac{3}{2}\\ \ns \frac{3}{2}&\frac{16}{11}\\\end{pmatrix}, \\
\ns\ds  S_1=\begin{pmatrix}0&0\\\sqrt{\frac{5}{2}}&0\end{pmatrix}, \q S_2=\begin{pmatrix}0&0\\0&\sqrt{\frac{5}{2}}\end{pmatrix}, \q \bar S_1 = \bar S_2=\begin{pmatrix}0&0\\0&0\end{pmatrix}, \\
\ns\ds R_1=\begin{pmatrix}1&0\\0&1\end{pmatrix}, \q R_2=\begin{pmatrix}1&0\\0&\frac{3}{2}\end{pmatrix},\q \bar R_1=\begin{pmatrix}1&0\\0&1\end{pmatrix}, \q \bar R_2=\begin{pmatrix}\frac{1}{2}&0\\0&0\end{pmatrix}.
\ea\right.$$
From \rf{hat}, we have
$$\left\{\2n\ba{ll}
\ds\h A=\begin{pmatrix}-2&-1\\-1&-2\\\end{pmatrix}, \q \h B=\begin{pmatrix}0&0\\0&0\\\end{pmatrix}, \q \h C=\begin{pmatrix}0&0\\0&0\\\end{pmatrix},\q \h D=\begin{pmatrix}2&0\\0&2\\\end{pmatrix}, \\
\ns\ds \h Q_1=\begin{pmatrix}\frac{171}{44}&2\\2&\frac{171}{44}\\\end{pmatrix}, \q \h Q_2=\begin{pmatrix}4&2\\2&\frac{49}{11}\\\end{pmatrix}, \q
\h S_1=\begin{pmatrix}0&0\\\sqrt{\frac{5}{2}}&0\end{pmatrix}, \q \h S_2=\begin{pmatrix}0&0\\0&\sqrt{\frac{5}{2}}\end{pmatrix},\\
\ns\ds \h R_1=\begin{pmatrix}2&0\\0&2\end{pmatrix}, \q \h R_2=\begin{pmatrix}\frac{3}{2}&0\\0&\frac{3}{2}\end{pmatrix}.
\ea\right.$$

To look at the closed-representation of open-loop Nash equilibria, we solve the corresponding (\ref{AE-LQ-game-component})--(\ref{Th-Th*-component}) to get
$$
P_1=\begin{pmatrix}1&-\frac{1}{2}\\0&1\\\end{pmatrix}, \qq P_2=\begin{pmatrix}\frac{1}{2}&0\\0&1\\\end{pmatrix}, \qq
\h P_1=\begin{pmatrix}1&-\frac{5}{44}\\0&1\\\end{pmatrix}, \qq \h P_2=\begin{pmatrix}1&0\\0&1\\\end{pmatrix}.
$$
We see that $P_1$ and $\h P_1$ are not symmetric. Since
$$\begin{pmatrix}R_{111}+D_1^\top P_1D_1&R_{112}+D_1^\top P_1D_2\\R_{221}+D_2^\top P_2D_1&R_{222}+D_2^\top P_2D_2\end{pmatrix}
=\begin{pmatrix}2&-\frac{1}{2}\\0&\frac{5}{2}\end{pmatrix}$$
and
$$\begin{pmatrix}\h R_{111}+\h D_1^\top P_1\h D_1&\h R_{112}+\h D_1^\top P_1\h D_2\\\h R_{221}+\h D_2^\top P_2\h D_1&\h R_{222}+\h D_2^\top P_2\h D_2\end{pmatrix}
=\begin{pmatrix}6&-2\\0&\frac{11}{2}\end{pmatrix}$$
are invertible, it follows that 
\bel{e5.3-open-loop-Th-Th*}\ba{ll}
\ds \Th^{**}=\begin{pmatrix}0&-0.1581\\ \ns 0&-0.6325\end{pmatrix}, \q
\bar\Th^{**}=\begin{pmatrix}0&-0.0958\\0&-0.2875\end{pmatrix}.
\ea\ee
Then, similar to Example 5.6, making use of Theorem \ref{T-3.3}, we see that the problem admits an open-loop Nash equilibrium for any initial state $x\in\dbR^n$. And it has the closed-loop representation of form \rf{closed-loop control}.

\ms

Next, we consider the closed-loop Nash equilibria.
According to Proposition A.5 in \cite{Huang-Li-Yong2015}, the system $[A,\bar A,C,\bar C;$ $B,\bar B,D,\bar D]$ above is MF-$L^2$-stabilizable, and
$(\Th,\bar\Th)\in\sS[A,\bar A,C,\bar C;B,\bar B,D,\bar D]$ if and only if there exist $P_0>0$ and $\bar P_0>0$, such that
\bel{e5.7-static-stabilizing-judgement-condition}\left\{\2n\ba{ll}
\ds (A+B\Th)P_0+P_0(A+B\Th)^\top+(C+D\Th)P_0(C+D\Th)^\top+(\h C+\h D\bar\Th)\bar P_0(\h C+\h D\bar\Th)^\top\\
\ns\ds =\left[\begin{pmatrix}-1&-1\\0&-1\\\end{pmatrix}+\Th\right]P_0+P_0\left[\begin{pmatrix}-1&-1\\0&-1\\\end{pmatrix}+\Th\right]^\top
+\Th P_0\Th^\top+\begin{pmatrix}2&0\\0&2\\\end{pmatrix}\bar\Th\bar P_0\bar\Th^\top\begin{pmatrix}2&0\\0&2\\\end{pmatrix}<0, \\
\ns\ds (\h A+\h B\bar\Th)\bar P_0+\bar P_0(\h A+\h B\bar\Th)^\top=\begin{pmatrix}-2&-1\\-1&-2\\\end{pmatrix}\bar P_0
+\bar P_0\begin{pmatrix}-2&-1\\-1&-2\\\end{pmatrix}^\top<0.
\ea\right.\ee
For $i=1,2$, we solve corresponding \rf{AE-LQ-game-closed}--\rf{Sigma} to get
$$\left\{\2n\ba{ll}
\ds P_1=\begin{pmatrix}1&-0.4955\\-0.4955&1.7645\\\end{pmatrix}, \qq P_2=\begin{pmatrix}0.5&0\\0&1.0226\\\end{pmatrix}, \\
\ns\ds\h P_1=\begin{pmatrix}1.0647&-0.1861\\-0.1861&1.2327\\\end{pmatrix}, \qq \h P_2=\begin{pmatrix}1.0016&-0.0032\\-0.0032&1.0111\\\end{pmatrix},
\ea\right.$$
$$\left\{\2n\ba{ll}
\Si_1\deq R_{111}+D_1^\top P_1D_1=2, \qq \Si_2\deq R_{222}+D_2^\top P_2D_2=2.5226, \\
\bar\Si_1\deq \h R_{111}+\h D_1^\top P_1\h D_1=6, \qq \bar\Si_2\deq \h R_{222}+\h D_2^\top P_2\h D_2=5.5902.
\ea\right.$$
and
\bel{e5.7-closed-loop-Th-Th*}\ba{ll}
\ds \Th^*=\begin{pmatrix}0&-0.1553\\
\ns 0&-0.6268\end{pmatrix}, \qq 
\bar\Th^*=\begin{pmatrix}0&-0.0934\\0&-0.2828\end{pmatrix}.
\ea\ee
Substituting the above $\Th^*$ and $\bar\Th^*$ of (\ref{e5.7-closed-loop-Th-Th*}) into (\ref{e5.7-static-stabilizing-judgement-condition}), we can get
$$P_0=\begin{pmatrix}0.3082&-0.0910\\-0.0910&0.3148\end{pmatrix}>0, \qq
\bar P_0=\begin{pmatrix}0.1901 &-0.0721\\-0.0721&0.1870\end{pmatrix}>0.$$
such that the negative definiteness conditions of (\ref{e5.7-static-stabilizing-judgement-condition}) hold.
Then $(\Th^*,\bar\Th^*)$ above is an MF-$L^2$-stabilizer of the system $[A,\bar A,C,\bar C;B,\bar B,D,\bar D]$.
Hence, by Theorem \ref{T-3.6}, the problem admits a closed-loop Nash equilibrium.

\ms

Finally, similar to Example 5.6, we easily check that the closed-loop representation of the optimal open-loop strategy and the optimal closed-loop strategy are intrinsically different.

\ms

From this example, we see that the solutions to system \rf{AE-LQ-game-component}--\rf{Th-Th*-component}, which come from the open-loop saddle point, may be asymmetric, but the solutions to system (\ref{AE-LQ-game-closed})--(\ref{Sigma}), which come from the closed-loop saddle point, are still symmetric. It is obvious that the closed-loop representation of the open-loop Nash equilibrium is different from the closed-loop Nash equilibrium.
\ex

\section{Proof of Theorem 2.9}

This section is denoted to a proof of Theorem 2.9.

\ms

(ii) $\Ra$ (i) is obvious from the statements after Definition \ref{D-2.5}, and (iii) $\Ra$ (ii) is direct from Corollary \ref{C-4.6}.

\ms
We only need to prove (i) $\Ra$ (iii). We first consider the case ${\bf 0}\equiv(0,0)\in\sS[A,\bar A,C,\bar C;B,\bar B,D,\bar D]$. By Definition \ref{feedback}, system $[A,\bar A,C,\bar C]$ is $L^2$-globally integrable. From Proposition \ref{P-2.4}, we have
$\mathcal{U}_{ad}(x)=L^2_\dbF(\dbR^m),\forall x\in\dbR^n.$
This allows us to represent the cost functional $J(x;u(\cd))$ of \rf{cost-LQ-control} as a quadratic functional on the Hilbert space $L^2_\dbF(\dbR^m)$.
\bl{L-6.1} \sl
Suppose system $[A,\bar A,C,\bar C]$ is $L^2$-globally integrable. Then there exists a bounded self-adjoint linear operator $\mathcal{M}_2:L^2_\dbF(\dbR^m)\rightarrow L^2_\dbF(\dbR^m)$, a bounded linear operator $\mathcal{M}_1:\dbR^n\rightarrow L^2_\dbF(\dbR^m)$, an $\mathcal{M}_0\in\dbS^n$, and $\h u(\cd)\in L^2_\dbF(\dbR^m)$, $\h x\in\dbR^n$, $c\in\dbR$ such that
\bel{representation}\ba{ll}
J(x;u(\cd))=\lan\mathcal{M}_2u,u\ran+2\lan\mathcal{M}_1x,u\ran+\lan\mathcal{M}_0x,x\ran+2\lan\h u,u\ran+2\lan\h x,x\ran+c,\\
J^0(x;u(\cd))=\lan\mathcal{M}_2u,u\ran+2\lan\mathcal{M}_1x,u\ran+\lan\mathcal{M}_0x,x\ran,\qq \forall (x,u(\cd))\in\dbR^n\times L^2_\dbF(\dbR^m).\ea\ee\el
\proof
It is similar as Proposition 5.1 of \cite{Sun-Yong2018}, and the proof can be obtained by the method in \cite{Yong2013}. We omit the details here.
\endpf

\ms
With the help of representation \rf{representation}, we have the following results concerning with the open-loop solvability of Problem (MF-SLQ) whose proof is classical.
\bl{L-6.2} \sl
Suppose system $[A,\bar A,C,\bar C]$ is $L^2$-globally integrable. We have the following results:

\ms

{\rm(i)}\ Problem (MF-SLQ) is open-loop solvable at $x\in\dbR^n$ if and only if $\cM_2\ges0$, or equivalently,
\bel{convex}J^0(0;u(\cd))\ges0,\qq \forall u(\cd)\in L^2_\dbF(\dbR^m),\ee
and $\cM_1x+\h u\in\sR(\cM_2)$, where $\h u(\cd)\in L^2_\dbF(\dbR^m)$ is defined in \rf{representation}. In this case, $u^*(\cd)$ is an open-loop optimal control for the initial state $x$ if and only if $\mathcal{M}_2u^*+\mathcal{M}_1x+\h u=0$.

\ms

{\rm(ii)}\ If there exists a constant $\d>0$ such that $\mathcal{M}_2\ges\d I$, or equivalently,
\bel{uniformly convex}J^0(0;u(\cd))\ges\delta\dbE\int_0^\i|u(t)|^2dt,\qq \forall u(\cd)\in L^2_\dbF(\dbR^m),\ee
then Problem (MF-SLQ) is uniquely open-loop solvable.

\ms

{\rm(iii)}\ If Problem (MF-SLQ) is (uniquely) open-loop solvable, then Problem (MF-SLQ)$^0$ is (uniquely) open-loop solvable.
\el
We need the following lemma about Problem (MF-SLQ)$^0$.
\bl{L-6.3} \sl
Suppose system $[A,\bar A,C,\bar C]$ is $L^2$-globally integrable and \rf{uniformly convex} holds for some $\d>0$. Then the following system of coupled AREs
\bel{ARE-LQ-control-convertible}\left\{\2n\ba{ll}
\ds PA+A^\top P+C^\top PC+Q-(PB+C^\top PD+S^\top)\Si^{-1}(B^\top P+D^\top PC+S)=0,\\
\ns\ds\h P\h A+\h A^{\,\top}\h P+\h C^\top P\h C+\h Q-(\h P\h B+\h C^\top P\h D+\h S^\top)\bar\Si^{-1}(\h B^\top\h P+\h D^\top P\h C+\h S)=0,\\
\ns\ds\Si\equiv R+D^\top PD>0,\q\bar\Si\equiv\h R+\h D^\top P\h D>0,\ea\right.\ee
admits a solution pair $(P,\h P)\in\dbS^n\times\dbS^n$, and $(\Th,\bar\Th)\in\dbR^{m\times n}\times\dbR^{m\times n}$ defined by
\bel{Th-barTh-convertible}
\Th=-\Si^{-1}(B^\top P+D^\top PC+S),\qq \bar\Th=-\bar\Si^{-1}(\h B^\top\h P+\h D^\top P\h C+\h S\,),\ee
is an MF-$L^2$-stabilizer of system $[A,\bar A,C,\bar C;B,\bar B,D,\bar D]$. Furthermore, the unique open-loop optimal control $u_x^*(\cd)$ of Problem (MF-SLQ)$^0$ at $x$ is given by
\bel{unique open-loop optimal control}
u_x^*(\cd)=\Th\big\{X_x^*(\cd)-\dbE[X_x^*(\cd)]\big\}+\bar\Th\dbE[X_x^*(\cd)],
\ee
where $X_x^*(\cd)\in\sX[0,\i)$ is the solution to the following closed-loop system:
\bel{closed-loop system-A}\left\{\2n\ba{ll}
\ds dX_x^*(t)=\big\{A_\Th X_x^*(t)+\bar A_\BTh\dbE[X_x^*(t)]\big\}dt+\big\{C_\Th X_x^*(t)+\bar C_\BTh\dbE[X_x^*(t)]\big\}dW(t),\q t\ges0,\\
\ns\ds X_x^*(0)=x.\ea\right.\ee
Moreover, the value function is given by
\bel{V0(x)-control} V^0(x)=J^0(x;u_x^*(\cd))=\blan\h P x,x\bran,\q \forall x\in\dbR^n.\ee
\el
\proof
For $T>0$, we consider the state equation on $[0,T]$:
\begin{equation*}\left\{\2n\ba{ll}
\ds dX_T(t)=\big[AX_T(t)+\bar A\dbE[X_T(t)]+Bu(t)+\bar B\dbE[u(t)]\big]dt\\
\ns\ds\qq\qq\qq +\big[CX_T(t)+\bar C\dbE[X_T(t)]+Du(t)+\bar D\dbE[u(t)]\big]dW(t),\q t\in[0,T],\\
\ds X_T(0)=x,
\ea\right.\end{equation*}
and the cost functional
$$\ba{ll}
&\ds J_T^0(x;u(\cd))\triangleq\dbE\int_0^T\Bigg[\Bigg\langle\2n
\begin{pmatrix}Q&S^\top\\ S&R\end{pmatrix}
\begin{pmatrix}X_T\\ u\end{pmatrix},\begin{pmatrix}X_T\\ u\end{pmatrix}\2n\Bigg\rangle+\Bigg\langle\2n
\begin{pmatrix}\bar Q&\bar S^\top\\ \bar S&\bar R\\\end{pmatrix}
\begin{pmatrix}\dbE[X_T]\\ \dbE[u]\end{pmatrix},\begin{pmatrix}\dbE[X_T]\\ \dbE[u]\end{pmatrix}
\2n\Bigg\rangle\Bigg]dt.\ea$$
We claim that
\bel{uniformly convex-T}
J_T^0(0;u(\cd))\ges\d\dbE\int_0^T|u(t)|^2dt,\q \forall u(\cd)\in L^2_\dbF(\dbR^m),\q \mbox{for some }\d>0.
\ee
To prove this, choose any $u(\cd)\in L^2_\dbF(\dbR^m)$ and let $X_T(\cd)$ be the corresponding solution to the above state equation with initial state $x$. Define the zero-extension of $u(\cd)$ as follows:
$$\bar v(t)\triangleq[u(\cd)\oplus0{\bf 1}_{(T,\i)}](t)=\left\{\begin{array}{ll}u(t),&t\in[0,T],\\0,&t\in(T,\i).\end{array}\right.$$
Thus $\bar v(\cd)\in L^2_\dbF(\dbR^m)$, and the solution $X(\cd)$ to
\begin{equation*}\left\{\2n\ba{ll}
\ds dX(t)=\big\{AX(t)+\bar A\dbE[X(t)]+B\bar v(t)+\bar B\dbE[\bar v(t)\}\big]dt\\
\ns\ds\qq\qq\qq +\big\{CX(t)+\bar C\dbE[X(t)]+D\bar v(t)+\bar D\dbE[\bar v(t)]\big\}dW(t),\q t\geq0,\\
\ds X(0)=x,\ea\right.\end{equation*}
satisfies
$X(t)=X_T(t)\chi_{[0,T]}(t)$, $t\in[0,\i)$. It is obvious that
$J_T^0(x;u(\cd))=J^0(x;\bar v(\cd)).$
In particular, taking $x=0$, by \rf{uniformly convex} we get
$$J^0_T(0;u(\cd))=J^0(0;\bar v(\cd))\ges\d\dbE\int_0^\i|\bar v(t)|^2dt=\d\dbE\int_0^T|u(t)|^2dt.$$
This proves our claim. The fact \rf{uniformly convex-T} allows us to apply Theorems 4.2, 4.4 and 5.2 of \cite{Sun2017} to conclude that for any $T>0$, the following system of two coupled differential Riccati equations
\begin{equation*}\left\{\2n\ba{ll}
\ds \dot{P}(t;T)+P(t;T)A+A^\top P(t;T)+C^\top P(t;T)C+Q\\
\ns\ds\qq-\big[P(t;T)B+C^\top P(t;T)D+S^\top\big]\Si(t;T)^{-1}\big[B^\top P(t;T)+D^\top P(t;T)C+S\big]=0,\q t\in[0,T],\\
\ns\ds \dot{\h P}(t;T)+\h P(t;T)\h A+\h A^\top\h P(t;T)+\h C^\top P(t;T)\h C+\h Q\\
\ns\ds\qq-\big[\h P(t;T)\h B+\h C^\top P(t;T)\h D+\h S^\top\big]\bar\Si(t;T)^{-1}\big[\h B^\top\h P(t;T)+\h D^\top P(t;T)\h C+\h S\big]=0,\q t\in[0,T],\\
\ns\ds\Si(t;T)\equiv R+D^\top P(t;T)D>0,\q \bar\Si(t;T)\equiv \h R+\h D^\top P(t;T)\h D>0,\q P(T;T)=0,\q \h P(T;T)=0,\\
\ea\right.\end{equation*}
admits a unique solution pair $(P(\cd\,;T),\h P(\cd\,;T))\in C([0,T];\dbS^n)\times C([0,T];\dbS^n)$, and
$$V_T^0(x)\triangleq\inf\limits_{u(\cd)\in L^2_\dbF([0,T];\dbR^m)}J_T^0(x;u(\cd))=\blan\h P(0;T)x,x\bran,\q \forall x\in\dbR^n.$$
Further, if we define
\begin{equation*}\left\{\2n\ba{ll}
\ds \Th(t;T)=-\Si(t;T)^{-1}\big[B^\top P(t;T)+D^\top P(t;T)C+S\big],\\
\ns\ds \bar\Th(t;T)=-\bar\Si(t;T)^{-1}\big[\h B^\top\h P(t;T)+\h D^\top P(t;T)\h C+\h S\big],\q t\in[0,T],
\ea\right.\end{equation*}
then the unique open-loop optimal control $u_T^*(\cd)$ is given by
$$u_T^*(\cd)=\Th(t;T)(X_T^*(\cd)-\dbE[X_T^*(\cd)])+\bar\Th(t;T)\dbE[X_T^*(\cd)],$$
where $X_T^*(\cd)\in L^2_\dbF([0,T];\dbR^n)$ is the solution to the following MF-SDE:
\begin{equation*}\left\{\2n\ba{ll}
\ds dX_T^*(t)=\big\{\big(A+B\Th(t;T)\big)X_T^*(t)+\big(\bar A+\bar B\bar\Th(t;T)+B(\bar\Th(t;T)-\Th(t;T))\big)\dbE[X_T^*(t)]\big\}dt\\
\ns\ds\qq\qq +\big\{\big(C+D\Th(t;T)\big)X_T^*(t)+\big(\bar C+\bar D\bar\Th(t;T)+D(\bar\Th(t;T)-\Th(t;T))\big)\dbE[X_T^*(t)]\big\}dW(t),\ t\in[0,T],\\
\ds X_T^*(0)=x.\ea\right.\end{equation*}

Similar to Theorem 5.2 of \cite{Huang-Li-Yong2015}, we can show that $\lim\limits_{T\rightarrow\i}P(t;T)=P$ and $\lim\limits_{T\rightarrow\i}\h P(t;T)=\h P$,
for all $t\ges0$, where $(P,\h P)$ satisfies \rf{ARE-LQ-control-convertible}. Thus $\lim\limits_{T\rightarrow\i}\Th(t;T)=\Th$
and $\lim\limits_{T\rightarrow\i}\bar\Th(t;T)=\bar\Th$, for all $t\ges0$, which satisfies \rf{Th-barTh-convertible}.
And $\lim\limits_{T\rightarrow\i}X_T^*(t)=X_x^*(t)$, for all $t\ges0$, which satisfies \rf{closed-loop system-A}.
Further, applying It\^{o}'s formula to $\blan P(X_x^*(\cd)-\dbE[X_x^*(\cd)]),X_x^*(\cd)-\dbE[X_x^*(\cd)]\bran+\blan\h P\dbE[X_x^*(\cd)],\dbE[X_x^*(\cd)]\bran$, we have
\begin{equation*}\ba{ll}
\ds \blan\h Px,x\bran=J^0(x;u_x^*(\cd))\\
\ns\ds =\blan\h P x,x\bran+\dbE\int_0^\i\Big[\blan\big(PA+A^\top P+C^\top PC+Q\big)(X_x^*-\dbE[X_x^*]),X_x^*-\dbE[X_x^*]\bran\\
\ns\ds \qq\qq\qq\qq +2\blan\big(PB+C^\top PD+S^\top\big)(u_x^*-\dbE[u_x^*]),X_x^*-\dbE[X_x^*]\bran+\blan\Si(u_x^*-\dbE[u_x^*],u_x^*-\dbE[u_x^*]\bran\\
\ns\ds \qq\qq\qq\qq +\blan P\big(\h C\dbE[X_x^*]+\h D\dbE[u_x^*]\big),\h C\dbE[X_x^*]+\h D\dbE[u_x^*]\bran\Big]dt\\
\ds \q +\dbE\int_0^\i\Big[\blan\big(\h PA+A^\top\h P+\h Q\big)\dbE[X_x^*],\dbE[X_x^*]\bran
   +2\blan\big(\h P\h B+\h S^\top\big)\dbE[u_x^*],\dbE[X_x^*]\bran+\blan\h R\dbE[u_x^*],\dbE[u_x^*]\bran\Big]dt\\
\ds =\blan\h P x,x\bran+\dbE\int_0^\i\Big[\blan\Si\big(u_x^*-\dbE[u_x^*]-\Th(X_x^*-\dbE[X_x^*])\big),u_x^*-\dbE[u_x^*]-\Th(X_x^*-\dbE[X_x^*])\bran\\
\ns\ds \qq\qq\qq\qq +\blan\bar\Si\big(\dbE[u_x^*]-\bar\Th\dbE[X_x^*]\big),\dbE[u_x^*]-\bar\Th\dbE[X_x^*]\bran\Big]dt.
\ea\end{equation*}
Since $\Si>0$ and $\bar\Si>0$, we must have \rf{unique open-loop optimal control} and \rf{V0(x)-control}. Since
$X_x^*(t)=X_T^*(t)\chi_{[0,T]}(t)$, $t\in(0,\i)$, we have $X_x^*(\cd)\in\sX[0,\i)$ for all $x\in\dbR^n$. We conclude that $(\Th,\bar\Th)$ is an MF-$L^2$-stabilizer of system $[A,\bar A,C,\bar C;B,\bar B,D,\bar D]$. The proof is complete.
\endpf

\ms

Now, we are in the position to continue the proof of (i) $\Ra$ (iii). By Lemma \ref{L-6.2} \rm (iii), Problem (MF-SLQ)$^0$ is open-loop solvable. For any $\e>0$, let us consider the state equation
\begin{equation*}\left\{\2n\ba{ll}
\ds dX(t)\1n=\1n\big\{AX(t)\1n+\1n\bar A\dbE[X(t)]\1n+\1n Bu(t)\1n+\1n\bar B\dbE[u(t)]\big\}dt\\
\ns\ds\qq\qq+\big\{CX(t)\1n+\1n\bar C\dbE[X(t)]\1n+\1n Du(t)\1n+\1n\bar D\dbE[u(t)]\big\}dW(t),\q t\1n\ges\1n0,\\
\ns\ds X(0)=x,\ea\right.\end{equation*}
and the cost functional
$$\ba{ll}
\ds J_\e^0(x;u(\cd))\triangleq J^0(x;u(\cd))+\e\dbE\int_0^\i|u(t)|^2dt\\
\ns\ds =\dbE\int_0^\i\Bigg[\Bigg\langle\2n
\begin{pmatrix}Q&S^\top\\ S&R\end{pmatrix}
\begin{pmatrix}X\\ u\end{pmatrix},\begin{pmatrix}X\\ u\end{pmatrix}\2n\Bigg\rangle
+\Bigg\langle\2n
\begin{pmatrix}\bar Q&\bar S^\top\\ \bar S&\bar R\\\end{pmatrix}
\begin{pmatrix}\dbE[X]\\ \dbE[u]\end{pmatrix},\begin{pmatrix}\dbE[X]\\ \dbE[u]\end{pmatrix}
\2n\Bigg\rangle\Bigg]dt+\e\dbE\int_0^\i|u|^2dt\\
\ns\ds =\dbE\int_0^\i\Bigg[\Bigg\langle\2n
\begin{pmatrix}Q&S^\top\\ S&R+\e I\end{pmatrix}
\begin{pmatrix}X\\ u\end{pmatrix},\begin{pmatrix}X\\ u\end{pmatrix}\2n\Bigg\rangle
+\Bigg\langle\2n
\begin{pmatrix}\bar Q&\bar S^\top\\ \bar S&\bar R\\\end{pmatrix}
\begin{pmatrix}\dbE[X]\\ \dbE[u]\end{pmatrix},\begin{pmatrix}\dbE[X]\\ \dbE[u]\end{pmatrix}\2n\Bigg\rangle\Bigg]dt.\ea$$
Denote the above problem by Problem (MF-SLQ)$^0_\e$, and the corresponding value function by $V^0_\e(\cd)$. By \rf{representation}, we have
$$J^0_\e(0;u(\cd))=\blan(\mathcal{M}_2+\e I)u,u\bran\ges\e\dbE\int_0^\i|u(t)|^2dt,\q \forall u(\cd)\in L^2_\dbF(\dbR^m).$$
Then by Lemma \ref{L-6.3}, Problem (MF-SLQ)$^0_\e$ admits a unique open-loop optimal control $u^*_\e(\cd;x)\in L^2_\dbF(\dbR^m)$ at $x$, which is given by
$$u^*_\e(\cd;x)=\Th_\e(\Psi_\e(\cd)-\dbE[\Psi_\e(\cd)])x+\bar\Th_\e\dbE[\Psi_\e(\cd)]x,$$
where $(\Th_\e,\bar\Th_\e)\in\dbR^{m\times n}\times\dbR^{m\times n}$ defined below is an MF-$L^2$-stabilizer of system $[A,\bar A,C,\bar C;B,\bar B,D,\bar D]$:
\bel{Th_e-barTh_e-convertible}
\Th_\e\triangleq-\Si_\e^{-1}(B^\top P_\e+D^\top P_\e C+S),\qq
\bar\Th_\e\triangleq-\bar\Si_\e^{-1}\big[\h B^\top\h P_\e+\h D^\top P_\e\h C+\h S\big],
\ee
and $\Psi_\e(\cd)\in\sX[0,\i)$ is the solution to the following $\dbR^{n\times n}$-matrix-valued closed-loop system:
\begin{equation*}\left\{\2n\ba{ll}
\ds d\Psi_\e(t)=\big\{A_{\Th_\e}\Psi_\e(t)+\bar A_{\BTh_\e}\dbE[\Psi_\e(t)]\big\}dt+\big\{C_{\Th_\e}\Psi_\e(t)+\bar C_{\BTh_\e}\dbE[\Psi_\e(t)]\big\}dW(t),\q t\ges0,\\
\ns\ds \Psi_\e(0)=I,\ea\right.\end{equation*}
where $A_{\Th_\e},\bar A_{\BTh_\e},C_{\Th_\e},\bar C_{\BTh_\e}$ are defined as in \rf{A-Th} such that
$$V^0_\e(x)=\blan\h P_\e x,x\bran,\q \forall x\in\dbR^n,$$
where $\h P_\e\in\dbS^n$ is the solution to the following ARE:
\bel{ARE-e}\left\{\2n\ba{ll}
\ns\ds\h P_\e\h A+\h A^{\,\top}\h P_\e+\h C^\top P_\e\h C+\h Q-\big(\h P_\e\h B+\h C^\top P_\e\h D+\h S^{\,\top}\big)\bar\Si_\e^{-1}\big(\h B^\top\h P_\e
+\h D^\top P_\e\h C+\h S\big)=0,\\
\ns\ds\bar\Si_\e\equiv\h R+\e I+\h D^\top P_\e\h D>0,\ea\right.\ee
with $P_\e\in\dbS^n$ satisfying
\bel{ARE-e-SunYong2018}\left\{\2n\ba{ll}
\ds P_\e A+A^\top P_\e+C^\top P_\e C+Q-(P_\e B+C^\top P_\e D+S^\top)\Si_\e^{-1}(B^\top P_\e+D^\top P_\e C+S)=0,\\
\ns\ds \Si_\e\equiv R+\e I+D^\top P_\e D>0.
\ea\right.\ee

We observe that \rf{ARE-e-SunYong2018} coincides with (5.16) in \cite{Sun-Yong2018} for the problem without mean fields. By the proof of Theorem 4.5 of \cite{Sun-Yong2018}, we know that along a sequence $\{\e_k\}_{k=1}^\i\subseteq(0,\i)$ with $\lim_{k\rightarrow\i}\e_k=0$, both
$$P=\lim_{k\rightarrow\i}P_{\e_k}\q\mbox{ and \q }\Th=\lim_{k\rightarrow\i}\Th_{\e_k}$$
exist, which solves the following ARE:
\bel{ARE-SunYong2018}\left\{\2n\ba{ll}
\ds PA+A^\top P+C^\top PC+Q-(PB+C^\top PD+S^\top)\Si^\dag(B^\top P+D^\top PC+S)=0,\\
\ns\ds\sR\big(B^\top P+D^\top PC+S\big)\subseteq\sR(\Si),\q\lim_{k\to\i}\Si_{\e_k}=\Si\ges0.\ea\right.\ee
Hence $\ds\lim_{k\to\i}\bar\Si_{\e_k}=\bar\Si\ges0$ and
\bel{Th-control}
\Th=-\Si^\dag(B^\top P+D^\top PC+S)+(I-\Si^\dag\Si)\th,\q\hb{for some\ }\th\in\dbR^{m\times n}.\ee

Next, let us consider the limits of ARE \rf{ARE-e} and $\bar\Th_\e$ when $\e\rightarrow0$. For this target, consider the state equation
\bel{state-det}
 \dot y(t)=\h Ay(t)+\h Bw(t),\q t\ges0,\q y(0)=x,\ee
and the cost functional
\bel{cost-det}
\bar J(x;w(\cd))\deq\int_0^\i\Big[\blan\big(\h Q+\h C^\top P\h C\big)y(t),y(t)\bran+2\blan\big(\h D^\top P\h C+\h S\,\big)y(t),w(t)\bran+\blan\bar\Si w(t),w(t)\bran\Big]dt.\ee
We pose the following deterministic LQ optimal control problem:

\ms

{\bf Problem (DLQ).}\q For any initial state $x\in\dbR^n$, to find a control $w^*(\cd)\in L^2(\dbR^m)$ such that
$$\bar V(x)\triangleq\bar J(x;w^*(\cd))=\inf\limits_{w(\cd)\in L^2(\dbR^m)}\bar J(x;w(\cd)),$$
and have the following lemma.
\bl{L-6.4} \sl
Suppose system $[A,\bar A,C,\bar C]$ is $L^2$-globally integrable and \rf{uniformly convex} holds for some $\d>0$. Then the map $w(\cd)\mapsto\bar J(x;w(\cd))$ is uniformly convex. Consequently, the following ARE:
\bel{ARE-det-convertible}\left\{\2n\ba{ll}
\ds \bar{\h P}\h A+\h A^\top\bar{\h P}+\h C^\top P\h C+\h Q-\big(\bar{\h P}\h B+\h C^\top P\h D+\h S^\top\big)\bar\Si^{-1}\big(\h B^\top\bar{\h P}+\h D^\top P\h C+\h S\big)=0,\\
\ns\ds \bar\Si\equiv \h R+\h D^\top P\h D>0,\ea\right.\ee
admits a unique solution $\bar{\h P}\in\dbS^n$, with $P\in\dbS^n$ satisfying the first equation of \rf{ARE-LQ-control-convertible}, such that $\bar V(x)=\blan\bar{\h P}x,x\bran$ for any $x\in\dbR^n$. Further,
\bel{Gamma}\ba{ll}
\bar\G\triangleq-\bar\Si^{-1}\big(\h B^\top\bar{\h P}+\h D^\top P\h C+\h S\big)\in\dbR^{m\times n}\ea\ee
is a stabilizer of system $\big[\h A;\h B\big]$. And the unique open-loop optimal control $w_x^*(\cd)$ at $x$ is given by
$$w_x^*(\cd)=\bar\G\bar\Psi(\cd)x,$$
where $\bar\Psi(\cd)\in L^2(\dbR^{n\times n})$ is the solution to the following matrix-valued closed-loop system:
$$d\bar\Psi(t)=\big(\h A+\h B\bar\G\big)\bar\Psi(t)dt,\q t\ges0,\q \bar\Psi(0)=I.$$
\el
\proof
Let $P\in\dbS^n$ be the solution to the first equation of \rf{ARE-LQ-control-convertible} and set $\Th=-\Si^{-1}(B^\top P+D^\top PC+S)$. Thus \rf{ARE-SunYong2018} is equivalent to
\bel{AE-SunYong2018}
P(A+B\Th)+(A+B\Th)^\top P+(C+D\Th)^\top P(C+D\Th)+\Th^\top R\Th+S^\top\Th+\Th^\top S+Q=0.\ee
We claim that
\bel{J0-barJ}
J^0(0;\Th X^w(\cd)+w(\cd))=\bar J(0;\Th y(\cd)+w(\cd)),\q \forall w(\cd)\in L^2(\dbR^m),
\ee
where $X^w(\cd)\in\sX(0,\i)$ satisfies
\begin{equation*}\left\{\2n\ba{ll}
\ds dX^w(t)=\big\{AX^w(t)+\bar A\dbE[X^w(t)]+B\big(\Th X^w(t)+w(t)\big)+\bar B\dbE[\Th X^w(t)+w(t)]\big\}dt\\
\ns\ds\qq\qq\qq +\big\{CX^w(t)+\bar C\dbE[X(t)]+D\big(\Th X^w(t)+w(t)\big)+\bar D\dbE[\Th X^w(t)+w(t)]\big\}dW(t),\q t\geq0,\\
\ds X^w(0)=0,
\ea\right.\end{equation*}
and $y(\cd)\in L^2(\dbR^n)$ is the solution to
$$\dot y(t)=\h Ay(t)+\h B\big[\Th y(t)+w(t)\big],\q t\geq0,\q y(0)=0.$$
In fact, noting that $w(\cd)$ is deterministic, it implies that
\begin{equation*}\left\{\2n\ba{ll}
\ds d\dbE[X^w(t)]=\big\{\h A\dbE[X^w(t)]+\h B\big(\Th\dbE[X^w(t)]+w(t)\big)\big\}dt, \q t\geq0,\\
\ns\ds \dbE[X^w(0)]=0.\ea\right.\end{equation*}
Thus by the uniqueness of solutions, $\dbE[X^w(t)]=y(t),\ t\ges0$. Now let $z(\cd)=X^w(\cd)-\dbE[X^w(\cd)]$, we have
\begin{equation*}\left\{\2n\ba{ll}
\ds dz(t)=(A+B\Th)z(t)dt+\big[(C+D\Th)z(t)+\h Cy(t)+\h D\big(\Th y(t)+w(t)\big)\big]dW(t),\q t\geq0,\\
\ns\ds z(0)=0.
\ea\right.\end{equation*}
Applying It\^{o}'s formula to $\blan Pz(\cd),z(\cd)\bran$, noting \rf{AE-SunYong2018} and $\dbE[z(\cd)]=0$, we have
\begin{equation*}\ba{ll}
\ds J^0(0;\Th X^w(\cd)+w(\cd))\\
%
%
%
%
%
%
%
\ns\ds =\dbE\int_0^\i\Bigg[\Bigg\langle\2n
\begin{pmatrix}Q&S^\top\\ S&R\end{pmatrix}
\begin{pmatrix}z\\ \Th z\end{pmatrix},\begin{pmatrix}z\\ \Th z\end{pmatrix}\2n\Bigg\rangle+\Bigg\langle\2n
\begin{pmatrix}\h Q& \h S^\top\\ \h S& \h R\\\end{pmatrix}
\begin{pmatrix}y\\ \Th y+w\end{pmatrix},\begin{pmatrix}y\\ \Th y+w\end{pmatrix}\2n\Bigg\rangle\Bigg]dt\\
\ns\ds =\dbE\int_0^\i\Big[\blan P(A+B\Th)z,z\bran+\blan Pz,(A+B\Th)z\bran\\
\ns\ds\qq\qq +\blan P\big[(C+D\Th)z+\h Cy+\h D\big(\Th y+w\big)\big],(C+D\Th)z+\h Cy+\h D\big(\Th y+w\big)\bran\\
\ns\ds\qq\qq +\blan(\Th^\top R\Th+S^\top\Th+\Th^\top S+Q)z,z\bran+\blan\h Qy,y\bran+2\blan\h Sy,\Th y+w\bran+\blan\h R(\Th y+w),\Th y+w\bran\Big]dt\\
\ns\ds =\int_0^\i\Big[\blan\big(\h Q+\h C^\top P\h C\big)y,y\bran+2\blan\big(\h D^\top\h C+\h S\big)y,\Th y+w\bran+\blan\bar\Si(\Th y+w),\Th y+w\bran\Big]dt\\
\ns\ds =\bar J(0;\Th y(\cd)+w(\cd)).
\ea\end{equation*}
Thus \rf{J0-barJ} holds. Consequently, by \rf{uniformly convex} and Jensen's inequality, there exists some $\d>0$, such that
\begin{equation*}
\begin{aligned}
&\bar J(0;\Th y(\cd)+w(\cd))=J^0(0;\Th X^w(\cd)+w(\cd))\ges\d\dbE\int_0^\i|\Th X^w(t)+w(t)|^2dt\\
&\ges\d\int_0^\i\big|\dbE[\Th X^w(t)+w(t)]\big|^2dt=\d\int_0^\i\big|\Th y(t)+w(t)\big|^2dt,\q \forall w(\cd)\in L^2(\dbR^m).
\end{aligned}
\end{equation*}
This implies the uniform convexity of $w(\cd)\mapsto\bar J(x;w(\cd))$. Repeating the finite time interval approximation method in the proof of Lemma \ref{L-6.3}, the rest of the theorem follows and the proof is complete.
\endpf

\ms
The following additional lemma is a direct consequence of Proposition 5.2, (iii) of \cite{Sun-Yong2018}.
\bl{L-6.5} \sl
If Problem (DLQ) is open-loop solvable, then there exists a $U^*(\cd)\in L^2(\dbR^{m\times n})$ such that for any $x\in\dbR^n$, $U^*(\cd)x$ is an open-loop optimal control for the initial state $x$.
\el

With the above two lemmas in hand, we can continue to consider the limitation of \rf{ARE-e} when $\e\rightarrow0$, without \rf{uniformly convex}. For any $\e>0$, consider the state equation \rf{state-det} and the cost functional
$$
\ba{ll}
\ds \bar J_\e(x;w(\cd))\triangleq\bar J(x;w(\cd))+\e\int_0^\i|w(t)|^2dt\\
\ns\ds =\int_0^\i\Big[\blan\big(\h Q+\h C^\top P\h C\big)y(t),y(t)\bran+2\blan\big(\h D^\top P\h C+\h S\big)y(t),w(t)\bran+\blan\bar\Si_\e w(t),w(t)\bran\Big]dt.
\ea$$
Denote the above problem by Problem (DLQ)$_\e$, and the corresponding value function by $\bar V_\e(\cd)$. Since Problem (MF-SLQ)$^0$ is open-loop solvable, by \rf{J0-barJ} and Lemma \ref{L-6.1}, we must have
$$\bar J_\e(0;w(\cd))=\blan(\wt\cM+\e I)w,w\bran\ges\e\int_0^\i|w(t)|^2dt,\q \forall w(\cd)\in L^2(\dbR^m),$$
where $\wt\cM$ is a bounded self-adjoint linear operator from $L^2(\dbR^m)$ to itself. That is to say, Problem (DLQ)$_\e$ is uniquely open-loop solvable, for each $\e>0$. Consequently, by Lemma \ref{L-6.4}, ARE \rf{ARE-e} admits a unique solution $\h P_\e\in\dbS^n$ such that $\bar V_\e(x)=\blan\h P_\e x,x\bran,\forall x\in\dbR^n$. Moreover, $\bar\Th_\e\in\dbR^{m\times n}$ is a stabilizer of system $\big[\h A;\h B\big]$, and the unique open-loop optimal control $w^*_\e(\cd;x)$ of Problem (DLQ)$_\e$ at $x$ is given by
\bel{w*-e}
w^*_\e(\cd;x)=\bar\Th_\e\bar\Psi_\e(\cd)x,
\ee
where $\bar\Psi_\e(\cd)\in L^2(\dbR^{n\times n})$ is the solution to the following closed-loop system:
$$d\bar\Psi_\e(t)=\big(\h A+\h B\bar\Th_\e\big)\bar\Psi_\e(t)dt,\q t\ges0,\q \bar\Psi_\e(0)=I.$$
Now let $U^*(\cd)\in L^2(\dbR^{m\times n})$ be a deterministic function with the property in Lemma \ref{L-6.5}. By the definition of value function, we have for any $x\in\dbR^n$ and $\e>0$,
\bel{U*-1}\ba{ll}
\ds \bar V(x)+\e\int_0^\i\big|w^*_\e(t;x)\big|^2dt\les \bar J(x;w^*_\e(\cd;x))+\e\int_0^\i\big|w^*_\e(t;x)\big|^2dt\\
\ns\ds =\bar J_\e(x;w^*_\e(\cd;x))=\bar V_\e(x)\les \bar J_\e(x;U^*(\cd)x)=\bar V(x)+\e\int_0^\i\big|U^*(t)x\big|^2dt,
\ea\ee
which implies
\bel{U*-2}
\bar V(x)\les\bar V_\e(x)=\blan\h P_\e x,x\bran\les\bar V(x)+\e\int_0^\i\big|U^*(t)x\big|^2dt,\q \forall x\in\dbR^n,\q \forall \e>0,
\ee
and
\bel{U*-3}
0\les\int_0^\i\bar\Psi_\e(t)^\top\bar\Th_\e^\top\bar\Th_\e\bar\Psi_\e(t)dt\les\int_0^\i U^*(t)^\top U^*(t)dt,\q \forall \e>0.\\
\ee
It is clear that $\bar V_{\e_1}(x)\les\bar V_{\e_2}(x)$, for any $0<\e_1\les\e_2$, $\forall x\in\dbR^n$. Thus $\h P_{\e_1}\les\h P_{\e_2}$ for any $0<\e_1\les\e_2$. Then, noting \rf{U*-2}, $\h P\equiv\lim_{\e\rightarrow0}\h P_\e$ exists and $\bar V(x)=\blan\h Px,x\bran$.

Moreover, noting that $\big[\h A+\h B\bar\Th_\e\big]$ is exponentially stable (see Theorem 3.7 of \cite{Huang-Li-Yong2015}), it is classical that the following Lyapunov equation for $\D_\e$:
$$\D_\e\big(\h A+\h B\bar\Th_\e\big)+\big(\h A+\h B\bar\Th_\e\big)^\top\D_\e+\bar\Th_\e^\top\bar\Th_\e=0.$$
admits a unique solution
$$\D_\e=\int_0^\i\bar\Psi_\e(t)^\top\bar\Th_\e^\top\bar\Th_\e\bar\Psi_\e(t)dt>0,$$
since $\bar\Th_\e^\top\bar\Th_\e>0$. It follows from \rf{U*-3} that $\{\D_\e\}_{\e>0}$ is bounded. Thus, there exists some constant $K>0$, such that
$$0<|\bar\Th_\e|^2=-\big\{\D_\e\big(\h A+\h B\bar\Th_\e\big)+\big(\h A+\h B\bar\Th_\e\big)^\top\D_\e\big\}\les K(1+|\bar\Th_\e|),\q\forall\e>0,$$
which implies the boundedness of $\{\bar\Th_\e\}_{\e>0}$. Without loss of generality, $\bar\Th\equiv\lim_{k\rightarrow\i}\bar\Th_{\e_k}$ exists. Then
\bel{Sigma-Gamma}\ba{ll}
\ns\ds\bar\Si\bar\Th=\lim_{k\to\i}\bar\Si_{\e_k}\bar\Th_{\e_k}=-\lim_{k\to\i}
\big(\h B^\top\h P_{\e_k}+\h D^\top P_{\e_k}\h C+\h S\big)=-\big(\h B^\top\h P+\h D^\top P\h C+\h S\,\big),\ea\ee
and thus
\bel{range-LQ-hat P}\left\{\2n\ba{ll}
\ds\sR\big(\h B^\top\h P+\h D^\top P\h C+\h S\,\big)\subseteq\sR\big(\bar\Si\big),\\
\ns\ds\bar\Th=-\bar\Si^\dag\big(\h B^\top\h P+\h D^\top P\h C+\h S\,\big)+\big(I-\bar\Si^\dag\bar\Si\big)\bar\th,\q\hb{for some\ }\bar\th\in\dbR^{m\times n}.\ea\right.\ee
Noting that by \rf{Th_e-barTh_e-convertible}, we have
$\h P_\e\h B+\h C^\top P\h D+\h S^\top=-\bar\Th_\e^\top\bar\Si_\e.$
Thus \rf{ARE-e} can be written as
\begin{equation*}\ba{ll}
\h P_\e\h A+\h A^\top\h P_\e+\h C^\top P\h C+\h Q-\bar\Th_\e^\top\bar\Si_\e\bar\Th_\e=0,\q \bar\Si_\e\equiv\h R+\e I+\h D^\top P_\e\h D>0.\ea\end{equation*}
Now, passing to the limit along $\{\e_k\}_{k=1}^\i$ in the above yields
\bel{ARE}\ba{ll}
\ds \h P\h A+\h A^\top\h P+\h C^\top P\h C+\h Q-\bar\Th^\top\bar\Si\bar\Th=0,\q \bar\Si\equiv\h R+\h D^\top P\h D\ges0,\ea\ee
which, together with \rf{range-LQ-hat P} and \rf{ARE-SunYong2018}, implies that $(P,\h P)$ solves \rf{ARE-LQ-control}.

Next, we show the pair $(P,\h P)$ is a static stabilizing solution to \rf{ARE-LQ-control}. Since $\bar\Th_{\e_k}\rightarrow\bar\Th$ as $k\rightarrow\i$, we have $\bar\Psi_{\e_k}(t)\rightarrow\bar\Psi(t)$ for all $t\ges0$, which satisfies
$$d\bar\Psi(t)=\big(\h A+\h B\bar\Th\big)\bar\Psi(t)dt,\q t\ges0,\q \bar\Psi(0)=I.$$
By Fatou's lemma, we have
$$
\int_0^\i\big|\bar\Th\bar\Psi(t)x\big|^2dt\les\liminf\limits_{k\rightarrow\i}\int_0^\i\big|\bar\Th_{\e_k}\bar\Psi_{\e_k}(t)x\big|^2dt
\les\int_0^\i\big|U^*(t)x\big|^2dt<\i,\q \forall x\in\dbR^n,$$
which implies $\bar\Th\bar\Psi(\cd)\in L^2(\dbR^{m\times n})$. Thus $\bar\Th\in\dbR^{m\times n}$ is a stabilizer of system $\big[\h A;\h B\big]$ and $\bar\Psi(\cd)\in L^2(\dbR^n)$. Next, $\Th_{\e_k}\rightarrow\Th$ as $k\rightarrow\i$ leads to $\Psi_{\e_k}(t)\to\Psi(t)$ for all $t\ges0,\as$, which solves
\begin{equation*}\left\{\2n\ba{ll}
\ds d\Psi(t)=\big\{A_\Th\Psi(t)+\bar A_\BTh\dbE[\Psi(t)]\big\}dt+\big\{C_\Th\Psi(t)+\bar D_\BTh\dbE[\Psi(t)]\big\}dW(t),\q t\ges0,\\
\ns\ds \Psi(0)=I.\ea\right.\end{equation*}
Since
$$d\dbE[\Psi(t)]=\big[\h A+\h B\bar\Th\big]\dbE[\Psi(t)]dt,\q t\ges0,\q \dbE[\Psi(0)]=I,$$
by the uniqueness of the solutions, we get $\bar\Psi(\cd)\equiv\dbE[\Psi(\cd)]\in L^2(\dbR^n)$. By Lemma 2.3 of \cite{Sun-Yong2018}, $\Psi(\cd)\in L^2_\dbF(\dbR^{n\times n})$. Then system $[A_\Th,C_\Th]$ is $L^2$-globally integrable, and $\Th$ is an $L^2$-stabilizer of system $[A,C;B,D]$. Moreover, $[A_\Th,\bar A_\BTh,C_\Th,\bar C_\BTh]$ is $L^2$-asymptotically stable by Proposition 2.3 of \cite{Huang-Li-Yong2015} since it is $L^2$-globally integrable. Therefore, $(\Th,\bar\Th)$ is an MF-$L^2$-stabilizer of system $[A,\bar A,C,\bar C;B,\bar B,D,\bar D]$ by Definition \ref{feedback}.

Now, we consider the BSDE \rf{BSDE-eta-control} and ODE \rf{BODE-bar-eta-control} on $[0,\i)$. By \rf{BSDE-eta-control} we have
\bel{BSDE-eta-control-Th}
-d\eta(t)=\big[A_\Th^\top\eta(t)+C_\Th^\top\z(t)+C_\Th^\top P\si(t)+\Th^\top\rho(t)+Pb(t)+q(t)\big]dt-\z(t)dW(t),\q t\ges0.
\ee
Since system $[A_\Th,C_\Th]$ is $L^2$-globally integrable, by Lemma 2.5 of \cite{Sun-Yong2018}, $(\eta(\cd),\z(\cd))\in\sX[0,\i)\times L^2_\dbF(\dbR^n)$. Noting \rf{BODE-bar-eta-control} is equivalent to
\bel{BODE-bar-eta-control-Th}\ba{ll}
\ns\ds-d\bar\eta(t)=\big\{\big(\h A+\h B\bar\Th\big)^\top\bar\eta(t)+\big(\h C+\h D\bar\Th\big)^\top\dbE[\z(t)]+\big(\h C+\h D\bar\Th\big)^\top P\dbE[\si(t)]\\
\ns\ds\qq\qq\qq+\bar\Th^\top\dbE[\rho(t)]+\dbE[q(t)]+\h P\dbE[b(t)]\big\}dt,\q t\ges0.\ea\ee
Since
$\int_0^\i|\dbE[\z(t)]|^2dt\les\dbE\int_0^\i|\z(t)|^2dt<\i,$
thus $\bar\eta(\cd)\in L^2(\dbR^n)$ since system $\big[\h A+\h B\bar\Th\big]$ is exponentially stable. Let $(x,u(\cd))\in\dbR^n\times L^2_\dbF(\dbR^m)$, and $X^u(\cd)\equiv X(\cd;x,u(\cd))$ satisfy \rf{state-LQ-control}. Applying It\^{o}'s formula to
$$\blan P\big(X^u(\cd)-\dbE[X^u(\cd)]\big)+2\eta(\cd),X^u(\cd)-\dbE[X^u(\cd)]\bran+\blan\h P\dbE[X^u(\cd)]+2\bar\eta(\cd),\dbE[X^u(\cd)]\bran,$$
noting \rf{ARE-LQ-control}, we have
\bel{Ito-4}\ba{ll}
\ds J(x;u(\cd))-\blan\h Px+2\bar\eta(0),x\bran\\
\ns\ds =\dbE\int_0^\i\Big[\blan\Si\big(u-\dbE[u]-\Th(X^u-\dbE[X^u])\big),u-\dbE[u]-\Th(X^u-\dbE[X^u])\bran\\
\ns\ds\qq +2\blan B^\top(\eta-\dbE[\eta])+D^\top(\z-\dbE[\z])+D^\top P(\si-\dbE[\si])+\rho-\dbE[\rho],u-\dbE[u]-\Th(X^u-\dbE[X^u])\bran\\
\ns\ds\qq +\blan\bar\Si\big(\dbE[u]\1n-\1n\bar\Th\dbE[X^u]\big),\dbE[u]\1n-\1n\bar\Th\dbE[X^u]
\bran\1n+\1n2\blan\h B^\top\bar\eta\1n+\1n\h D^\top P\dbE[\si]\1n+\1n\h D^\top\dbE[\z]\1n+\1n\dbE[\rho],\dbE[u]\1n-\1n\bar\Th\dbE[X^u]\bran\\
\ns\ds\qq +\blan P\si,\si\bran+2\blan\eta,b-\dbE[b]\bran+2\blan\z,\si\bran+2\blan\bar\eta,\dbE[b]\bran\Big]dt.
\ea\ee
Let $u^*(\cd)$ be an open-loop optimal control of Problem (MF-SLQ) for the initial state $x$. By Proposition \ref{P-2.4}, it admits the following closed-loop representation:
$$u^*(\cd)=\Th(X^*(\cd)-\dbE[X^*(\cd)])+\bar\Th\dbE[X^*(\cd)]+v^*(\cd),$$
for some $v^*(\cd)\in L^2_\dbF(\dbR^m)$, where $X^*(\cd)\equiv X^{\Th,\bar\Th,v^*}(\cd)$ is the solution to \rf{closed-loop-system}. Hence
\bel{optimal control u*-1}
J(x;u^*(\cd))\les J(x;u(\cd))\equiv J\big(x;\Th(X^{\Th,\bar\Th,v}(\cd)-\dbE[X^{\Th,\bar\Th,v}(\cd)])+\bar\Th\dbE[X^{\Th,\bar\Th,v}(\cd)]+v(\cd)\big),\ \forall v(\cd)\in L^2_\dbF(\dbR^m),\ee
where $X^{\Th,\bar\Th,v}(\cd)$ is the solution to \rf{closed-loop-system} corresponding to $(\Th,\bar\Th,v(\cd))$. From \rf{Ito-4} and \rf{optimal control u*-1}, we have that for any $v(\cd)\in L^2_\dbF(\dbR^m)$,
\bel{optimal control u*-2}\ba{ll}
\ds \blan\h Px,x\bran+2\blan\bar\eta(0),x\bran+\dbE\int_0^\i\Big[\blan P\si,\si\bran+2\blan\eta,b-\dbE[b]\bran
+2\blan\z,\si\bran+2\blan\bar\eta,\dbE[b]\bran\Big]dt\\
\ns\ds\q +\dbE\int_0^\i\Big[\blan\Si(v^*-\dbE[v^*]),v^*-\dbE[v^*]\bran+2\blan B^\top(\eta-\dbE[\eta])+D^\top(\z-\dbE[\z])+D^\top P(\si-\dbE[\si])+\rho\\
\ns\ds\qq\qq -\dbE[\rho],v^*-\dbE[v^*]\bran+\blan\bar\Si\dbE[v^*],\dbE[v^*]\bran
+2\blan\h B^\top\bar\eta+\h D^\top P\dbE[\si]+\h D^\top\dbE[\z]+\dbE[\rho],\dbE[v^*]\bran\Big]dt\\
\ns\ds =J(x;u^*(\cd))\les J(x;u(\cd))\\
\ns\ds =\blan\h Px,x\bran+2\blan\bar\eta(0),x\bran+\dbE\int_0^\i\Big[\blan P\si,\si\bran+2\blan\eta,b-\dbE[b]\bran+2\blan\z,\si\bran+2\blan\bar\eta,\dbE[b]\bran\Big]dt\\
\ns\ds\q +\dbE\int_0^\i\Big[\blan\Si(v-\dbE[v]),v-\dbE[v]\bran+2\blan B^\top(\eta-\dbE[\eta])+D^\top(\z-\dbE[\z])+D^\top P(\si-\dbE[\si])+\rho\\
\ns\ds\qq\qq\q -\dbE[\rho],v-\dbE[v]\bran+\blan\bar\Si\dbE[v],\dbE[v]\bran+2\blan\h B^\top\bar\eta+\h D^\top P\dbE[\si]+\h D^\top\dbE[\z]+\dbE[\rho],\dbE[v]\bran\Big]dt.\ea\ee
The above shows that $(v^*(\cd),\dbE[v^*(\cd)])$ is a minimizing pair of the functional
$$\ba{ll}
\ns\ds F(v(\cd),\dbE[v(\cd)])\triangleq\dbE\int_0^\i\Big[\blan\Si(v-\dbE[v]),v-\dbE[v]\bran+2\blan B^\top(\eta-\dbE[\eta])+D^\top(\z-\dbE[\z])+D^\top P(\si-\dbE[\si])\\
\ns\ds\qq\qq\qq\qq +\rho-\dbE[\rho],v-\dbE[v]\bran+\blan\bar\Si\dbE[v],\dbE[v]\bran+2\blan\h B^\top\bar\eta+\h D^\top P\dbE[\si]+\h D^\top\dbE[\z]+\dbE[\rho],\dbE[v]\bran\Big]dt.
\ea$$
Therefore,
$$\left\{\2n\ba{ll}
\ds \Si(v(\cd)\1n-\1n\dbE[v(\cd)])\1n+\1n B^\top\1n(\eta(\cd)\1n-\1n\dbE[\eta(\cd)])\1n+\1n D^\top\1n(\z(\cd)\1n-\1n\dbE[\z(\cd)])\1n
+\1n D^\top P(\si(\cd)\1n-\1n\dbE[\si(\cd)])\1n+\1n\rho(\cd)\1n-\1n\dbE[\rho(\cd)]=0,\q\as,\\
\ns\ds \bar\Si\dbE[v(\cd)]+\h B^\top\bar\eta(\cd)+\h D^\top P\dbE[\si(\cd)]+\h D^\top\dbE[\z(\cd)]+\dbE[\rho(\cd)]=0,\ea\right.$$
which implies
\bel{range-Sigma-bar Sigma}\left\{\2n\ba{ll}
\ds B^\top(\eta(\cd)-\dbE[\eta(\cd)])+D^\top(\z(\cd)-\dbE[\z(\cd)])+D^\top P(\si(\cd)-\dbE[\si(\cd)])+\rho(\cd)-\dbE[\rho(\cd)]\in\sR(\Si),\q \as,\\
\ns\ds\h B^\top\bar\eta(\cd)+\h D^\top P\dbE[\si(\cd)]+\h D^\top\dbE[\z(\cd)]+\dbE[\rho(\cd)]\in\sR\big(\bar\Si\big),
\ea\right.\ee
and
\bel{v-Ev and Ev}
v^*(\cd)-\dbE[v^*(\cd)]=\f(\cd)+\big(I-\Si^\dag\Si\big)(\nu(\cd)-\dbE[\nu(\cd)]),\q \dbE[v^*(\cd)]=\bar\f(\cd)+\big(I-\bar\Si^\dag\bar\Si\big)\bar\nu(\cd),
\ee
where $\f(\cd),\bar\f(\cd)$ are defined in \rf{phi-phi}, $\nu(\cd)\in L^2_\dbF(\dbR^m)$ and $\bar\nu(\cd)\in L^2(\dbR^m)$. By \rf{v-Ev and Ev}, we obtain
\bel{v}
v^*(\cd)=\f(\cd)-\dbE[\f(\cd)]+\bar\f(\cd)+\big(I-\Si^\dag\Si\big)(\nu(\cd)-\dbE[\nu(\cd)])+\big(I-\bar\Si^\dag\bar\Si\big)\bar\nu(\cd).
\ee
Further, combining \rf{optimal control u*-2}-\rf{v}, \rf{V(x)-control} is clear.

\ms
In the end, we give the proof of (i) $\Ra$ (iii) for the general case $\sS[A,\bar A,C,\bar C;B,\bar B,D,\bar D]\ne\varnothing$, that is, under (H1). Take $\BTh\equiv(\Th,\bar\Th)\in\sS[A,\bar A,C,\bar C;B,\bar B,D,\bar D]$, and consider the state equation \rf{closed-loop-system} with the cost functional $J^\BTh(x;v(\cd))$ by \rf{cost-LQ-control-closed-loop}. By Proposition 3.6 and Theorem 3.7 of \cite{Huang-Li-Yong2015}, system $[A_\Th,\bar A_\BTh,C_\Th,\bar C_\BTh]$ is $L^2$-globally integrable. We denote by Problem (MF-SLQ)$_\BTh$ the corresponding mean-field LQ stochastic optimal control problem. The following lemma lists some facts about it, whose proofs are straightforward consequences of Proposition \ref{P-2.4}.
\bl{L-6.6} \sl We have the following statements.

\ms

{\rm(i)}\ Problem (MF-SLQ)$_\BTh$ is open-loop solvable at $x\in\dbR^n$ if and only if so
is Problem (MF-SLQ). In this case, $v^*(\cd)$ is an open-loop optimal control of Problem (MF-SLQ)$_\BTh$ if and only if
$$u^*(\cd)\deq v^*(\cd)+\Th\big\{X^{v^*}(\cd)-\dbE[X^{v^*}(\cd)]\big\}+\bar\Th\dbE[X^{v^*}(\cd)]$$
is an open-loop optimal control of Problem (MF-SLQ).

\ms

{\rm(ii)}\ Problem (MF-SLQ)$_\BTh$ is closed-loop solvable if and only if so is Problem (MF-SLQ). In this case, $(\Th^*,\bar\Th^*,v^*(\cd))$ is a closed-loop optimal strategy of Problem (MF-SLQ)$_\BTh$ if and only if $(\Th^*+\Th,\bar\Th^*+\bar\Th,v^*(\cd))$
is a closed-loop optimal strategy of Problem (MF-SLQ).
\el

By Lemma \ref{L-6.6}, and the result for the $L^2$-globally integrable case, the following system of generalized AREs (recall \rf{A-Th}, \rf{Q-Th} and \rf{Hat})
\bel{ARE-Th-Th}\left\{\2n\ba{ll}
\ds PA_\Th+A_\Th^\top P+C_\Th^\top PC_\Th+Q_\Th-(PB+C_\Th^\top PD+S_\Th^\top)\Si^\dag(B^\top P+D^\top PC_\Th+S_\Th)=0,\\
\ns\ds\h P\h A_\BTh+\h A_\BTh^{\,\top}\1n\h P+\h C_\BTh^{\,\top}P\h C_\BTh+\h Q_\BTh-\big(\h P\h B+\h C_\BTh^{\,\top}\1n P\h D
+\h S_\BTh^{\,\top}\big)\bar\Si^\dag\big(\h B^\top\h P+\h D^\top P\h C_\BTh+\h S_\BTh\big)=0,\\
\ns\ds\Si\equiv R+D^\top P D\ges0,\qq\sR\big(B^\top P+D^\top PC_\Th+S_\Th\big)\subseteq\sR(\Si),\\
\ns\ds\bar\Si\equiv\h R+\h D^\top P\h D\ges0,\qq\sR\big(\h B^\top\h P+\h D^\top P\h C_\BTh+\h S_\BTh\,\big)\subseteq\sR\big(\bar\Si\big),\ea\right.\ee
admits a unique static stabilizing solution pair $(P,\h P)\in\dbS^n\times\dbS^n$, and the BSDE on $[0,\i)$:
\bel{BSDE-Th-Th}\ba{ll}
\ns\ds -d\eta(t)=\big\{A_\Th^\top\eta(t)-\big(PB+C_\Th^\top PD+S_\Th^\top\big)\Si^\dag\big[B^\top\eta(t)+D^\top\z(t)+D^\top P\si(t)+\rho(t)\big]\\
\ns\ds\qq\qq\q +C_\Th^\top\big[\z(t)+P\si(t)\big]+Pb(t)+\wt q(t)\big\}dt-\z(t)dW(t),\qq t\ges0,\ea\ee
admits a solution $(\eta(\cd),\z(\cd))\in\sX[0,\i)\times L^2_\dbF(\dbR^n)$ such that
\bel{range-eta-Th}B^\top\big(\eta(\cd)-\dbE[\eta(\cd)]\big)+D^\top\big(\z(\cd)-\dbE[\z(\cd)]\big)+D^\top P\big(\si(\cd)-\dbE[\si(\cd)]\big)+\rho(\cd)-\dbE[\rho(\cd)]\in\sR(\Si),\q \as,\ee
and the ODE on $[0,\i)$:
\bel{BODE-Th-Th}\ba{ll}
\ns\ds \dot{\bar\eta}(t)+\h A_\BTh^\top\bar\eta(t)-\big(\h P\h B+\h C_\BTh^{\,\top}P\h D+\h S_\BTh^\top\big)\bar\Si^\dag\big\{\h B^\top
\bar\eta(t)+\h D^\top\big(P\dbE[\si(t)]+\dbE[\z(t)]\big)+\dbE[\rho(t)]\big\}\\
\ns\ds\qq +\h C_\BTh^{\,\top}\big(P\dbE[\si(t)]+\dbE[\z(t)]\big)+\dbE[\wt q(t)]+\h P\dbE[b(t)]=0,\qq t\ges0,\ea\ee
admits a solution $\bar\eta(\cd)\in L^2(\dbR^n)$ such that
\bel{range-bar-eta-Th}\h B^\top\bar\eta(\cd)+\h D^\top P\dbE[\si(\cd)]+\h D^\top\dbE[\z(\cd)]+\dbE[\rho(\cd)]\in\sR(\bar\Si).
\ee
By \rf{range-eta-Th} and noting that
$$\big(PB+C_\Th^\top PD+S_\Th^\top\big)\Si^\dag=\big(PB+C^\top PD+S^\top\big)\Si^\dag+\Th^\top\Si\Si^\dag,$$
it is clear that \rf{BSDE-Th-Th} is equivalent to \rf{BSDE-eta-control}. Similarly, by \rf{range-bar-eta-Th} and noting
$$\big(\h PB+\h C_\BTh^\top\h P\h D+\h S_\BTh^\top\big)\bar\Si^\dag=\big(\h PB+\h C^\top P\h D+\h S^\top\big)\bar\Si^\dag+\bar\Th^\top\bar\Si\bar\Si^\dag,$$
it is straightforward to show that \rf{BODE-Th-Th} is equivalent to \rf{BODE-bar-eta-control}.

The rest is to prove that $(P,\h P)$ is a static stabilizing solution pair to the system of generalized AREs \rf{ARE-LQ-control}. To this end, choose $(\L,\bar\L)\in\dbR^{m\times n}\times\dbR^{m\times n}$ such that
\begin{equation*}\left\{\2n\ba{ll}
\ds \Th^*\triangleq-\Si^\dag\big(B^\top P+D^\top PC_\Th+S_\Th\big)+(I-\Si^\dag\Si)\L,\\
\ns\ds \bar\Th^*\triangleq-\bar\Si^\dag\big(\h B^\top\h P+\h D^\top P\h C_\BTh+\h S_\BTh\big)+(I-\bar\Si^\dag\bar\Si)\bar\L,
\ea\right.\end{equation*}
is an MF-$L^2$-stabilizer of $[A_\Th,\bar A_\BTh,C_\Th,\bar C_\BTh;B,\bar B,D,\bar D]$. Since
\bel{Sigma-Sigma*}\ba{ll}
\ds \Si(\Th^*+\Th)=-\big(B^\top P+D^\top PC_\Th+S_\Th\big)+\Si\Th=-\big(B^\top P+D^\top PC+S\big),\\
\ns\ds \bar\Si(\bar\Th^*+\bar\Th)=-\big(\h B^\top\h P+\h D^\top P\h C_\BTh+\h S_\BTh\big)+\bar\Si\bar\Th
=-\big(\h B^\top\h P+\h D^\top P\h C+\h S\big),\ea\ee
we have
$$\ba{ll}
\sR(B^\top P+D^\top PC+S)\subseteq\sR(\Si),\q \sR\big(\h B^\top\h P+\h D^\top P\h C+\h S\big)\subseteq\sR\big(\bar\Si\big).\ea$$
Moreover, noting \rf{Sigma-Sigma*}, we have
$$\ba{ll}
\ds 0=\h PA_\Th+A_\Th^\top P+C_\Th^\top PC_\Th+Q_\Th-(PB+C_\Th^\top PD+S_\Th^\top)\Si^\dag(B^\top P+D^\top PC_\Th+S_\Th)\\
\ns\ds =\h PA+A^\top\h P+C^\top PC+Q+(PB+C^\top PD+S^\top)\Th+\Th^\top(B^\top P+D^\top PC+S)\\
 \q-(PB+C^\top PD+S^\top)\Si^\dag(B^\top P+D^\top PC+S)\\
 \q-(PB+C^\top PD+S^\top)\Si^\dag\Si\Th-\Th^\top\Si\Si^\dag(B^\top P+D^\top PC+S)\\
\ns\ds =PA+A^\top P+C^\top PC+Q-(PB+C^\top PD+S^\top)\Si^\dag(B^\top P+D^\top PC+S)\\
 \q+(PB+C^\top PD+S^\top)(I-\Si^\dag\Si)\Th+\Th^\top(I-\Si\Si^\dag)(B^\top P+D^\top PC+S)\\
\ns\ds =PA+A^\top P+C^\top PC+Q-(PB+C^\top PD+S^\top)\Si^\dag(B^\top P+D^\top PC+S)\\
 \q-(\Th^*+\Th)^\top\Si(I-\Si^\dag\Si)\Th-\Th^\top(I-\Si\Si^\dag)\Si(\Th^*+\Th)\\
\ns\ds =PA+A^\top P+C^\top PC+Q-(PB+C^\top PD+S^\top)\Si^\dag(B^\top P+D^\top PC+S),\ea$$
and similarly, we get
$$\ba{ll}
\ds 0=P\h A_\BTh+\h A_\BTh^{\,\top}\1n\h P+\h C_\BTh^{\,\top}P\h C_\BTh+\h Q_\BTh-\big(\h P\h B+\h C_\BTh^{\,\top}\1n P\h D
 +\h S_\BTh^{\,\top}\big)\bar\Si^\dag\big(\h B^\top\h P+\h D^\top P\h C_\BTh+\h S_\BTh\big)\\
%
%
%
%
\ns\ds =\h P\h A+\h A^\top\h P+\h C^\top P\h C+\h Q-(\h P\h B+\h C^\top P\h D+\h S^\top)\bar\Si^\dag(\h B^\top\h P+\h D^\top P\h C+\h S).\ea$$
Then we know that $(P,\h P)$ solves \rf{ARE-LQ-control}. By again \rf{Sigma-Sigma*}, we can find $(\L',\bar\L')\in\dbR^{m\times n}\times\dbR^{m\times n}$ such that
$$
\left\{
\ba{ll}
        \Th^*+\Th=-\Si^\dag\big(B^\top P+D^\top PC+S\big)+(I-\Si^\dag\Si)\L',\\
\bar\Th^*+\bar\Th=-\bar\Si^\dag\big(\h B^\top\h P+\h D^\top P\h C+\h S\big)+(I-\bar\Si^\dag\bar\Si)\bar\L'.
\ea\right.$$
Thus $(\Th^*+\Th,\bar\Th^*+\bar\Th)$ is an MF-$L^2$-stabilizer of $[A,\bar A,C,\bar C;B,\bar B,D,\bar D]$ and $(P,\h P)$ is static stabilizing. The rest of the proof is clear. The proof of Theorem \ref{T-2.9} is complete.
\endpf

\section{Concluding Remarks}

In this paper, we have presented a systematic theory for two-person non-zero sum differential games of mean-field SDEs with quadratic performance indexes in $[0,\i)$. The case of two-person zero-sum, which is also new, has been treated as a special case. Our results cover several existing ones in the literature for infinite horizon problems, including LQ optimal control problems of mean-field type (\cite{Huang-Li-Yong2015}), two-person zero-sum LQ stochastic differential games (without mean-field terms) (\cite{Sun-Yong-Zhang2016}), LQ optimal control problem (without mean-field terms), and the equivalence between the open-loop solvability and the closed-loop solvability for stochastic LQ problem in $[0,\i)$ (\cite{Sun-Yong2018}). Finally, we have to leave the following question open: Is the existence of the open-loop and closed-loop saddle points equivalent for the mean-field LQ two-person zero-sum stochastic differential game in an infinite horizon? We will research this topic in the future.

\section*{Acknowledgement}

The authors would like to thank Dr. Jingrui Sun of Department of Mathematics, Southern University of Science and Technology for quite a few discussions.

This work was carried out during the stay of Jingtao Shi at the University of Central Florida, from December 2019 to July 2020. He would like to thank the invitation of Professor Jiongmin Yong, the hospitality of the Department of Mathematics, University of Central Florida, and the financial support from the China Scholarship Council.

\end{document}